\documentclass[12pt]{amsart}

\usepackage{amssymb,latexsym}

\usepackage[ansinew]{inputenc}
\usepackage{graphicx}
\usepackage{stmaryrd}

\def\cal{\mathcal}

\def\Bbb{\mathbb}

\def\C{{\Bbb C}}

\def\R{{\Bbb R}}

\def\<{\left<}
\def\>{\right>}

\def\PSet{\mbox{\rm I\kern-.22em P}}

\def\S{{\cal S}}

\def\BMO{{\rm BMO}}

\def\({( \hspace{-0.335em}(}
\def\){) \hspace{-0.335em})}
\def\supp{{\rm supp\,}}

\def\fu2{\frac{n}{2}}

\def\l({\left(}
\def\r){\right)}
\def\be{\begin{enumerate}}
\def\ee{\end{enumerate}}

\def\sign{{\rm sign}}

\def \dint {\int \! \! \! \int }

\def \diam {{\rm diam}}

\newtheorem{theorem}{Theorem}[section]
\newtheorem{definition}[theorem]{Definition}
\newtheorem{lemma}[theorem]{Lemma}
\newtheorem{corollary}[theorem]{Corollary}
\newtheorem{proposition}[theorem]{Proposition}

\newtheorem{thm}{Theorem}[section]

\newtheorem{remark}[thm]{Remark}

\address{Centre for Mathematical Sciences, University of Lund, Lund, Sweden}\email{s.pott@maths.lth.se}
\address{Department of Mathematics, University of Glasgow, Glasgow, UK}\email{paco.villarroya@maths.gla.uk}

\author{Sandra Pott}
\author{Paco Villarroya}
\date{}
\title{A T(1) theorem on product spaces}
\thanks{The first author was supported by a Heisenberg Fellowship of the German Research Foundation (DFG). The second author was 
supported by a Marie-Curie Intra European Fellowship PIEF-GA-2008-220266 and Spanish MEC project
MTM2008-04594.}

\begin{document}



\date{\today}

\begin{abstract}

We prove a new T(1) theorem for multiparameter singular integral operators.

\end{abstract}

\maketitle

\section{Introduction}

\subsection{Historical introduction.}

In 1984 G. David and J.L. Journ\'e (see \cite{DJ}) published their celebrated $T(1)$ theorem, a result that characterizes the $L^2$-boundedness of non-convolution
integral operators with a Calder\'on-Zygmund kernel. In their theorem, the necessary and sufficient conditions for boundedness
are expressed by the behaviour of the operator when acting
over particular families of functions: the membership in $\BMO$ of properly defined $T(1)$, $T^{*}(1)$ functions and
the so-called weak boundedness property, which is the validity of $L^2$ bounds when duality is tested over bump functions with the same space localization.

Since then, many other proofs of this fundamental result in the theory of singular integration have appeared, while it has also been
extended to a large variety of settings. Actually, only one year later Journ\'e \cite{J1} established the extension to product spaces
when he proved an analogous result of $L^2$-boundedness for
multiparameter singular integrals. These are
operators whose class of kernels is homogeneous
with respect to non-isotropic dilations of the form $\rho_{\delta_1,\ldots ,\delta_n}(x_1,\ldots, x_n)=(\delta_1x_1,\ldots, \delta_nx_n)$
for $x_i\in \mathbb R^{d_i}$ and $\delta_i>0$, where the number of parameters of the problem coincides with the quantity of independent dilations.
The simplest examples of such operators are convolution type operators like the multiple Hilbert transform defined in $\mathbb R^n$ by
$$
H_1\cdots H_n(f)={\rm p.v.} f*\frac{1}{x_1\ldots x_n}
$$
or the multiple Riesz transforms defined in $\prod_{i=1}^{n}\mathbb R^{d_i}$ by
$$
R_{j_1}\cdots R_{j_n}(f)=f*(\frac{\pi_{j_1}(x_1)}{|x_1|^{d_1+1}}\cdots \frac{\pi_{j_n}(x_n)}{|x_n|^{d_n+1}})
$$
where $\pi_{j_{i}}$ is the orthogonal projection from $\mathbb R^{d_{i}}$ into $\mathbb R$ that \textquotedblleft keeps\textquotedblright \, the $j_{i}$-coordinate.


A direct application of Fubini's theorem shows that the multiple Hilbert transform is bounded in all spaces $L^p(\mathbb R^{n})$ for $1<p<\infty $. However,
the situation
is not so simple for more general multiparameter singular integrals, especially if they are of non-convolution type.
These multiparameter operators, even in the simplest cases, are very different from their classical
counterparts mainly because the singularities of their kernels lie not only at the origin as in the case of standard Calder\'on-Zygmund kernels, but instead,
they are spread over larger subspaces. For example, in the case of multiple Hilbert transform the set of singularities is the union of the coordinate axes $x_i=0$.
As a consequence, these operators are not in general weak type on $L^{1}(\mathbb R^n)$ and moreover, the
strong maximal operator does not control their boundedness properties.

The
main motivation to extend the theory of singular integration to operators that commute with multiparameter families of dilations
comes from their close relationship with multiplier operators in $\mathbb R^n$. Namely,
in the same way as the classical linear Hilbert transform is
closely related to the Fourier partial sum operator $S_{N}(f)(x)=\sum_{|k|<N}\hat{f}(k)e^{2\pi ik x}$, different multiparameter singular integrals are related
to different Fourier partial sum operators in several variables. In particular, the rectangular partial sums operator defined in $\mathbb R^n$ by
$$
S_{N_1,\ldots ,N_n}(f)(x_1,\ldots ,x_n)=\sum_{j=1}^{n}\sum_{|k_j|<N_j}\widehat{f}(k_1,\ldots ,k_n)e^{2\pi ik_{j} x_{j}}
$$
is controlled by the multiple Hilbert transform. In section 3 we apply our main result to extend boundedness of product multiplier operators to the
non-convolution setting.

We notice that although multiparametric singular integrals were intensively studied more 
than two decades ago, due to recent developments in Harmonic Analysis this issue  
has experienced in the last years a renewed interest as it can be seen from the papers 
\cite{B}, \cite{DMT}, \cite{MPTT1}, \cite{MPTT2}, \cite{LM} and \cite{PT} among others.

\subsection{On Journ\'e's theorem}

Journ\'e's result is the first attempt to characterize $L^2$ boundedness of non-convolution multiparameter singular integral operators.
As stated before, many of the classical techniques, like for example a proper Calder\'on-Zygmund decomposition
and the control of singular integrals by means of
maximal functions, are no longer available in the multiparameter setting. So, the method Journ\'e chose to overcome this difficulty
was the use of vector valued Calder\'on-Zygmund theory.
In order to state his theorem in a simplified form, we require some notation.

Let $\Delta $ be the diagonal in $\mathbb R^2$ and $B$ be a Banach space.
A continuous function $K:\R^2 \setminus \Delta \to B$
is called a vector valued standard Calder\'on-Zygmund kernel, if for some $0<\delta\leq 1$ and some
constant $C>0$, we have
\begin{equation*}\label{e.e1Cal1}
\| K(x,t)\|_B\leq C |x-t|^{-1}
\end{equation*}
\begin{equation*}\label{e.e1Cal2}
\| K(x,t)-K(x',t')\|_B\leq C (|x-x'|+|t-t'|)^\delta |x-t|^{-1-\delta}
\end{equation*}
whenever $|x-x'|+|t-t'|\leq |x-t|/2$. In this context, $|K|$ usually denotes the best constant in both inequalities.

\begin{definition}\label{journesi}
A continuous linear mapping $T$ from $C^{\infty }_0(\mathbb R)\otimes C^{\infty }_0(\mathbb R)$ into
its algebraic dual is called
a singular integral operator, if
there are vector-valued C-Z kernels
$K^1,K^2:\mathbb R^2 \setminus \Delta \rightarrow {\cal L}(L^2(\mathbb R),L^2(\mathbb R))$  such that for
$f_1,f_2,g_1,g_2\in C^{\infty }_0(\mathbb R)$, we have
$$
\langle T(f_1\otimes f_2),g_1\otimes g_2\rangle =\int_{\mathbb R^2}\int_{\mathbb R^2}f_1(t_1)g_1(x_1)\langle K^1(x_1,t_1)f_2,g_2\rangle dt_1dx_1
$$
whenever $\supp f_1\cap \supp g_1=\emptyset $, and symmetrically for $K^2$.
\end{definition}

The definition of the weak boundedness property makes use of the restricted operators:
given $f_i,g_i\in C^{\infty }_0(\mathbb R)$ for $i=1,2$, let $\langle T^{i}(f_i),g_i\rangle : C^{\infty }_0(\mathbb R)\rightarrow
C^{\infty }_0(\mathbb R)'$ be defined by
$$
\langle \langle T^{1}(f_2),g_2\rangle f_1, g_1\rangle =\langle \langle T^{2}(f_1),g_1\rangle f_2, g_2\rangle
=\langle T(f_1\otimes f_2),g_1\otimes g_2\rangle
$$
Notice that the kernel of $T^1$ for
example is precisely $\langle K^1(x_1,t_1)(f_2),g_2\rangle $.


Then, a singular integral operator $T$ is said to satisfy the weak boundedness property if for any bounded subset $A$ of $C^{\infty }_0(\mathbb R)$
there is a constant $C>0$, that may depend on $A$, such that for any $f,g\in A$ we have that
$$
\| \langle T^{i}(f_{x,R}),g_{x,R}\rangle \|_{CZ}:=\| \langle T^{i}(f_{x,R}),g_{x,R}\rangle \|_{L^2(\mathbb R)\rightarrow L^2(\mathbb R)}+|K^{i}|
\leq C
$$
where $f_{x,R}(y)=R^{-1/2}f(R^{-1}(y-x))$, and the same for $g_{x,R}$.

Finally, also associated with $T$ we can define its partial adjoints as the adjoint operators with respect to each variable, that is, the
operator given by
$$
\langle T_1(f_1\otimes f_2),g_1\otimes g_2\rangle =\langle T(g_1\otimes f_2),f_1\otimes g_2\rangle
$$
and analogously for $T_2$.
Notice that $T_2=T_1^{*}$.

With all these definitions we can state Journ\'e's result (see \cite{J1}):
\begin{theorem}
Let $T$ be a singular integral operator on $\mathbb R\times \mathbb R$ as described in definition \ref{journesi} satisfying the weak boundedness property and
$T(1),T^{*}(1),T_1(1),T_1^{*}(1)\in \BMO_{\rm prod}(\mathbb R^2)$. Then $T$ extends boundedly on $L^2(\mathbb R^2)$.
\end{theorem}

We would like to stress here how restrictive these conditions are, in particular the
definitions of singular integral operator and of the weak boundedness property. When written in the language of
vector valued Calder\'on-Zygmund theory they look quite simple, but a more detailed description reveals all their complexity.
The sufficient hypotheses for $T$ to be bounded on
$L^{2}(\mathbb R^2)$ are the following ones:
\begin{itemize}
\item[a)] The $K^{i}$ are vector valued C-Z kernels.
This condition implies that
$K^1(x_1,t_1)$ are C-Z operators bounded on $L^2(\mathbb R)$
and that, moreover, their C-Z norms defined by
$\| K^1(x_1,t_1)\|_{CZ}:=\| K^{1}(x_1,t_1)\|_{L^2(\mathbb R)\rightarrow L^2(\mathbb R)}+|K^{1}_{x_1,t_1}|$
satisfy
$$
\| K^{1}(x_1,t_1)\|_{CZ} \leq C |x_i-t_i|^{-1}
$$
$$
\| K^{1}(x_1,t_1)-K^{1}(x_{1}',t_{1}')\|_{CZ}\leq C (|x_1-x_{1}'|+|t_1-t_{1}'|)^\delta |x_1-t_1|^{-1-\delta}
$$
whenever $|x_1-x_1'|+|t_1-t_1'|\leq |x_1-t_1|/2$, and the same for $K^{2}(x_2,t_2)$.

\item [b)] Weak boundedness property.
This condition implies that $\langle T^{1}(f_{s,R}),g_{s,R}\rangle $
are also C-Z operators bounded on
$L^{2}(\mathbb R)$ and moreover, their C-Z norms defined as
$
\| \langle T^{1}(f_{s,R}),g_{s,R}\rangle \|_{CZ}:=\| \langle T^{1}(f_{s,R}),g_{s,R}\rangle \|_{L^2(\mathbb R)\rightarrow L^2(\mathbb R)}+|K^{1}|
$
satisfy
$$
\| \langle T^{1}(f_{s,R}),g_{s,R}\rangle \|_{CZ}\leq C
$$
and the same for $T^{2}$.

\item[c)] $T(1),T^{*}(1),\tilde{T}_1(1),\tilde{T}_1^{*}(1)\in \BMO_{\rm prod}(\mathbb R^2)$, the latter space being much more complex that its one
variable counterpart.

\end{itemize}

So, in order to conclude that the product operator is bounded, Journ\'e's theorem assumes that \textquotedblleft some parts\textquotedblright \, of the operator, in particular the vector valued
kernels and the restricted operators, are known to be bounded a priori. This is quite a different situation from the original $T(1)$ theorem in which nothing
is assumed to be bounded a priori.
However, for the same above mentioned reasons,
the use of vector valued theory was also adopted by other authors in later developments of singular integration in product spaces
(see \cite{F1}, \cite{F2}, \cite{F3}, \cite{J1} and \cite{J2}).

Our purpose in the present paper is to state and prove a new $T(1)$ theorem for product spaces in which any hypothesis related to operators which
need to be bounded a priori disappear. Therefore, we avoid the use of vector valued Calder\'on-Zygmund theory. Instead,
we seek other sufficient hypotheses for $L^2$-boundedness which are
much closer to the spirit of the classical $T(1)$ theorem of David and Journ\'e: conditions related to scalar
decay estimates of the kernel and to the behaviour of the operator over special
families of functions. To get such new hypotheses, we combine the three classical conditions (kernel estimates, weak boundedness condition and
$T(1)\in \BMO$) according to their separate action over different parameters
to generate a range of new mixed conditions. For example, in the bi-parameter case we consider new properties by combining two classical ones,
namely kernel decay estimates in one parameter and weak boundedness property in the other parameter, to get what we call the mixed \textquotedblleft kernel\textquotedblright-\textquotedblleft weak boundedness\textquotedblright \,
condition. As a result, we obtain nine different conditions that cover all possible combinations. This procedure better preserves the symmetry given by
the product structure of the kernels and therefore, it is better suited for the general multi-parameter situation.

Moreover, in lemma \ref{representation} we obtain a decomposition of the operators under consideration which
shows that the quantity and the statement of our conditions are the right ones in the sense that they entirely describe the boundedness properties of the terms of the decomposition.
We would like to highlight the role played by some of those conditions that give sense to a new class of paraproducts,
which do not appear in previous developments of the theory. We plan a deeper study of such operators in forthcoming papers.


The main advantage of our approach is that, as we said before, Journ\'e's technique needs to assume that the operator when restricted to the one-parameter case is known to be bounded a-priori. However, in our 
approach, we only require the restricted operators to be weakly bounded. This might look like a very small gain, specially because, since we also demand the restricted operators to satisfy the cancellation conditions, the classical one-parameter theory proves our restricted operators to be also bounded. Actually the following bound is known to hold
$$
\| T^{i}_{x_{i},t_{i}}\|_{2\rightarrow 2}\leq C_{x_{i},t_{i}}+ \| T^{i}_{x_{i},t_{i}}(1)\|_{\BMO}+\| (T^{i}_{x_{i},t_{i}})^{*}(1)\|_{\BMO}
$$
with $C_{x_{i},t_{i}}>0$ the constant in the weak boundedness condition satisfied by $T^{i}_{x_{i},t_{i}}$.
And yet, notice that unlike Journ\'e's work which essentially assumes the whole bound 
$\| T^{i}_{x_{i},t_{i}}\|_{2\rightarrow 2}\leq C|x_{i}-t_{i}|^{-1}$, we only require
$C_{x_{i},t_{i}}\leq C|x_{i}-t_{i}|^{-1}$.

Another advantage is that, 
at least in principle, the result
can be applied to a larger family of operators since in our hypotheses no operator is ever assumed to be bounded a priori.
Actually, none of the examples treated in section 3 are under the scope of Journ\'e's Theorem.
Moreover, those new conditions should, again in principle, be
easier to test since there is no need to calculate operators norms.


On the other hand, the price we have to pay for adopting this new point of view is a larger number of hypotheses, growing rapidly with the number of parameters. In the
case of $n$-parameter operators we have to deal with $3^n$ hypotheses to ensure that the operator is bounded.
However, although Journ\'e's Theorem only
requires three conditions and so its statement remains as concise as
in the uni-parameter case, when the number of parameters grows, these three hypotheses need to be applied iteratively. Then, one might also
say that the number of conditions also increases exponentially.
From this perspective, the vector valued formulation turns out to be a clever way to encode the complicated structure of the problem, and
when one unfolds all the information, the complexity always grows accordingly.

Finally, it has to be said that either Journ\'e's Theorem and our result exhibit a common weak point: the given sufficient conditions for
boundedness of product singular integrals are not necessary.
This was first shown by Journ\'e (see the same paper \cite{J1}) when he constructed an example of a bounded operator for which the partial adjoint
$T_{1}(1)$ is not in $\BMO_{\rm prod}(\mathbb R^2)$.
The problem is that either in his theorem and in ours, the stated conditions imply not only boundedness of $T$ but also of $T_{1}$
(and so in such case $T_{1}(1)$ will have
to be in $\BMO_{\rm prod}(\mathbb R^2)$). The underlying reason for this is that the partial adjoint of a bounded operator on $L^2(\mathbb R^2)$ is not necessarily bounded. In the
language of operator spaces, taking adjoints is not a completely bounded map. In fact, the boundedness properties of the so-called mixed multiparameter paraproducts, which is a necessary part of any full characterization of boundedness of such product singular integrals operators, remain to be completely understood. This is the subject of ongoing research. 

The paper is organized as follows. In section  \ref{main} we define all the sufficient hypotheses for $L^2$-boundedness of biparameter singular integral operators and
state our $T(1)$ theorem. We also state without proof the analogous results for multiparameter operators in several variables.
In section \ref{apply} we apply our main result to prove boundedness of non-convolution operators previously studied by R. Fefferman and E. Stein in the
convolution setting.

We start the proof of our result in section \ref{T(1)} by the rigorous definition of the functions $T(1)$ and $T(\phi_{I}\otimes 1)$.
In section \ref{bump} we obtain an appropriate estimate for the rate of decay of the action
of the operator over bump functions when special cancellation properties are assumed. Sections \ref{L2} and \ref{Lp} focus on
the proof of $L^2$ boundedness and the extension to $L^{p}$ spaces respectively, both of them under the special cancellation hypotheses.
The latter case makes use of new bi-parameter modified square functions whose boundedness is a direct consequence of
analogous uni-parameter modified square functions. The proof of boundedness of these
new square functions is provided in an appendix at the end of the paper.
Finally, in section \ref{para} we construct the necessary paraproducts to prove the result
in the general case, that is, in absence of the cancellation assumptions.

In a sequel to the present paper, we plan to deal the endpoint case of boundedness from  
$L^\infty(\mathbb R^2)$ into $\BMO_{\rm prod}(\mathbb R^2)$.

We would like to thank Anthony Carbery and Jim Wright for valuable conversations and helpful comments. We would also like to acknowledge
the School of Mathematics of the
University of Edinburgh for the stimulating research environment
provided
which so positively influenced the development of this work.


\section{Definitions and statement of the main theorem}\label{main}

\begin{definition}\label{prodCZ}
Let $\Delta $ be the diagonal in $\mathbb R^2$.
A function $K:(\R^2 \setminus \Delta )\times (\R^2 \setminus \Delta )\to \mathbb R$ is called a
product Calder\'on-Zygmund kernel, if for some $0<\delta\le 1$ and some
constant $C>0$ we have
$$
\begin{array}{l}
|K(x,t)|\le C {\displaystyle \prod_{i=1,2}\frac{1}{|x_i-t_i|}} \\
|K(x,t)-K((x_1,x_2'),(t_1,t_2'))-K((x_1',x_2),(t_1',t_2))+K(x',t')|
\le C \displaystyle{\prod_{i=1,2}\frac{(|x_i-x_i'|+|t_i-t_i'|)^\delta}{|x_i-t_i|^{1+\delta}}}
\end{array}
\vspace{-.5cm}
$$
\hspace{.5cm} whenever $2(|x_i-x_i'|+|t_i-t_i'|)\leq |x_i-t_i|$.
\end{definition}

\vskip10pt
\begin{remark}
For $\delta =1$, the second hypothesis is satisfied if we assume the stronger smoothness condition
$$
|\partial_{t_1}\partial_{t_2}K(x,t)|
+|\partial_{t_1}\partial_{x_2}K(x,t)|+|\partial_{x_1}\partial_{t_2}K(x,t)|
+|\partial_{x_1}\partial_{x_2}K(x,t)|\leq C \prod_{i=1,2}|x_i-t_i|^{-2}
$$
(notice that the derivatives $\partial_{t_1}\partial_{x_1}K(x,t)$ and $\partial_{t_2}\partial_{x_2}K(x,t)$ do not appear 
in this condition).
This is due to the trivial inequality
\begin{eqnarray*}
&&|K(x,t)-K((x_1,x_2'),(t_1,t_2'))-K((x_1',x_2),(t_1',t_2))+K(x',t')|\\
&\leq &|K(x,t)-K((x_1,x_2),(t_1,t_2'))-K((x_1,x_2),(t_1',t_2))+K(x,t')|\\
&+&|K((x_1,x_2),(t_1,t_2'))-K((x_1,x_2'),(t_1,t_2'))-K((x_1,x_2),(t_1',t_2'))+K((x_1,x_2'),(t_1',t_2'))|\\
&+&|K((x_1,x_2),(t_1',t_2))-K((x_1,x_2),(t_1',t_2'))-K((x_1',x_2),(t_1',t_2))+K((x_1',x_2),(t_1',t_2'))|\\
&+&|K(x,t')-K((x_1,x_2'),(t_1',t_2'))-K((x_1',x_2),(t_1',t_2'))+K(x',t')|
\end{eqnarray*}

\end{remark}

\begin{definition}\label{intrep}
A bilinear form
$\Lambda:\S(\R^2)\times \S(\R^2)\to \C$
is said to be associated with a product Calder\'on-Zygmund kernel $K$ if it satisfies the following 
three integral representations:
\begin{enumerate}
\item for all Schwartz functions $f,g\in {\cal S(\R^2)}$
such that $f(\cdot ,t_2), g(\cdot ,x_2)$ and $f(t_1,\cdot ), g(x_1,\cdot )$ have respectively disjoint supports
$$
\Lambda(f,g)=\int_{\R^2}\int_{\R^2} f(t)g(x) K(x,t)\, dt \, dx
$$

\item for all Schwartz functions 
$f_1,f_2,g_1,g_2\in {\cal S(\R )}$
such that $f_1$ and $g_1$ have disjoint supports, we have
there is a bilinear form $\Lambda_{t_1,x_1}^1:\S(\R)\times \S(\R)\to \C$ 
$$
\Lambda(f,g)=\int_{\R^2}\int_{\R^2} f_1(t_1)g_1(x_1) \Lambda_{t_1,x_1}^1(f_2,g_2)\, dt_1 \, dx_1
$$

\item[another (2)] for all Schwartz functions 
$f,g\in {\cal S(\R^2)}$
such that $f(\cdot ,t_{2}), g(\cdot ,x_{2})$ have disjoint supports for all $t_{2},x_{2}$, 
there is a bilinear form $\Lambda_{t_1,x_1}^1:\S(\R)\times \S(\R)\to \C$ 
$$
\Lambda(f,g)=\int_{\R}\int_{\R} \Lambda_{t_1,x_1}^1(f_1(t_{1},\cdot ),g_1(x_{1},\cdot ))\, dt_1 \, dx_1
$$

\item analogous representation with $\Lambda_{t_2,x_2}^2$.
\end{enumerate}

Here, we have defined the restricted bilinear forms $\Lambda^1$, $\Lambda^2$ by
$$
\langle \Lambda^1(f_2,g_2)f_1,g_1\rangle =\langle \Lambda^2(f_1,g_1)f_2,g_2\rangle 
=\Lambda(f_1\otimes f_2,g_1\otimes g_2),
$$
and  $\Lambda_{t_1,x_1}^1(f_2,g_2)$ is a kernel associated to $\Lambda^1(f_2,g_2)$ in the usual sense in one parameter. (NO NEED)

If the form $\Lambda$ is continuous on ${\cal S}(\mathbb R^2) \times {\cal S}(\mathbb R^2)$, then it will be called a bilinear Calder\'on-Zygmund form.
\end{definition}

With a small abuse of notation, we will say that a bilinear form is bounded on $L^p(\mathbb R^2)$ if there is a constant $C>0$ such that
$|\Lambda (f,g)|\leq C\| f\|_{L^{p}(\mathbb R^2)}\| g\|_{L^{p'}(\mathbb R^2)}$ for all $f,g\in {\cal S}(\mathbb R^2)$.

\vskip15pt
Given a bilinear form $\Lambda:\S(\R^2)\times \S(\R^2)\to \C$,
we define associated linear operators $T$, adjoint bilinear forms $\Lambda_i$, and restricted linear forms $\Lambda^{i}$, in the following way:

\begin{definition}{(Dual operators).}
Given a bilinear form $\Lambda $, we define linear operators $T$, $T^*$ through duality:
$$
\langle T(f),g\rangle =\langle f,T^*(g)\rangle =\Lambda(f,g)
$$
\end{definition}

\begin{definition}{(Adjoint bilinear forms).} We define the adjoint bilinear forms $\Lambda_{i}$ such that for $f=f_1\otimes f_2$, $g=g_1\otimes g_2$
functions of tensor product type, we have
$$
\Lambda_1(f,g)=\Lambda (g_1\otimes f_2,f_1\otimes g_2),
\hskip 10pt
\Lambda_2(f,g)=\Lambda (f_1\otimes g_2,g_1\otimes f_2)
$$
and then extended by linearity and continuity.
\end{definition}
These new bilinear forms are also associated with linear operators  $T_1$, $T_2$ via duality
$$
\langle T_i(f),g\rangle =\langle f,T_{i}^{*}(g)\rangle =\Lambda_i(f,g)
$$
which in the case of tensor products, $f=f_1\otimes f_2$, $g=g_1\otimes g_2$, satisfy
$$
\langle T_1(f_1\otimes f_2),g_1\otimes g_2\rangle =\Lambda_1(f,g)
=\Lambda(g_1\otimes f_2,f_1\otimes g_2)=\langle T(g_1\otimes f_2),f_1\otimes g_2\rangle
$$
Notice that $T_{2}=T_{1}^{*}$ and $T_{2}^{*}=T_{1}$. 

From now we will sometimes write $\Lambda_0=\Lambda $ and $T_{0}$ associated to $\Lambda_{0}$.

\begin{definition}{(Restricted bilinear forms).}
As above, we define the restricted bilinear forms by
$$
\langle \Lambda^1(f_2,g_2)f_1,g_1\rangle =\langle \Lambda^2(f_1,g_1)f_2,g_2\rangle 
=\Lambda(f_1\otimes f_2,g_1\otimes g_2)
$$
We will denote by $T^{i}$ the linear operators associated with the restricted bilinear form $\Lambda^{i}$ through duality, 
$\Lambda^{i}(f_j,g_j)=\langle T^{i}(f_j),g_j\rangle $.
\end{definition}


We also notice that we use subindexes to denote the partial adjoint operators or forms while we use superindexes to denote the restricted ones.

As mentioned in the introduction, boundedness of the bilinear form $\Lambda $ implies boundedness of the dual linear operators $T$ and $T^{*}$,
but it does not imply boundedness of any of
the adjoint bilinear forms $\Lambda_i$ nor their corresponding associated adjoint operators $T_{i}$, $T_{i}^{*}$. In other words,
boundedness of $\Lambda $ on $L^2(\mathbb R^2)$
$$
\Big| \Lambda (\sum_n f_1^{n_1}\otimes f_2^{n_2},\sum_m g_1^{m_1}\otimes g_2^{m_2})\Big|
\leq C\Big\|\sum_n f_1^{n_1}\otimes f_2^{n_2}\Big\|_2 \big\| \sum_m g_1^{m_1}\otimes g_2^{m_2}\Big\|_2
$$
implies boundedness of $\Lambda_1$ only on $L^2(\mathbb R)\hat{\otimes }L^2(\mathbb R)$:
$$
\Big|\Lambda_1 (\sum_n f_1^{n_1}\otimes f_2^{n_2},\sum_m g_1^{m_1}\otimes g_2^{m_2})\Big|
=\Big|\sum_{n,m} \Lambda (g_1^{m_1}\otimes f_2^{n_2},f_1^{n_1}\otimes g_2^{m_2})\Big|
$$
$$
\leq C\sum_{n,m}\| g_1^{m_1}\|_2 \| f_2^{n_2}\|_2  \| f_1^{n_1}\|_2 \| g_2^{m_2}\|_2
=C\sum_{n}\| f_1^{n_1}\otimes  f_2^{n_2}\|_2 \sum_{m} \| g_1^{m_1}\otimes  g_2^{m_2}\|_2
$$

\begin{definition}
For every interval $I\subset \mathbb R$ we denote its centre by $c(I)$ and its length by $|I|$.
Then, a $L^2(\mathbb R)$-normalized bump function adapted to $I$ with constant $C>0$ and order $N\in \mathbb N$,
is a Schwartz function $\phi $ such that
$$|\phi^{(n)}(x)|\le C|I|^{-1/2-n}(1+|I|^{-1}|x-c(I)|)^{-N}, \ \ \ 0\leq  n \le N$$

We will denote the cube in $\mathbb R^d$ of measure one centered at the origin as the standard cube.
A bump function $\phi $ adapted to the standard cube of order $N$
is a Schwartz function satisfying
$$|\partial^\alpha \phi (x)|\le C(1+|x|)^{-N}, \ \ \ 0\leq |\alpha|\le N$$

A bump function $\phi $ is said to be adapted to a box $R$ in $\R^d$, if
for any affine linear transformation $A$ transforming the
standard cube into the box $R$, the function
$$|\det(A)|^{1/2}\phi(Ax)$$
is a bump function adapted to the standard cube. This definition
does not depend on the choice of the linear transformation.

\end{definition}

Observe that all these bump functions are normalized to be uniformly bounded in $L^2(\mathbb R^d)$. The order of the bump functions
will always be denoted by $N$, even though its value might change from line to line.
It is also worth saying that we usually reserve the greek letter $\phi $, $\varphi $ for general bump functions while
we reserve the use of $\psi $ to denote bump functions with mean zero.

\begin{definition}\label{WB}
We say that a bilinear form $\Lambda $ satisfies the weak boundedness condition, if for any rectangle $R$ and every pair $\phi_R, \varphi_R$ of
$L^2$-normalized bump functions adapted to $R$ with constant $C$, we have
$$
|\Lambda (\phi_R,\varphi_R)|\leq C
$$
\end{definition}

\begin{definition}\label{WB-CZ}
We say that a bilinear form $\Lambda $ satisfies the mixed weak boundedness-Calder\'on Zygmund condition, if 
whenever $2(|x_i-x'_i|+|t_i-t'_i|)<|t_i-x_i|$ for $i\in \{ 1,2\}$, we have 
$$
|\Lambda^i_{t_i,x_i}(\phi_I,\varphi_I)|\leq C|t_i-y_i|^{-1}
$$
$$
|(\Lambda^i_{t_i,x_i}-\Lambda^i_{t'_i,x'_i}) (\phi_I,\varphi_I)|\leq C(|x_i-x'_i|+|t_i-t'_i|)^{\delta }|t_i-x_i|^{-(1+\delta )}
$$
for any interval $I$ and
every pair $\phi_I, \varphi_I$ of $L^2$-normalized bump functions adapted to $I$ with constant $C$, 
\end{definition}

Obviously, the second condition is implied by the smoothness condition
$$
|\partial_{t_i}\Lambda^i_{t_i,x_i}(\phi_I,\varphi_I)|+|\partial_{x_i}\Lambda^i_{t_i,x_i}(\phi_I,\varphi_I)|
\leq C|t_i-x_i|^{-(1+\delta )}
$$

\vskip10pt
Finally, we notice that, in order to simplify notation, from now being the space product $\BMO$, that is, the dual of 
$H^1(\mathbb R^2)$, will be simply denoted by $\BMO(\mathbb R^2)$. 

We can now state our main result,

\begin{theorem}{(bi-parameter T(1) theorem).}\label{biparT1}
Let $\Lambda$ be a bilinear Calder\'on-Zygmund form satisfying the mixed WB-CZ conditions.
Then, the following are equivalent:
\begin{enumerate}
\item $\Lambda_i$ are bounded bilinear forms on $L^2(\mathbb R^2)$ for all $i=0,1,2$,

\item
$\Lambda $ satisfies
the weak boundedness condition
and the special cancellation conditions:
\begin{enumerate}
\item[a)]
$
T(1),T^{*}(1),T_1(1),T_1^{*}(1)\in \BMO(\mathbb R^2)
$,
\item[b)]
$
\langle T(\phi_{I}\otimes 1),\varphi_{I}\otimes \cdot \rangle ,\langle T(1\otimes \phi_{I}),\cdot \otimes \varphi_{I}\rangle ,
\langle T^{*}(\phi_{I}\otimes 1),\varphi_{I}\otimes \cdot \rangle ,
\langle T^{*}(1\otimes \phi_{I}),\cdot \otimes \varphi_{I}\rangle \in \BMO (\mathbb R)
$
for all $\phi_{I}$, $\varphi_{I}$ bump functions adapted to $I$
with norms uniformly bounded in $I$.
\end{enumerate}
\end{enumerate}
\end{theorem}

We remark that
boundedness of those operators $T_i$ and $T_j$ for $i\neq j$ are not equivalent.
A way to show this is by considering Carleson's example that proves $\BMO_{\rm rec}(\mathbb R^2)\subsetneq \BMO(\mathbb R^2)$.
In \cite{C}, he described a recursive process to construct a sequence of functions such that
$\| b_{n}\|_{{\BMO}_{\rm rec}(\mathbb R^2)}=1$ while $\| b_{n}\|_{\BMO_{\rm prod}(\mathbb R^2)}\geq C_n$ where $(C_n)_{n\in \mathbb N}$ was an unbounded sequence of
positive numbers. Then, we can consider paraproducts associated with that sequence of functions
$$
T_n(f)
=\sum_{I}\langle b_{n},\psi_R\rangle \langle f,\psi_R^2\rangle \psi_R
$$
in such a way that $\| T_n\|_{L^{2}(\mathbb R^2)\rightarrow L^{2}(\mathbb R^2)}\approx \|b_{n}\|_{\BMO_{\rm prod}(\mathbb R^2)}\geq C_n$ while
$\| T_n^*\|_{L^{2}(\mathbb R^2)\rightarrow L^{2}(\mathbb R^2)}\approx \|b_{n}\|_{\BMO_{\rm rect}(\mathbb R^2)}\leq C$.
This shows again that none of the conditions
$T_{1}(1)\in \BMO(\mathbb R^2)$,
$T_{1}^{*}(1)\in \BMO(\mathbb R^2)$
are necessary for boundedness of $T$.

\vskip10pt
We end this section by stating the analogous result in the multiparameter case. We simplify the notation as much as possible.

Let $m\leq n$ and $n_1,\ldots, n_m$ such that $n=\sum_{i=1}^mn_i$. Let $K:\prod_{i=1}^{m}(\mathbb R^{2n_i}\backslash \Delta_{2n_i})\rightarrow \mathbb R$  be such that
$$
|K(x,t)|\leq C \prod_{i=1}^m|x_i-t_i|^{-n_i}
$$
$$
|\nabla_{t_{i_1}}\cdots \nabla_{t_{i_m}}K(x,t)|\leq C \prod_{i=1}^m|x_i-t_i|^{-(n_i+\delta )}
$$
where $x_i,t_i\in \mathbb R^{n_i}$, $0<\delta \leq 1$.

\begin{definition}{(Restricted bilinear forms).}
Let $N_1,N_2\subset \{ 1,\ldots ,m\}$ such that $N_1\cup N_2= \{ 1,\ldots ,m\}$ disjointly.
Given a bilinear form $\Lambda $, we define the restricted bilinear forms by
$$
\langle \Lambda^{N_1}(\otimes_{j\in N_2}f_{j},\otimes_{j\in N_2}g_{j})\otimes_{i\in N_1}f_{i},\otimes_{i\in N_1}g_{i}\rangle
=\Lambda(f,g)
$$
for $f=\otimes_{i=1}^{m}f_i$, $g=\otimes_{i=1}^{m}g_{i}$ with $f_{i},g_{i}\in {\mathcal S}(\mathbb R^{n_i})$,
and then extended by linearity and continuity.

We will call restricted operators to the linear operators associated with the restricted bilinear forms by duality.
\end{definition}

Notice that the kernels of the forms $\Lambda^{N_1}$ depend on the variables of the functions $\otimes_{j\in N_2}f_{j}$, $\otimes_{j\in N_2}g_{j}$
and so we can write $\Lambda^{N_1}_{t_{j},x_{j}}$.

\begin{definition}
A bilinear form
$\Lambda:\S(\R^n)\times \S(\R^n)\to \C$
is said to be associated with a product Calder\'on-Zygmund kernel $K$ if it satisfies the following integral representations
$$
\Lambda(f,g)=\int_{\R^n}\int_{\R^n} f(t)g(x) K(x,t)\, dx \, dt
$$
for all Schwartz functions $f,g\in {\cal S(\R^n)}$
with disjoint support. 

\end{definition}

\begin{definition}
We say that a bilinear form $\Lambda $ satisfies the weak boundedness condition if for any box $R\subset \mathbb R^{n}$ and every $\phi_R, \varphi_R$
$L^2(\mathbb R^n)$-normalized bump functions adapted to $R$ with constant $C$ we have
$$
|\Lambda (\phi_R,\varphi_R)|\leq C
$$
\end{definition}

\begin{definition}
Let $W\!\!B,C\!Z\subset \{ 1,\ldots ,m\}$ such that $W\!\!B\cup C\!Z= \{ 1,\ldots ,m\}$ disjointly.
We say that a bilinear form $\Lambda $ satisfies the mixed weak boundedness-Calder\'on Zygmund condition
if for any $R=\prod_{i\in W\!\!B}R_{i}$ and every $\phi_R, \varphi_R$ bump functions $L^2(\mathbb R^{\sum_{i\in W\!\!B}n_i})$-normalized and
adapted to $R$ with constant $C>0$, we have
$$
|\Lambda^{W\!\!B}_{t_{i},x_{i}}(\phi_R,\varphi_R)|\leq C\prod_{i\in W\!B}|t_{i}-y_{i}|^{-1}
$$
$$
|(\Lambda^{W\!\!B}_{t_{i},x_{i}}-\Lambda^{W\!\!B}_{t'_{i},x'_{i}}) (\phi_R,\varphi_R)|
\leq C\prod_{i\in W\!B}(|t_{i}-t'_{i}|+|x_{i}-x'_{i}|)^{\delta }|t_{i}-x_{i}|^{-(n_j+\delta )}
$$
\end{definition}

\begin{theorem}{($m$-parameter $T(1)$ theorem).}\label{multiT1}
Let $\Lambda$ be a bilinear Calder\'on-Zygmund form with associated
kernel $K$ satisfying the mixed WB-CZ conditions.

Then, the following are equivalent:
\begin{enumerate}
\item $\Lambda_i$ are bounded bilinear forms on $L^2(\mathbb R^n)$ for all $i$,

\item $\Lambda $ satisfies
the weak boundedness condition and the following sequence of special cancellation conditions:
for every $k\in \{ n_1,n_1+n_2,\ldots, n\}$ and all bump functions $\phi_{R_i},\varphi_{R_i}\in {\mathcal S}(\mathbb R^{n_i})$ both adapted to $R_i\subset \mathbb R^{n_i}$
we have
$$
\langle T(1\otimes \ldots \otimes 1\otimes \Phi_{R^{n-k}}),\tilde{\Phi}_{R^{n-k}}\rangle \in \BMO (\mathbb R^{k})
$$
where $\Phi_{R^{n-k}}=\phi_{R_1}\otimes \ldots \otimes \phi_{R_{n-k}}$ and $\tilde{\Phi}_{R^{n-k}}=\varphi_{R_1}\otimes \ldots \otimes \varphi_{R_{n-k}}$
with norms uniformly bounded when varying over the boxes $R^{n-k}=\prod_{i=1}^{n-k}R_i$.

The same condition applies for all possible permutations of the entries in 
$1\otimes \ldots \otimes 1\otimes \phi_{R_1}\otimes \ldots \otimes \phi_{R_{n-k}}$
(in total ${\tiny (\begin{array}{c}n \\ n-k\end{array}) }\, 2^{k}$ conditions).
\end{enumerate}
\end{theorem}

Notice that the stated bi-parameter case corresponds with $n=2$, $m=2$ and $n_{i}=1$ and therfore 
$k\in \{ 1,2\}$.

\section{Application}\label{apply}

\vskip10pt

We give now an example of how our results can be applied to the study of boundedness of operators defined by product kernels.

In \cite{FS}, R. Fefferman and E. Stein explain that in some boundary-value problems, in particular in the $\bar{\partial}$-Neumann problem, one faces convolution operators defined
in $\mathbb R^{n+1}$ with kernels like
$$
K_{k}(t,t_{n+1})=\frac{t_k}{(|t|^2+t_{n+1}^2)^{(n+1)/2}}\frac{1}{|t|^2+it_{n+1}}
$$
with $t\in \mathbb R^n$ and $t_{n+1}\in \mathbb R$, which are product of two kernels with different types of homogeneity.
With this motivation in mind they prove the following

\begin{theorem}\label{convFS}
Let $K$ be a kernel defined in $\mathbb R^{n}\times \mathbb R^{m}$
by $K(t)=K_1(t)K_2(t)$ such that $K_1$ is homogeneous of degree $-n$ with respect the family of dilations $t\rightarrow (\delta t_1,\delta^at_2)$ for all
$\delta >0$ and fixed $a>0$, while $K_2$ is homogeneous of degree $-m$ with respect the family of dilations $t\rightarrow (\delta^b t_1,\delta t_2)$ for all
$\delta >0$ and fixed $b>0$.

It is also assumed that $K_1(t_1,0)$ has mean zero on the unit sphere of $\mathbb R^{n}$, $K_2(0,t_2)$ has mean zero on the unit sphere of $\mathbb R^{m}$ and
$$
\Big| \int_{\alpha_1<|t_1|<\beta_1,\alpha_2<|t_2|<\beta_2}K(t)dt\Big| \leq A
$$
for all $0<\alpha_i<\beta_i$.
Then, for all $1<p<\infty $,
$$
\| K*f\|_{L^p(\mathbb R^n\times \mathbb R^m)}\leq C_p\| f\|_{L^p(\mathbb R^n\times \mathbb R^m)}
$$
where the constant $C_p$ depends on $A$ and $p$.
\end{theorem}

In their paper, Theorem \ref{convFS} appears as a corollary of the following more general result:
\begin{theorem}\label{techconvFS}
Let $K:\mathbb R^n\times \mathbb R^m\rightarrow \mathbb R$ be an integrable function that satisfies
\begin{enumerate}
\item the kernel conditions: for $t=(t_1,t_2)$, $h=(h_1,h_2)$,
\begin{enumerate}
\item $|K(t)|\leq A|t_1|^{-n}|t_2|^{-m}$
\item $|K(t_1+h_1,t_2)-K(t)|\leq A|h|^{\delta_1}|t_1|^{-n-\delta_1}|t_2|^{-m}$ whenever $2|h_1|<|t_1|$
\item $|K(t_1,t_2+h_2)-K(t)|\leq A|h|^{\delta_2}|t_1|^{-n}|t_2|^{-m-\delta_2}$ whenever $2|h|<|t_2|$
\item $|K(t+h)-K(t_1+h_1,t_2)-K(t_1,t_2+h_2)+K(t)|
\leq A|h_1|^{\delta_1}|h_2|^{\delta_2}|t_1|^{-n-\delta_1}|t_2|^{-m-\delta_2}$ whenever $2|h_1|<|t_1|$, $2|h_2|<|t_2|$
\end{enumerate}

\item the cancellation condition: $\Big| \int_{\alpha_i<|t_i|<\beta_i}K(t)dt\Big| \leq A$ for all $0<\alpha_i <\beta_i$

\item the mixed kernel-cancellation conditions:
\begin{enumerate}
\item if $K_1(t_1)= \int_{\alpha_2<|t_2|<\beta_2}K(t_1,t_2)dt_2$ then
\begin{enumerate}
\item $|K_1(t_1)|\leq A|t_1|^{-n}$
\item $|K_1(t_1+h_1)-K_1(t_1)|\leq A|h|^{\delta_1}|t_1|^{-n-\delta_1}$ whenever $2|h_1|<|t_1|$
\end{enumerate}
\item similar conditions for $K_2(t_2)=\int_{\alpha_1<|t_1|<\beta_1}K(t_1,t_2)dt_1$.
\end{enumerate}
\end{enumerate}

Then, for all $1<p<\infty $,
$$
\| f*K\|_{L^{p}(\mathbb R^n\times \mathbb R^m)}\leq A_p \| f\|{L^p(\mathbb R^n\times \mathbb R^m)}
$$
with $A_p$ depending only on $A$ and $p$.
\end{theorem}

It is not difficult to see that that conditions of Theorem \ref{techconvFS} imply the hypotheses of Theorem \ref{biparT1}.
Conditions $(1)-a)$ and $(1)-d)$ imply $K$ is a product Calder\'on-Zygmund standard kernel while $(2)$ implies that the convolution operator $T$
with kernel $K$ satisfies weak boundedness condition and $T(1),T^{*}(1),T_1(1),T_{1}^{*}(1)\in \BMO(\mathbb R^{n}\times \mathbb R^m)$.
On the other hand, the mixed type hypotheses
of Theorem \ref{biparT1}, that is, the mixed
WB-CZ condition and mixed T(1)-CZ conditions follow from $(1)-b)$, $(1)-c)$; and $(3)$  respectively.

Even more, Theorem \ref{biparT1} allow us two extend Theorem \ref{techconvFS} to the case of non-convolution kernels, a result that is
stated below.
\begin{definition}\label{FSnonconvdef}
Let
$K(x,t)$ with $x,t\in \mathbb R^{n}\times \mathbb R^{m}$ be an integrable function that satisfies
\begin{enumerate}
\item the kernel conditions:
\begin{enumerate}
\item $|K(x,t)|\leq A|t_1|^{-n}|t_2|^{-m}$
\item $|K((x'_1,x_2),(t'_1,t_2))-K(x,t)|\leq A(|x_1-x'_1|+|t_1-t'_1|)^{\delta_1}|t_1|^{-n-\delta_1}|t_2|^{-m}$ whenever $2(|x_1-x'_1|+|t_1-t'_1|)<|t_1|$
\item $|K((x_1,x'_2),(t_1,t'_2))-K(x,t)|\leq A(|x_2-x'_2|+|t_2-t'_2|)^{\delta_2}|t_1|^{-n}|t_2|^{-m-\delta_2}$ whenever $2(|x_2-x'_2|+|t_2-t'_2|)|<|t_2|$
\item $|K(x',t')-K((x'_1,x_2),(t'_1,t_2))-K((x_1,x'_2),(t_1,t'_2))+K(x,t)|\\
\leq A(|x_1-x'_1|+|t_1-t'_1|)^{\delta_1}(|x_2-x'_2|+|t_2-t'_2|)^{\delta_2}|t_1|^{-n-\delta_1}|t_2|^{-m-\delta_2}$\\
whenever $2(|x_i-x'_i|+|t_i-t'_i|)<|t_i|$ for $i=1,2$
\end{enumerate}

\item the cancellation condition: $\Big| \int_{\alpha_1<|t_1|<\beta_1,\alpha_2<|t_2|<\beta_2}K(x,t)dt\Big| \leq A$

\item the mixed kernel-cancellation conditions:
\begin{enumerate}
\item if $K_1(x,t_1)= \int_{\alpha_2<|t_2|<\beta_2}K(x,t_1,t_2)dt_2$ then
\begin{enumerate}
\item $|K_1(x,t_1)|\leq A|t_1|^{-n}$
\item $|K_1((x'_1,x_2),t'_1)-K_1(x,t_1)|\leq A(|x_1-x'_1|+|t_1-t'_1|)^{\delta_1}|t_1|^{-n-\delta_1}$\\ whenever $2(|x_1-x'_1|+|t_1-t'_1|)<|t_1|$
\end{enumerate}
\item similar conditions for $K_2(x,t_2)=\int_{\alpha_1<|t_1|<\beta_1}K(x,t_1,t_2)dt_1$.
\end{enumerate}
\end{enumerate}
\end{definition}

\begin{definition}
We say that an operator $T$ is associated with $K$ if
$$
T(f)(x)=\int_{\mathbb R^n\times \mathbb R^m} f(x-t)K(x,t)dt
$$
whenever $x\notin \supp{(f)}$.
\end{definition}

Then, we have
\begin{theorem}\label{technonconvFS}
Let $T$  be an operator associated with $K$ satisfying all the conditions of definition \ref{FSnonconvdef}. Then,
for all $1<p<\infty $,
$$
\| T(f)\|_{L^{p}(\mathbb R^n\times \mathbb R^m)}\leq A_p \| f\|_{L^p(\mathbb R^n\times \mathbb R^m)}
$$
with $A_p$ depending only on $A$ and $p$.
\end{theorem}

{\it Sketch of proof.} To give a flavour of the ideas involved in dealing with non convolution kernels in the product seeting,
we outline how the 
hypotheses of Theorem \ref{technonconvFS} imply the ones in Theorem \ref{biparT1}.
In particular, we partially show how the mixed weak boundedness Calder\'on-Zygmund condition
$|\langle T^{1}_{x_1,t_1}(\phi_{I}),\varphi_{I}\rangle | \leq C|x_1-t_1|^{-n}$
and the cancellation property
$\langle T(\phi_{I}\otimes 1),\varphi_{I}\otimes \cdot \rangle \in \BMO(\mathbb R)$ are checked.

Let $I$ be a fixed interval. We consider
$\phi_{I}
=|I|^{-1/2}\sum_{k} a_{k}\chi_{I_{k}}$
to be an approximation by step functions of a general bump function adapted to $I$, where
the intervals $I_k$ are pairwise disjoint and of the same arbitrary small lenght $|I_k|<\epsilon $. We consider a similar description for $\varphi_{I}$.
Then,
$$
\langle T^{1}_{x_1,t_1}(\phi_{I}),\varphi_{I}\rangle
=\lim_{\epsilon \rightarrow 0}|I|^{-1}\sum_{k,j\in \mathbb Z}a_{k}b_{j}\int_{|x_2-t_2|>\epsilon} \chi_{I_{k}}(x_2)\chi_{I_{j}}(t_2)K(x,x-t)dt_2dx_2
$$
so we just need to bound
$$
|I|^{-1}\sum_{k,j\in \mathbb Z}a_{k}b_{j}T_{k,j}
$$
where $T_{k,j}$ denotes the integral in the sum,
independently of $\epsilon >0$ and $|I|$.

When $k=j$,
$$
T_{k,k}=\int_{\tiny \begin{array}{c}|x_2-c(I_{k})|<|I_{k}|/2\\|t_2-c(I_{k})|<|I_{k}|/2\\|x_2-t_2|>\epsilon \end{array}}K(x,x-t)dt_2dx_2
=\int_{\tiny \begin{array}{c}|x_2|<|I_{k}|/2\\|t_2|<|I_{k}|/2\\|x_2-t_2|>\epsilon \end{array}}K(x_1,x_2-c(I_k),x-t)dt_2dx_2
$$
$$
=\int_{\tiny \begin{array}{c}|x_2|<|I_{k}|/2\end{array}}
\int_{\tiny \begin{array}{c}|x_2-t_2|<|I_{k}|/2\\|t_2|>\epsilon\end{array}}K(x_1,x_2-c(I_k),x_1-t_1,t_2)dt_2dx_2
$$
$$
=\int_{\tiny \begin{array}{c}0<x_2<|I_{k}|/2\end{array}}
\Big( \int_{\tiny \begin{array}{c}|x_2-t_2|<|I_{k}|/2\\ |t_2|>\epsilon\end{array}}K(x_1,x_2-c(I_k),x_1-t_1,t_2)dt_2
$$
$$
+\int_{\tiny \begin{array}{c}|-x_2-t_2|<|I_{k}|/2\\|t_2|>\epsilon\end{array}}K(x_1,-x_2-c(I_k),x_1-t_1,t_2)dt_2\Big) dx_2
$$
By addition and substraction of constants, the tems between brackets in the last expression equal
{\small
\begin{equation}\label{adsubs}
\int_{\tiny \begin{array}{c}|x_2-t_2|<|I_{k}|/2\\ |t_2|>\epsilon\end{array}}K(x_1,x_2-c(I_k),x_1-t_1,t_2)dt_2
+\int_{\tiny \begin{array}{c}|-x_2-t_2|<|I_{k}|/2\\|t_2|>\epsilon\end{array}}K(x_1,x_2-c(I_k),x_1-t_1,t_2)dt_2
\end{equation}
$$
+\int_{\tiny \begin{array}{c}|-x_2-t_2|<|I_{k}|/2\\|t_2|>\epsilon\end{array}}
\big(K(x_1,-x_2-c(I_k),x_1-t_1,t_2)-K(x_1,x_2-c(I_k),x_1-t_1,t_2)\big)dt_2
$$}

The first two terms can be rewritten with a symmetric domain of integration as
$$
\int_{\tiny \begin{array}{c}|x_2-|I_{k}|/2|<|t_2|<|x_2+|I_{k}|/2|\end{array}}K(x_1,x_2-c(I_k),x_1-t_1,t_2)dt_2
$$
$$
+2\int_{\tiny \begin{array}{c}\epsilon <|t_2|<|x_2-|I_{k}|/2|\end{array}}K(x_1,x_2-c(I_k),x_1-t_1,t_2)dt_2
$$
where the second integral is zero if $|x_2-|I_{k}|/2|\leq \epsilon $. Then,
by the hypothesis 3.a.i) with $\alpha_2=|x_2-|I_{k}|/2|$, $\beta_2=|x_2+|I_{k}|/2|$ for the first integral and
$\alpha_2=\epsilon $, $\beta_2=|x_2-|I_{k}|/2|$ for the second one, we can bound them by $3A|x_1-t_1|^{-n}$.

Meanwhile, the last term in expression (\ref{adsubs}) can be treated by condition 1.c) and bounded by
$$
\int_{\tiny \begin{array}{c}|-x_2-t_2|<|I_{k}|/2\\|t_2|>\epsilon\end{array}}2|x_2|^{\delta_2}|x_1-t_1|^{-n}|t_2|^{-m-\delta_2}dt_2
\leq C |x_2|^{\delta_2}|x_1-t_1|^{-n}
$$

With both estimates we get
$$
\Big| \int_{\tiny \begin{array}{c}|x_2-c(I_{k})|<|I_{k}|/2\\|t_2-c(I_{k})|<|I_{k}|/2\\|x_2-t_2|>\epsilon \end{array}}K(x-t)dt_2dx_2\Big|
\leq \int_{\tiny \begin{array}{c}0<x_2<|I_{k}|/2\end{array}}3A(1+|x_2|^{\delta_2})|x_1-t_1|^{-n}dx_2
$$
$$
\leq 3A|x_1-t_1|^{-n}(1+|I_k|^{\delta })|I_{k}|/2
\leq 3A|x_1-t_1|^{-n}(1+\epsilon )|I_{k}|/2
$$

Finally, since we may assume $|b_k|\leq 1$ and $\sum_{k\in \mathbb Z}|a_{k}||I_{k}|\leq C|I|$, we finish this case with the following bounds:
$$
|I|^{-1}\sum_{k\in \mathbb Z}|a_{k}||b_{k}|T_{k,k}|
\leq C|I|^{-1}\sum_{k\in \mathbb Z}|a_{k}||b_{k}|CA|x_1-t_1|^{-n}|I_{k}|
\leq CA|x_1-t_1|^{-n}
$$

The case $k\neq j$ is technically more complex since
we need to consider several terms together in order to get the same kind of symmetry in the domain of integration.
Despite this, the same type of ideas apply: the kernel decay estimates allow to obtain a similar result and prove this way
the mixed WB-CZ condition.

\vskip10pt
On the other hand, let $\psi_J$ be an atom. Then,
$$
\langle T(\phi_{I}\otimes 1),\varphi_{I}\otimes \psi_{J} \rangle
$$
$$
=\lim_{\lambda \rightarrow \infty }\lim_{\epsilon \rightarrow 0}
\dint_{\tiny \begin{array}{c}|x_2-t_2|>\epsilon \end{array}}\chi_{\lambda I}(x_2)\psi_{J}(t_2) |I|^{-1}\sum_{k,j\in \mathbb Z}a_{k}b_{j}
\dint_{\tiny \begin{array}{c}|x_1-t_1|>\epsilon \end{array}} \chi_{I_{k}}(x_1)\chi_{I_{j}}(t_1)K(x,x-t)dtdx
$$
$$
=\lim_{\lambda \rightarrow \infty }\lim_{\epsilon \rightarrow 0}
\dint_{\tiny \begin{array}{c}|x_2-t_2|>\epsilon \end{array}}\chi_{\lambda I}(x_2)\psi_{J}(t_2)
|I|^{-1}\sum_{k,j\in \mathbb Z}a_{k}b_{j}T_{k,j}(x_2,t_2)dx_2dt_2
$$
and we bound the last expression independently of $\lambda >0$.
By using the mean zero of $\psi_{J}$ this is equal to
$$
\lim_{\epsilon \rightarrow 0}\dint_{\tiny \begin{array}{c}|x_2-t_2|>\epsilon \end{array}}\chi_{\lambda I}(x_2)\psi_{J}(t_2)
\Big( |I|^{-1}\sum_{k,j\in \mathbb Z}a_{k}b_{j}T_{k,j}(x_2,t_2)dx_2
-|I|^{-1}\sum_{k,j\in \mathbb Z}a_{k}b_{j}T_{k,j}(x_2,0)\Big) dx_2dt_2
$$
now by a similar argument as before but using the smoothness condition instead the decay we can bound by
$$
\lim_{\epsilon \rightarrow 0}\int |\psi_{J}(t_2)|
\int_{\tiny \begin{array}{c}|x_2-t_2|>\epsilon \end{array}}|\chi_{\lambda I}(x_2)|
CA|x_2-t_2|^{-(n+\delta )}dx_2dt_2
\leq CA\| \psi_{J}\|_{L^1(\mathbb R)}
$$

\section{Definition of $T(1)$, $\langle T(\phi_I\otimes 1),\varphi_I\rangle $ and $\langle T(\phi_I\otimes 1),\varphi_J\rangle $}\label{T(1)}

In this section we give a rigorous definition of $T(1)$, $\langle T(\phi_I\otimes 1),\varphi_I\rangle $ and
$\langle T(\phi_I\otimes 1),\varphi_J\rangle $ as distributions modulo constants.
The approach is similar to the uni-parametric case and so we will follow some of the arguments in \cite{ST}.

We start with the technical lemma that gives meaning to $T(1)$ (and also the partial adjoints $T_i(1)$).
The condition $T(1)\in \BMO(\mathbb R^2)$ means
that the following inequality
$$
|\langle T(1), f\rangle |\leq C\| f\|_{H^1(\mathbb R^2)}
$$
holds for all $f$ that belong to a dense subset of $H^{1}(\mathbb R^2)$. 

We recall that an atom in $\mathbb R^{2}$ is any function $f(x_1,x_2)$ supported on an open set $\Omega$ of finite measure, of the form
$$
   f= \sum_{ R \subset \Omega, R \text{ dyadic rectangle }}c_R f_R ,
$$
where each $f_R$ is a pre-atom on $R=I \times J$, that means, $f_R$ is $\mathcal{C}^1$, supported on $2R$,
$$
   \int f_R(x_1, x_2) dx_1 =0 \text{ for all }  x_2 \in \R, \quad \int f_R(x_1, x_2) dx_2 =0   \text{ for all }  x_1 \in \R,
$$
and
$$
   \|\partial_{x_1}^\alpha \partial_{x_2}^\beta f_R \|_\infty \le \frac{1}{|I|^{\alpha} |J|^{\beta} |R|^{1/2}}   \quad (0 \le \alpha, \beta \le 1).
$$
and such that
$$
    \sum_{ R \subset \Omega, R \text{ dyadic rectangle }}|c_R|^2  \le  \frac{1}{|\Omega|}.
$$
(see \cite{CF80, CF}).

However, in our case, such a dense subset will be
the family of Schwartz functions $f$ compactly supported with mean zero in each variable, meaning
$
\int_{\mathbb R}f(x,t)dx=\int_{\mathbb R}f(x,t)dt=0
$.

The way to apply such reduction is the following. We approximate any pre-atom $f$  by a Schwartz function $f_{n}$
with the described properties such that $f_{n}$ converges to $f$ in the topology of ${\mathcal S}(\mathbb R^{2})$
and $\| f_{n}\|_{L^{1}(\mathbb R)}\leq \| f\|_{L^{1}(\mathbb R)}$. The, we prove the below stated lemmata for $f_{n}$ but without the use of any smoothness property of $f_{n}$. Finally, the continuity of 
$\Lambda $
allow to conclude the same results for $f$. We will not give further details. 


In order to give a proper meaning to the left hand side of the previous inequality we
use Lemma \ref{definebmo1}:

\begin{lemma}\label{definebmo1}
Let $\Lambda $ be a Calder\'on-Zygmund bilinear form with associated kernel $K$.
Let $S$ be a rectangle and $f\in \S(\R^{2})$ with compact support in $S$, $L^{1}$-adapted to $S$
and with mean zero in each variable.
Let also $\Phi \in {\cal S}(\R^2)$ 
such that $\Phi (x)=1$ for $\| x\| \leq 1$ and $\Phi (x)=0$ for $\| x\| \geq 2$.

Then, the limit
$$L(f)=\lim_{k\to \infty} \Lambda(\mathcal T_{c(S)} D_{2^{k}|S_1|,2^{k}|S_2|}\Phi ,f)$$
exists.
Moreover, we have the error bound
$$|L(f)- \Lambda(\mathcal T_{c(S)} D_{2^{k}|S_1|,2^{k}|S_2|}\Phi ,f)|\le C 2^{-\delta k}
$$
where $\delta$ is the parameter in the Calderon-Zygmund property of
the kernel $K$ and the constant depends only on $\Phi$ and $\Lambda$.


\end{lemma}

\proof
For simplicity of notation we shall assume that $S$ is centered at the origin. 


For $k\in \mathbb N$, 
we set
$\psi_k=D_{2^{k}|S_1|,2^{k}|S_2|}\Phi -D_{2^{k-1}|S_1|,2^{k-1}|S_2|}\Phi$.
Let $C=\{ t\in \mathbb R^{2}: \min_{i}(|S_{i}|^{-1}|t_{i}|)\leq 1\}$. 
We decompose this function into two parts: 
$\psi_{k}^{in}=\psi_k \chi_{C}$ and $\psi_{k}^{out}=\psi -\psi_{k}^{1}$ and we bound
$|\Lambda(\psi_{k}, f)|$ by estimating each part separately. 

Since the supports of $\psi_{k}^{out}$ and $f$
are disjoint, we use the kernel
representation
$$
\Lambda(\psi_k^{out}, f)
=\int \psi_k^{out}(t)f(x) K(x,t)\, dtdx
$$
Due to the support of $f$ we may restrict the domain of integration to $|x_i|<|S_i|/2$. On the other hand, 
due to the support of $\psi_k$ we have that $2^{k-1}|S_i|<|t_i|\leq 2^{k+1}|S_i|$ 
for at least one coordinate $i\in \{ 1,2\}$ and $|S_i|<|t_i|\leq 2^{k+1}|S_i|$ for both coordinates. 

Using the mean zero of $f$ in each $x_i$ variable we rewrite the above integral as
$$
\int \psi_{k}^{out}(t)f(x)(K(x,t)-K((x_1,0),t)-K((0,x_2),t)+K(0,t))\, dtdx
$$
Since $2|x_i|<|S_i|<|t_i|$ we use the properties of product C-Z kernel to estimate the last display by
$$
C \int |\psi_{k}^{out}(t)||f(x)| \frac{|x_1|^\delta }{|t_1|^{1+\delta}}\frac{|x_2|^\delta }{|t_2|^{1+\delta}}\, dtdx
$$
Using again the restriction on the variables, we estimate now by
$$
C \int_{S}\int_{\tiny \begin{array}{l}|S_i|<|t_i|<2^{k+1}|S_i|\\ 2^{k-1}|S_{j}|<|t_j|<2^{k+1}|S_j|\end{array}} 
|f(x)|  \frac{|x_1|^\delta }{|t_1|^{1+\delta}}\frac{|x_2|^\delta }{|t_2|^{1+\delta}}
\, dtdx
$$
$$
=C \sum_{r=1}^{k}\int_{S}\int_{\tiny \begin{array}{l}2^{r-1}|S_i|<|t_i|<2^{r+1}|S_i|\\ 2^{k-1}|S_{j}|<|t_j|<2^{k+1}|S_j|\end{array}} 
|f(x)| \frac{|x_1|^\delta }{|t_1|^{1+\delta}}\frac{|x_2|^\delta }{|t_2|^{1+\delta}}\, dtdx
$$
$$
\leq C \sum_{r=1}^{k} \int_{S}
\int_{\tiny \begin{array}{l}2^{r-1}|S_i|<|t_i|<2^{r+1}|S_i|\\ 2^{k-1}|S_j|<|t_j|<2^{k+1}|S_j|\end{array}} 
|f(x)| \frac{1}{2^{r(1+\delta)}}\frac{1}{|S_1|}\frac{1}{2^{k(1+\delta)}}
\frac{1}{|S_2|}\, dtdx
$$
$$
\leq C \sum_{r=1}^{k} 2^{-r\delta }2^{-k\delta }\| f\|_{L^1(\mathbb R^2)}
\leq C 2^{-k\delta }\| f\|_{L^1(\mathbb R^2)}
$$

To estimate $|\Lambda(\psi_{k}^{in}, f)|$, we write $\psi_{k}^{in}=\psi_{k}^{in,1}\otimes \psi_{k}^{in,2}$. 
By symmetry, we need only to deal with the case 
when $|t_{2}|\leq |S_{2}|$. 
Then, by disjointness of the supports of $\psi_k^{in,1}$ and $f(\cdot ,t_{2})$ for all $t_{2}$, 
we can use the kernel representation of the restricted operator $T^1_{t_1,x_1}$
$$
\Lambda(\psi_k^{in}, f)
=\int \psi_k^{in,1}(t_1)\langle T^1_{t_1,x_1}(\psi_{k}^{in,2}),f(x_{1},\cdot )\rangle \, dt_1dx_1
$$
Due to the supports of $f$ and $\psi_k$ we have $|x_1|<|S_1|/2$ and 
$2^{k-1}|S_1|<|t_1|<2^{k+1}|S_1|$ respectively.
Using the mean zero of $f$ with respect the variable $x_{1}$, we write the above integral as
$$
\int \psi_k^{in,1}(t_1) \langle (T^1_{t_1,x_1}-T^1_{t_1,0})(\psi_{k}^{in,2}),f(x_{1},\cdot )\rangle \, dt_1dx_1
$$
Now $\psi_{k}^{in,2}$ is supported on $S_{2}$ and $L^{2}$-adapted to $S_{2}$ with constant 
$C|S_{2}|^{1/2}$. On the other hand, for every $x_{1}$, 
$f(x_1,\cdot ) $ is also supported on $S_{2}$ and $L^{2}$-adapted to $S_{2}$ with constant
$C\frac{|S_{2}|^{1/2}}{|S|}$. 
Then, since we also have $2|x_1|<|S_1|<|t_1|$, by the mixed WB-CZ properties we can estimate the last expression by
$$
C\frac{1}{|S_{1}|}\int_{|x_{1}|\leq |S_{1}|/2}\int_{2^{k-1}|S_1|<|t_1|\leq 2^{k+1}|S_1|}  |\psi_k^{in,1}(t_1)| 
\frac{|x_1|^\delta }{|t_1|^{1+\delta}}\, dt_1dx_1
$$
Finally, due to the restriction on the variables, we can estimate by
$$
C \frac{1}{|S_{1}|}\int_{|x_1|\leq |S_1|/2} 
\int_{|t_1|\leq 2^{k+1}|S_1|}\frac{1}{2^{k(1+\delta)}}\frac{1}{|S_1|}\, dt_1dx_1
\leq C 2^{-k\delta }
$$

These estimates prove that the sequence $(\Lambda(D_{2^{k}|S_1|,2^{k}|S_2|}\Phi ,f))_{k>0}$ is Cauchy and so the existence of the
limit $L(f)$.

Now the explicit rate of convergence stated in
the lemma follows by summing a geometric series: for every $k\in \mathbb N^2$, and every $0<\epsilon <2^{-k\delta }\| f\|_{L^1(\mathbb R^2)}$
let $m\in \mathbb N^2$ be with $|m|$ sufficiently large such that 
$|L(f)-\Lambda(D_{2^{m}|S_1|,2^{m}|S_2|}\Phi ,f)|<\epsilon $; then,
$$
|L(f)-\Lambda(D_{2^{m}|S_1|,2^{m}|S_2|}\Phi ,f)|
\leq |L(f)-\Lambda(D_{2^{m}|S_1|,2^{m}|S_2|}\Phi  ,f)|
+
\sum_{k'=k}^{m}|\Lambda(\psi_{k'} ,f)|
$$
$$
\leq \epsilon +C
\sum_{k'=k}^{m} 2^{-k'\delta }
\leq C2^{-k\delta }
$$
and the proof is finished. 

\vskip10pt
It can be easily proved that the definition of $T(1)$ is independent of the translation selected proving that $L$ is invariant under scaling and translation.
Moreover, it can also be shown that the definition is independent of the chosen cutoff function $\Phi $.

We notice that, since we have only worked with smooth atoms, strictly speaking we haven't finished the definition of $T(1)$. 
To do it rigorously, we should prove that the sequence $(T(D_{2^{k_1},2^{k_2}}\Phi ))_{k\in \mathbb Z^2}$
is uniformly bounded in $\BMO(\mathbb R^2)$. Then, using that the unit ball of the dual of Banach space is weak$^{*}$-compact, 
we can extract a subsequence of previous sequence which converges to $L(f)$ for  
functions $f$ in $C^{\infty }(\mathbb R^2)$ with compact 
support. Finally, since these functions are dense in $H^1(\mathbb R^2)$, we can deduce that previous functional can properly 
been extended to all $H^1(\mathbb R^2)$ and that $T(1\otimes 1)$ is the unique limit in $\BMO(\mathbb R^2)$ of the 
previous sequence. We will not get into any further detail about this issue.

\vskip10pt
We move now to the definition of $\langle T(\phi_I\otimes 1),\varphi_I\otimes \cdot\rangle $ following the previous schedule.
The condition $\langle T(\phi_I\otimes 1),\varphi_I\otimes \cdot \rangle \in \BMO(\mathbb R)$ means
the validity of the following inequality
$$
|\langle T(\phi_I\otimes 1),\varphi_I\otimes f\rangle |\leq C\| f\|_{H^1(\mathbb R)}
$$
for all $f$ that belong to a dense subset of $H^{1}(\mathbb R)$. In this case, such dense subset will be
the family of Schwartz functions $f$ compactly supported with mean zero.
Then, in order to give a proper meaning to the left hand side of previous inequality we
use the following Lemma:

\begin{lemma}\label{definebmo2}
Let $\Phi\in {\cal S}(\R )$ such that $\Phi(x)=1$ for $|x|\le 1$ and $\Phi(x)=0$ for $| x| \geq 2$.
Let $S$ be a rectangle and $\phi_{S_1}$, $\varphi_{S_1}$ be two $L^2$-normalized bump functions adapted to $S_1$.
Let $f\in \mathcal S(\mathbb R)$ be supported in $S_2$ with mean zero.
Then, the limit
$$
L_{\phi_{S_1},\varphi_{S_1}}(f)=\lim_{k\to \infty} \Lambda(\phi_{S_1}\otimes {\mathcal T}_{c(S_2)}D_{2^k|S_2|}\Phi ,\varphi_{S_1}\otimes f)
$$
exists.
Moreover, we have the error bound
$$|L_{\phi_{S_1},\varphi_{S_1}}(f)- \Lambda(\phi_{S_1}\otimes {\mathcal T}_{c(S_{2})}D_{2^k|S_2|}\Phi ,\varphi_{S_1}\otimes f)|\le C 2^{-\delta k}\| f\|_{L^1(\mathbb R)}$$
where $\delta$ is the parameter in the Calderon-Zygmund property of
the kernel $K$ and $C$ depends only on $\Phi$ and $\Lambda$.


\end{lemma}

\proof We mimic the proof of previous lemma and
for simplicity of notation we assume that $S_2$ is centered at the origin.
For $k\geq 1$, we set
$\psi_k=D_{2^k|S_2|}\Phi-D_{2^{k-1}|S_2|}\Phi$. We will estimate
$|\Lambda(\phi_{S_1}\otimes\psi_{k},\varphi_{S_1}\otimes f)|$.

Since the supports of $\psi_k$ and $f$ are disjoint
we use the kernel representation of the restricted operator $T^2_{t_2,x_2}$
$$
\Lambda(\phi_{S_1}\otimes\psi_{k},\psi_{S_1}\otimes \psi_{S_2})
=\int \psi_k(t_2)f(x_2) \langle T^2_{t_2,x_2}(\phi_{S_1}),\varphi_{S_1}\rangle \, dt_2dx_2
$$
Due to the supports of $f$ and $\psi_k$ we have $|x_2|<|S_2|/2$ and $2^{k-1}|S_2|<|t_2|<2^{k+1}|S_2|$ respectively.
Using the mean zero of $f$ we write the above integral as
$$
\int \psi_k(t_2)f(x_2) \langle (T^2_{t_2,x_2}-T^2_{t_2,0})(\phi_{S_1}),\varphi_{S_1}\rangle \, dt_2dx_2
$$
Since $2|x_2|<|S_2|<|t_2|$, by the mixed WB-CZ properties we can estimate the last expression by
$$
\int |\psi_k(t_2)||f(x_2)| C\frac{|x_2|^\delta }{|t_2|^{1+\delta}}\, dt_2dx_2
$$
and finally due to the restriction on the variables we can estimate by
$$
C \int_{|t_2|<2^k|S_2|} |f(x_2)|\frac{1}{2^{k(1+\delta)}}\frac{1}{|S_2|}\, dt_2dx_2
\leq C 2^{-k\delta }\| f\|_{L^1(\mathbb R)}
$$
As before this estimate is summable in $k$, which proves that
the sequence $(\Lambda(\phi_{S_1}\otimes D_{2^k|S_2|}\Phi ,\varphi_{S_1}\otimes f))_{k>0}$ is Cauchy and so, the
existence of the
limit $L(f)$. The explicit rate of convergence stated in
the lemma follows again by summing a geometric series.

\vskip10pt
Notice that the functional $L_{\phi_{S_1},\varphi_{S_1}}$
may also be denoted by $\Lambda(\phi_{S_1}\otimes 1,\varphi_{S_1}\otimes \cdot )$,
$\langle T(\phi_{S_1}\otimes 1),\varphi_{S_1}\otimes \cdot \rangle $,
$\langle \langle T^1(1),\cdot \rangle \phi_{S_1},\phi_{S_1}\rangle$ or $\langle T^2(\phi_{S_1}),\varphi_{S_1}\rangle 1$, since
\begin{eqnarray*}
\Lambda(\phi_{S_1}\otimes 1,\varphi_{S_1}\otimes f)
&=&\langle T(\phi_{S_1}\otimes 1),\varphi_{S_1}\otimes f\rangle\\
&=&\langle \langle T^1(1),f\rangle \phi_{S_1},\varphi_{S_1}\rangle\\
&=&\langle \langle T^2(\phi_{S_1}),\varphi_{S_1}\rangle 1,f\rangle
\end{eqnarray*}

Notice that in this way, the condition $L_{\phi_{S_1},\varphi_{S_1}}\equiv 0$
turns into $\langle T^2(\phi_{S_1}),\varphi_{S_1}\rangle \equiv 0$ for all $\phi_{S_1},\varphi_{S_1}$ adapted to $S_1$. On the other hand, the condition
$L_{\phi_{S_1},\varphi_{S_1}}\in \BMO (\mathbb R)$
turns into $\langle T^2(\phi_{S_1}),\varphi_{S_1}\rangle 1\in \BMO (\mathbb R)$ or $\langle T(\phi_{S_1}\otimes 1),\varphi_{S_1}\rangle_{x_1}(x_2)\in \BMO (\mathbb R)$
for all $\phi_{S_1},\varphi_{S_1}$ adapted to $S_1$.

\vskip15pt
Finally we define $\langle T(1\otimes \phi_I),\cdot \otimes \psi_J\rangle $ when $\phi_{I},\psi_{J}$ have disjoint support and $\psi_{J}$ has mean zero.
We follow a similar schedule as before by mixing the two previous cases.
The condition $\langle T(1\otimes \phi_I),\cdot \otimes \psi_J\rangle \in \BMO(\mathbb R)$ means
that
$$
|\langle T(1\otimes \phi_I),f\otimes \psi_J\rangle |\leq C\| f\|_{H^1(\mathbb R)}
$$
for all $f$ that belong to a dense subset of $H^{1}(\mathbb R)$. Again, the dense subset will be
the family of Schwartz functions $f$ compactly supported with mean zero.
Then, in order to give a proper meaning to the left hand side of previous inequality we
use the following lemma:

%
%
%

\begin{lemma}\label{definebmo3}
Let $\Phi\in {\cal S}(\R )$ such that $\Phi(x)=1$ for $|x|\le 1$ and $\Phi(x)=0$ for $| x| \geq 2$.
Let $\phi_{R_2}$, $\psi_{S_2}$ be $L^2$-normalized bump functions adapted and supported to the dyadic intervals $R_2$, $S_2$ respectively, such that
$|R_2|\geq |S_2|$, $|R_2|<\diam (R_2,S_2)$ and $\psi_{S_2}$ has mean zero.

Let $f\in \mathcal S(\mathbb R)$ be supported in a dyadic interval $S_1$ with mean zero.
Then, the limit
$$
L(f)=\lim_{k\to \infty} \Lambda({\mathcal T}_{c(S_1)}D_{2^k|S_1|}\Phi \otimes \phi_{R_2},f\otimes \psi_{S_2})
$$
exists.
Moreover, we have the error bound
$$
|L(f)- \Lambda(T_{c(S_1)}D_{2^k|S_1|}\Phi \otimes \phi_{R_2},f\otimes \psi_{S_2})|
\le C 2^{-\delta k}\Big(\frac{|S_2|}{|R_2|}\Big)^{1/2+\delta}(|R_2|^{-1}\diam(R_2\cup S_2))^{-(1+\delta)}
\| f\|_{L^1(\mathbb R)}
$$
where $\delta$ is the parameter in the Calderon-Zygmund property of
the kernel $K$ and $C$ depends only on $\Phi$ and $\Lambda$.


\end{lemma}

\proof Again
for simplicity of notation we assume that $S_1$ is centered at the origin.

For $k\geq 1$, we set
$\psi_k=D_{2^k|S_1|}\Phi-D_{2^{k-1}|S_1|}\Phi$. We will estimate
$|\Lambda(\psi_{k}\otimes \phi_{R_2},f\otimes \psi_{S_2})|$.

The supports of $\psi_k$ and $f$ and the supports of $\phi_{R_2}$ and $\psi_{S_2}$ are respectively disjoint (the latter because of the condition $|R_2|<\diam (R_2,S_2)$). Then,  the supports of $\psi_k\otimes \phi_{R_2}$  and $f\otimes \psi_{S_2}$  are also disjoint and so, 
we can use the kernel representation
$$
\Lambda(\psi_{k}\otimes \phi_{R_2},f\otimes \psi_{S_2})
=\int \psi_{k}(t_1)\phi_{R_2}(t_2)f(x_1)\psi_{S_2}(x_2) K(x,t)dtdx
$$
Due to the supports of the functions $\psi_{k}$ and $f$
we may restrict the domain of integration to $2^{k-1}|S_1|<|t_1|<2^{k+1}|S_1|$,
$|x_1|<|S_1|/2$ while, by hypothesis, we have  $|t_2-c(S_2)|>\diam (R_2\cup S_2) $, $|x_2-c(S_2)|<|S_2|/2$.

Using the mean zero of $f$ and $\psi_{S_2}$ we write the above integral as
$$
\int \psi_{k}(t_1)\phi_{R_2}(t_2)f(x_1)\psi_{S_2}(x_2) (K(x,t)-K((x_1,c(S_2)),t)-K((0,x_2),t)+K((0,c(S_2)),t)dtdx
$$
Since $2|x_1|<|S_1|<|t_1|$ and $2|x_2-c(S_2)|<|S_2|<|R_2|<\diam(R_2\cup S_2)<|t_2-c(S_2)|$, by the kernel properties we can estimate the last expression by
$$
\int_{\tiny{\begin{array}{l}|x_1|<|S_1|/2\\|t_1|<2^{k}|S_1|\end{array}}}
\int_{\tiny{\begin{array}{l}|x_2-c(S_2)|<|S_2|/2\\|x_2-t_2|>\diam(R_2\cup S_2)\end{array}}}
|\psi_{k}(t_1)||\phi_{R_2}(t_2)||f(x_1)||\psi_{S_2}(x_2)| C\frac{|x_1|^\delta }{|t_1|^{1+\delta}}
\frac{|x_2-c(S_2)|^\delta }{|t_2|^{1+\delta}}\, dtdx
$$
which, due to the restriction on the variables, we can estimate by
$$
C \int_{|t_1|<2^k|S_1|} |f(x_1)|\frac{1}{2^{k(1+\delta)}}\frac{1}{|S_1|}\, dt_1dx_1
\frac{|S_2|^{\delta}}{\diam(R_2\cup S_2)^{1+\delta}}\int |\phi_{R_2}(t_2)||\psi_{S_2}(x_2)|dt_2dx_2
$$
$$
= C 2^{-k\delta }\| f\|_{L^1(\mathbb R)}\frac{|S_2|^{\delta}}{\diam(R_2\cup S_2)^{1+\delta}}\| \phi_{R_2}\|_{L^1(\mathbb R)}\| \psi_{S_2}\|_{L^1(\mathbb R)}
$$
$$
\leq C 2^{-k\delta }\| f\|_{L^1(\mathbb R)}\frac{|S_2|^{\delta}}{\diam(R_2\cup S_2)^{1+\delta}}|R_2|^{1/2}|S_2|^{1/2}
$$
$$
\leq C 2^{-k\delta }\| f\|_{L^1(\mathbb R)}\Big(\frac{|S_2|}{|R_2|}\Big)^{1/2+\delta}
(|R_2|^{-1}\diam(R_2\cup S_2))^{-(1+\delta)}
$$

As before this estimate is summable in $k$, which proves that
the sequence $(\Lambda(D_{2^k|S_1|}\Phi \otimes \phi_{R_2},f\otimes \psi_{S_2}))_{k>0}$ is Cauchy and so, the
existence of the
limit $L(f)$. The explicit rate of convergence stated in
the Lemma follows again by summing a geometric series.

\section{$\Lambda$ applied to bump functions}\label{bump}

In this section we study the action of $\Lambda$ on bump functions
to obtain good estimates of the dual pair in terms of the space and
frequency (or scale) localization of the bump functions.

We start with several lemmata about localization properties of bump functions. 
Since these are essentially uni-parameter results, we will not include here their proofs, which can be otherwise found in \cite{PPV}. These lemmata
will be frequently used in the proof of 
proposition \ref{twobumplemma}, the main result of this section.
In particular, Lemma \ref{bump1} will be used when we apply the weak boundedness condition away
from the origin. On the other hand, Lemmata \ref{bump3} and \ref{bump4} will be mostly used when we need to use the cancellation condition $T(1)=0$ and the weak boundedness condition close to the origin.

\begin{lemma}\label{bump1}
Let $I,J$ be two intervals such that $|I|\geq |J|$. Let $0<\theta <1$,
$\lambda =(|J|^{-1}\diam (I\cup J))^{\theta }\geq 1$ and $\lambda J$ the interval with same center as $J$
and lenght $\lambda |J|$.

Let $\Phi_{\lambda J}$ be the usual $L^\infty$-normalized function adapted to $\lambda J$.
Let $\phi_{J}$ be a $L^2$-normalized bump function adapted to $J$ with constant $C$ and order $N$.

Then, $\phi_{J}(1-\Phi_{\lambda J})$ is a $L^2$-normalized bump function adapted to $I$ with constant
$$
C\Big(\frac{|J|}{|I|}\Big)^{\theta N/4-1/2}(|I|^{-1}\diam(I\cup J))^{-\theta N/2}
$$
and order $[\theta N/4]$.
\end{lemma}

\begin{lemma}\label{bump3}
Let $I,J$ be two intervals such that $|J|\leq |I|$. Let $N \in \mathbb N$, $0<\theta \leq 1$, 
$R\geq 3$ and 
$\lambda \geq R^{-1}(|J|^{-1}\diam (I\cup J))^{\theta }$.


Then, $|\lambda J|^{-1/2}\Phi_{\lambda J}$ is an $L^2$-normalized bump function adapted to $I$ with constant 
$$
CR^{2N}\Big(\frac{|J|}{|I|}\Big)^{\theta N/4-1/2}\lambda^{(N-1)/2}
$$
and order $[\theta N/4]$.
\end{lemma}

\begin{lemma}\label{bump4}
Let $I,J$ be two intervals such that $|c(I)-c(J)|<\max(|I|,|J|)$. Let $\phi_{J}$ be a bump function adapted to $J$ with 
constant $C>0$ and order $N$. 

Then,  $\phi_{J}$ is a bump function adapted to $I$ with constant 
$C\Big(\displaystyle{\frac{\max(|I|,|J|)}{\min(|I|,|J|)}}\Big)^{N+1/2}$ and order $N$.
\end{lemma}


\vskip10pt
Now we state and prove the result describing the action of the operator when it satisfies the special cancellation properties.
\begin{proposition}{(Bump lemma)}\label{twobumplemma}
Let $K$ be a product Calderon-Zygmund kernel with parameter $\delta$. Let $\Lambda$ be a bilinear Calderon-Zygmund form with associated
kernel $K$ which satisfies the mixed WB-CZ conditions.

Assume that $\Lambda $ also satisfies
the weak boundedness condition
and the special cancellation conditions $\Lambda (1\otimes 1,\psi_{S})=0$
for all $\psi \in {\mathcal S}(\mathbb R^2)$ with mean zero and
$\Lambda (f_1\otimes 1,g_1\otimes \psi )=\Lambda (1\otimes f_2,\psi \otimes g_2)=0$ 
for all $f_i,g_i,\psi \in {\mathcal S}(\mathbb R)$ with $\psi$ of mean zero.

For every $0<\theta <\min(2^{-15}, \delta /(1+\delta ))$, we choose $\frac{8}{\theta }<N<\frac{1}{16\theta^{2}}$
and define $0<\delta'=\delta-\theta (1+\delta )<\delta $.

Then,  
there exists $C_{\delta'}>0$  such that
for all rectangles $R$, $S$ such that $|R_i|\geq |S_i|$ for $i=1,2$ and all $\phi_1 $, $\psi_2$  bump functions $L^2$-adapted to
$R$ and $S$ respectively with order $N$ and constant $C>0$
such that $\psi_2$ has mean zero, we have
$$
|\Lambda(\phi_1,\psi_2)|\le C_{\delta'}\left(\frac{|S|}{|R|}\right)^{1/2+\delta'}
\prod_{i=1}^{2}\left( |R_i|^{-1}{\rm diam}(R_i\cup S_i) \right)^{-(1+\delta')}
$$
\end{proposition}

Notice that with some abuse of notation, whenever we use this estimate we will simply write $\delta $ instead of $\delta'$.

By symmetry on the arguments one can prove the following
\begin{corollary}\label{symmetricspecialcancellation}
Let $\Lambda $ a bilinear form that satisfies all the previously requested properties
and the following special cancellation conditions:
$$
\Lambda (1\otimes 1,\psi)=\Lambda (\psi,1\otimes 1)=\Lambda (\psi_{1}\otimes 1,1\otimes \psi_{2})
=\Lambda (1\otimes \psi_{2},\psi_{1}\otimes 1)=0
$$
for all $\psi \in {\mathcal S}(\mathbb R^2)$, $\psi_i \in {\mathcal S}(\mathbb R)$ with mean zero; and
$$
\Lambda (f_{1}\otimes 1,g_{1}\otimes \psi )
=\Lambda (1\otimes f_{2},\psi \otimes g_{2})
=\Lambda (\psi \otimes f_{2},1\otimes g_{2})
=\Lambda (f_{1}\otimes \psi,g_{1} \otimes 1)=0
$$ 
for all smooth functions $f_{i},g_{i}\in {\mathcal S}(\mathbb R)$ and $\psi \in {\mathcal S}(\mathbb R)$ of mean zero.

Let $0<\theta<1$, $N\in \mathbb N$ and $0<\delta' <\delta $ be as before. 
Then,  
there exists $C_{\delta'}>0$  such that
for every rectangles $R$, $S$ and $\psi_1 $, $\psi_2$ bump functions $L^2$-adapted to
$R$ and $S$ respectively with order $N$, constant $C>0$ and mean zero, we have
\begin{equation}\label{twobump}
|\Lambda(\psi_1,\psi_2)|\leq C_{\delta'} \prod_{i=1}^{2}\left(\frac{\min(|R_i|,|S_i|)}{\max (|R_i|,|S_i|)}\right)^{1/2+\delta'}
\left( \max(|R_i|,|S_i|)^{-1}{\rm diam}(R_i\cup S_i) \right)^{-(1+\delta')}
\end{equation}
\end{corollary}

{\it Proof of Proposition \ref{twobumplemma}.}
For simplicity of notation we shall assume that $S_i$ are both centered at the origin.
For each rectangle $R$ and $\lambda \in \mathbb R^2$, we denote by $\lambda R$ the dilated rectangle $(\lambda_1R_1)\times (\lambda_2R_2)$ that
shares the same centre as $R$ and has measure $\lambda_{1} \lambda_{2}|R|$.

Let $\Phi_{R_i}$ be the usual $L^\infty$-normalized function adapted to the interval $R_i$ and let $\Phi_{R}=\Phi_{R_1}\otimes \Phi_{R_2}$.

We denote $\psi (t,x)=\phi_1(t)\psi_2(x)$ which is adapted to $R\otimes S$ and truncate the function as follows.
We start by splitting $\psi $ in the $x_i$ variables iteratively, first in $x_1$ and later in $x_2$.
Let $\lambda_{i}=(|S_i|^{-1}{\rm diam}(R_i\cup S_i))^{\theta }$. Then,
$\psi=\psi_{in,-}+\psi_{out,-}$
where
$$
\psi_{in,-}(t,x)= (\psi(t,x) -c_1(t,x_2))\Phi_{\alpha \lambda_1 S_1}(x_1)
$$
and
$$
\psi_{out,-}(t,x)=\psi(t,x) (1-\Phi_{\lambda_1 S_1}(x_1))+c_1(t,x_2)\Phi_{\alpha \lambda_1 S_1}(x_1)
$$
where $\alpha =\frac{1}{32}$ and $c_1(t,x_2)$ is chosen so that $\psi_{in,-}$ and $\psi_{out,-}$ have mean zero in the variable $x_1$. Notice that both functions have
mean zero in the variable $x_2$. Now we decompose
$\psi_{in,-}=\psi_{in,in}+\psi_{in,out}$,
where
\begin{equation}\label{psiinin}
\psi_{in,in}(t,x)=(\psi_{in,-}(t,x)-c_2(t,x_1))\Phi_{\alpha\lambda_2 S_2}(x_2)
\end{equation}
and
$$
\psi_{in,out}(t,x)=\psi_{in,-}(t,x)(1-\Phi_{\frac{1}{32}\lambda_2 S_2}(x_2))+c_2(t,x_1)\Phi_{\alpha\lambda_2 S_2}(x_2)
$$
where  $c_2(t,x_1)$ is chosen so that both $\psi_{in,in}$ and $\psi_{in,out}$ have mean zero in the variable $x_2$. Meanwhile we also decompose
$\psi_{out,-}=\psi_{out,in}+\psi_{out,out}$,
where
\begin{equation}\label{psioutin}
\psi_{out,in}(t,x)=(\psi_{out,-}(t,x)-c_3(t,x_1))\Phi_{\alpha \lambda_2 S_2}(x_2)
\end{equation}
and
$$
\psi_{out,out}(t,x)=\psi_{out,-}(t,x)(1-\Phi_{\alpha \lambda_2 S_2}(x_2))+c_3(t,x_1)\Phi_{\alpha\lambda_2 S_2}(x_2)
$$
where $c_3(t,x_1)$ is chosen so that both $\psi_{out,in}$ and $\psi_{out,out}$ have mean zero in the $x_2$ variable.
Notice that for example
\begin{equation}\label{c2}
c_2(t,x_1)=-c|S_2|^{-1}\int \psi_{in,-}(t,x)(1-\Phi_{\alpha\lambda_2 S_2}(x_2))dx_2
\end{equation}

We see now that the four functions just constructed 
have mean zero in each variable $x_i$. This is obvious in the variable $x_2$, since $c_2$ and $c_3$ have been chosen to accomplish this.
Moreover, we know that $\psi_{in,-}$ and $\psi_{out,-}$ have mean zero in the variable $x_1$. Because of this, we have for each $x_2,t$:
\begin{equation}\label{intpsiinin}
\int \psi_{in,in}(t,x)dx_1\Phi_{\alpha \lambda_2 S_2}(x_2)
=-\int c_2(t,x)dx_1\Phi_{\alpha \lambda_2 S_2}(x_2)
\end{equation}
$$
=-c|S_2|^{-1}\int \psi_{in,-}(t,x)dx_1\int (1-\Phi_{\alpha \lambda_2 S_2}(x_2))\Phi_{\alpha \lambda_2 S_2}(x_2)dx_2
=0
$$
An analogous argument also proves mean zero of $\psi_{out, out}$ in each variable $x_i$. Mean zero of both
$\psi_{in,in}$, $\psi_{out, out}$ imply the same for $\psi_{in,out}$ and $\psi_{out,in}$.

Now we split the four functions in the $t_i$ variables. For $\psi_{out}$ we only decompose the first two terms to get
$\psi_{out,in}(t,x)=\psi_{1,2}(t,x) +\psi_{1,3}(t,x)$,
$\psi_{in,out}(t,x)=\psi_{2,1}(t,x)+\psi_{3,1}(t,x)$ and
$\psi_{out,out}(t,x)=\psi_{1,1}(t,x)$
where
$$
\psi_{1,2}(t,x)=\psi_{out,in}(t,x) \Phi_{\beta \mu_2 \tilde{R}_2}(t_2),
\hskip30pt
\psi_{1,3}(t,x)=\psi_{out,in}(t,x) (1-\Phi_{\beta \mu_2 \tilde{R}_2}(t_2))
$$
$$
\psi_{2,1}(t,x)=\psi_{in,out}(t,x) \Phi_{\beta \mu_1 \tilde{R}_1}(t_1),
\hskip30pt
\psi_{3,1}(t,x)=\psi_{in,out}(t,x) (1-\Phi_{\beta \mu_1 \tilde{R}_1}(t_1))
$$
$$
\psi_{1,1}(t,x)=\psi_{out,out}(t,x)
$$
with $\tilde{R}_i$ the translate of $R_i$ centered at the origin and $\beta =1/4$, $\mu_i=|R_i|^{-1}{\rm diam}(R_i\cup S_i)$.
The reason for notation $\psi_{i,j}$ will become clear later.

Finally, for $\psi_{in,in}$ we repeat the first type of decomposition to get the following four terms:
\begin{eqnarray*}
\psi_{in}(t,x)&=&\psi_{in,in}(x,t)\Phi_{\beta \mu \tilde{R}}(t) \\
&+&\psi_{in,in}(x,t)\Phi_{\beta \mu_1 \tilde{R}_1}(t_1)(1-\Phi_{\beta \mu_2 \tilde{R}_2}(t_2))\\
&+&\psi_{in,in}(x,t)(1-\Phi_{\beta \mu_1 \tilde{R}_1}(t_1))\Phi_{\beta \mu_2 \tilde{R}_2}(t_2)\\
&+&\psi_{in,in}(x,t)(1-\Phi_{\beta \mu_1 \tilde{R}_1}(t_1))(1-\Phi_{\beta \mu_2 \tilde{R}_2}(t_2))\\
&=&\psi_{2,2}(t,x)+\psi_{2,3}(t,x)+\psi_{3,2}(t,x)+\psi_{3,3}(t,x)
\end{eqnarray*}

A careful look at all these terms reveals that they can be described by
\begin{eqnarray*}
\psi_{1,2}(t,x)&=&\psi(t,x) (1-\Phi_{\alpha \lambda_1 S_1}(x_1))\Phi_{\alpha \lambda_2 S_2}(x_2)\Phi_{\beta \mu_2 \tilde{R}_2}(t_2)+c^1_1+c^1_3\\
\psi_{1,3}(t,x)&=&\psi(t,x) (1-\Phi_{\alpha \lambda_1 S_1}(x_1))\Phi_{\alpha \lambda_2 S_2}(x_2)(1-\Phi_{\beta \mu_2 \tilde{R}_2}(t_2))+c^2_1+c^2_3\\
\psi_{2,1}(t,x)&=&\psi(t,x) \Phi_{\alpha \lambda_1 S_1}(x_1)(1-\Phi_{\alpha \lambda_2 S_2}(x_2))\Phi_{\beta \mu_1 \tilde{R}_1}(t_1)+c^3_1+c^3_2\\
\psi_{3,1}(t,x)&=&\psi(t,x) \Phi_{\alpha \lambda_1 S_1}(x_1)(1-\Phi_{\alpha \lambda_2 S_2}(x_2))(1-\Phi_{\beta \mu_1 \tilde{R}_1}(t_1))+c^4_1+c^4_2\\
\psi_{1,1}(t,x)&=&\psi(t,x) (1-\Phi_{\alpha \lambda_1 S_1}(x_1))(1-\Phi_{\alpha \lambda_2 S_2}(x_2))+c^5_1+c^5_3\\
\psi_{2,2}(t,x)&=&\psi (t,x)\Phi_{\alpha \lambda S}(x)\Phi_{\beta \mu \tilde{R}}(t)+c^6_1+c^6_2\\
\psi_{2,3}(t,x)&=&\psi(t,x)\Phi_{\alpha \lambda S}(x)\Phi_{\beta \mu_1 \tilde{R}_1}(t_1)(1-\Phi_{\beta \mu_2 \tilde{R}_2}(t_2))+c^7_1+c^7_2\\
\psi_{3,2}(t,x)&=&\psi(t,x)\Phi_{\alpha \lambda S}(x)(1-\Phi_{\beta \mu_1 \tilde{R}_1}(t_1))\Phi_{\beta \mu_2 \tilde{R}_2}(t_2)+c^8_1+c^8_2\\
\psi_{3,3}(t,x)&=&\psi(t,x)\Phi_{\alpha \lambda S}(x)(1-\Phi_{\beta \mu_1 \tilde{R}_1}(t_1))(1-\Phi_{\beta \mu_2 \tilde{R}_2}(t_2))+c^9_1+c^9_3\\
\end{eqnarray*}
where the functions $c^{i}_j$ are error terms that ensure that all functions $\psi_{i,j}$ have mean zero in the variables $x_1,x_2$.
We notice that $c^{i}_1=c^{i}_1(t,x_2)$, $c^{i}_2=c^{i}_2(t,x_1)$ and $c^{i}_3=c^{i}_3(t,x_2)$.
At the end we will prove that the functions
$c^{i}_j$ are small and have the right support to allow us to assume that the main terms have the stated
zero averages. We will denote the main terms again by $\psi_{i,j}$.
Also notice that they are of tensor product type. Moreover,
with a small abuse of notation, we will write the action of the dual pair over $\psi_{i,j}$ as $\Lambda (\psi_{i,j})$.

\vskip10pt
We call $(1)$ the use of {\it weak boundedness condition} away from the origin, the mean zero in the variables $x_1,x_2$ and rate of decay of $\psi $;
$(2)$ the use of {\it the special cancellation condition} $T^i(1)\equiv 0$, the weak boundedness condition close to the origin and the mean zero of $\psi $ in the $x_i$ variable;
and $(3)$ the use of {\it the integral representation, the properties of the Calder\'on-Zygmund kernel} and the mean zero of $\psi $ in the variable $x_i$.
We call $(i)\times (j)$ the combined use of $(i)$ in the variables $t_1,x_1$
and $(j)$ in the variables $t_2,x_2$.
Then, we plan to bound $\Lambda (\psi )$ dealing each term $\Lambda (\psi_{i,j})$ by means of $(i)\times (j)$.

\vskip15pt
{\bf a) WB-WB.} We start with
$
\psi_{1,1}(t,x)
=\psi(t,x) (1-\Phi_{\alpha \lambda_1 S_1}(x_1))(1-\Phi_{\alpha \lambda_2 S_2}(x_2))
$
with mean zero in variables $x_1,x_2$, for which we will prove the decay by using weak boundedness in these variables.

We know that $\psi $ is adapted to $R\times S$ and so, by Lemma \ref{bump1} in the variables $x_{1}, x_{2}$,
$\psi_{1,1}$ is adapted to $R\times R$ with a gain in the constant of at least
$$
C \left(\frac{|S|}{|R|}\right)^{\theta N/4-1/2}\prod_{i=1,2}\left(|R_i|^{-1}{\rm diam}(R_i\cup S_i)\right)^{-\theta N/2}
$$
$$
\leq C \left(\frac{|S|}{|R|}\right)^{3/2}\prod_{i=1,2}\left(|R_i|^{-1}{\rm diam}(R_i\cup S_i)\right)^{-4}
$$
because $\theta $ and $N$ have been chosen so that $\theta N>8$.

Then, by the weak boundedness condition we have
$$
|\Lambda (\psi_{1,1})|\leq C\left(\frac{|S|}{|R|}\right)^{3/2}\prod_{i=1,2}\left(|R_i|^{-1}{\rm diam}(R_i\cup S_i)\right)^{-4}
$$

\vskip 10pt
{\bf b) CZ-CZ.} Now we consider
$
\psi_{3,3}(t,x)=\psi(t,x)\Phi_{\alpha \lambda S}(x)(1-\Phi_{\beta \mu_1 \tilde{R}_1}(t_1))(1-\Phi_{\beta \mu_2 \tilde{R}_2}(t_2))
$
which will be bounded by the integral representation and the properties of the CZ kernel.

On the support of the $\psi_{3,3}$, we have that $|t_i|>1/4\mu_i|R_i|={\rm diam}(R_i\cup S_i)/4$ while
$|x_i|\leq 2/32\lambda_i|S_i|=2/32|S_i|^{1-\theta }{\rm diam}(R_i\cup S_i)^{\theta }$. Since
$|S_i|^{-1}{\rm diam}(R_i\cup S_i)\geq 1$, the previous two inequalities imply $2|x_i|<|t_i|$
and so,
the support of $\psi_{3,3}$ is disjoint with the diagonal. This allow us to use the Calder\'on-Zygmund kernel representation
$$
\Lambda(\psi_{3,3})=\int \psi_{3,3}(t,x)K(x,t)\, dtdx
$$
Now using the mean zero of $\psi_{3,3}$ in the variable $x$, the above integral equals
$$
\int \psi_{3,3}(t,x)(K(x,t)-K((x_1,0),t)-K((0,x_2),t)+K(0,t))\, dtdx
$$
Moreover, since $2|x_i|<|t_i|$ we can use the property of a product CZ kernel and estimate the last expression by
$$
C \int |\psi_{3,3}(t,x)| \frac{|x_1|^\delta }{|t_1|^{1+\delta}}\frac{|x_2|^\delta }{|t_2|^{1+\delta}}\, dtdx
$$
Finally, since $|t_i|>{\rm diam}(R_i\cup S_i)/4$
and $|x_i|<2/32|S_i|^{1-\theta }{\rm diam}(R_i\cup S_i)^{\theta }$, we estimate by
$$
C\prod_{i=1,2}|S_i|^{(1-\theta )\delta }{\rm diam}(R_i\cup S_i)^{\theta \delta} 
{\rm \diam}(R_i\cup S_i)^{-1-\delta}\|\psi_{3,3}\|_1
$$
$$
\le C\prod_{i=1,2}|S_i|^{(1-\theta )\delta } {\rm diam}(R_i\cup S_i)^{-1-\delta+\theta \delta}|R_i|^{1/2}|S_i|^{1/2}
$$
$$
=C \prod_{i=1,2}(|S_i|/|R_i|)^{1/2+(1-\theta )\delta } (\diam(R_i\cup S_i)/|R_i|)^{-1-(1-\theta )\delta }
$$
$$
\leq C \prod_{i=1,2}(|S_i|/|R_i|)^{1/2+\delta'} (\diam(R_i\cup S_i)/|R_i|)^{-(1+\delta') }
$$
since $(1-\theta )\delta >\delta -\theta (1+\delta )=\delta'$. 
This proves the desired estimate for $\Lambda(\psi_{3,3})$.

\vskip 10pt
{\bf c) T(1)-T(1).} To bound,
$\psi_{2,2}(t,x)=\psi (t,x)\Phi_{\alpha \lambda S}(x)\Phi_{\beta \mu \tilde{R}}(t)$
we first argue the fact that we can make the extra assumption $\psi_{2,2}(0,x)=0$ for any $x$.

The assumption comes from the substitution of $\psi_{2,2}(t,x)$ by
\begin{equation}\label{subtraction1}
\psi_{2,2}(t,x)-
{\mathcal D}_{\frac{|R|}{|S|}4\alpha \lambda_{1}|S_1|,\frac{|R|}{|S|}4\alpha \lambda_{2}|S_2|}\Phi(t)\psi_{2,2}(0,x)
\end{equation}
Actually, in order to prove $\tilde{\psi}_{2,2}$ is adapted to certain rectangle, we also need to subtract a term including some partial derivatives of $\psi_{2,2}(0,x)$. However, in order to simplify the exposition, we will obvious this fact. For more details, we refer to \cite{PPV}. 

We first need to prove that the subtracted term satisfies the same bounds we want to prove.
Denote $\tilde{\psi}_{2,2}(t,x)=\psi_{2,2}(0,x)$. 
Since $\psi_{2,2} $ is adapted to $R\times S$ with constant comparable to $C$, we have that 
$$
|\tilde{\psi}_{2,2}(t,x)|\leq C|R|^{-1/2}\prod_{i=1,2}(1+|R_{i}|^{-1}|c(R_{i})-c(S_{i})|)^{-N}|\varphi_{S}(x)|
$$
with $\varphi_{S}$ an $L^{2}$ normalized bump function adapted to $S$. Moreover, since 
$$
1+|R_{i}|^{-1}|c(R_{i})-c(S_{i})|=|R_{i}|^{-1}(|R_{i}|+|c(R_{i})-c(S_{i})|)
$$
$$
\geq |R_{i}|^{-1}\diam(R_{i}\cup S_{i})
$$
we have 
$$
|\tilde{\psi}_{2,2}(t,x)| \leq C|R|^{-1/2}\prod_{i=1,2}(|R_i|^{-1}\diam(R_i\cup S_i))^{-N}|\varphi_{S}(x)|
$$
Moreover, $\tilde{\psi}_{2,2}$ is supported in $\frac{4}{32}\lambda_{2}S$ with mean zero.
Then, by the special cancellation condition $\Lambda (1\otimes 1,\psi_{S_1}\otimes \psi_{S_2})=0$, the explicit error of Lemma \ref{definebmo1}
and the decay of $\tilde{\psi}$ just calculated, 
we can estimate the contribution of the term subtracted in (\ref{subtraction1}) by

$$
|\Lambda({\mathcal D}_{\frac{|R|}{|S|}\alpha \lambda_{1}|S_1|,\frac{|R|}{|S|}\alpha \lambda_{2}|S_2|}\Phi,
\tilde{\psi}_{2,2})|
=|\Lambda({\mathcal D}_{\frac{|R|}{|S|}\alpha \lambda_{1}|S_1|,\frac{|R|}{|S|}\alpha \lambda_{2}|S_2|}\Phi,
\tilde{\psi}_{2,2})-\Lambda (1\otimes 1,\tilde{\psi}_{2,2})|
$$
$$
\leq C
\left(\frac{|R|}{|S|}\right)^{-\delta}
\| \tilde{\psi}_{2,2} \|_{L^{1}(\mathbb R^{2})}
\leq C\left(\frac{|S|}{|R|}\right)^{\delta}|R|^{-1/2}\prod_{i=1,2}(|R_i|^{-1}\diam(R_i\cup S_i))^{-N}|S|^{1/2}
$$
$$
=C\left(\frac{|S|}{|R|}\right)^{\delta+1/2}\prod_{i=1,2}(|R_i|^{-1}\diam(R_i\cup S_i))^{-N}
$$
and the right hand side is no larger than the desired bound.

Now, with the assumption, we further decompose $\psi_{2,2}$ with respect the $t$ variables in the following way:
$$
\psi_{2,2}=\psi_{2,2,in}+\psi_{2,2,out}
$$
$$
\psi_{2,2,in}(t,x)= \psi_{2,2}(t,x)\Phi_{\frac{8}{32}\lambda S}(t)
$$

c1) We first prove that $\psi_{2,2,in}$ is adapted to $\lambda S\times \lambda S$ with constant
$$
C\left(\frac{|S|}{|R|}\right)^{3/4}
\prod_{i=1,2}(|R_i|^{-1}{\rm diam}(R_i\cup S_i))^{-7}
$$ 
and order $N_{2}=[N^{1/2}]$.
Actually, in order to shorten the proof, we will only show how to obtain the bound for the function. For the corresponding work to estimate also the derivatives we refer again to \cite{PPV}. 

On one side, we have, 
$$
|\Phi_{\frac{8}{32}\lambda S}(t)|
\lesssim |\lambda S|^{1/2}|\lambda S|^{-1/2}
(1+|\lambda S|^{-1}|t-c(\lambda S)|)^{-N_{2}}
$$
which shows that $\Phi_{\frac{8}{32}\lambda S}$ is $L^{2}$-adapted to $\lambda S$ with constant comparable to  $|\lambda S|^{1/2}$.

On the other hand, the support  of $\psi_{2,2,in}$ in the variables $t_{i}$ is in $\lambda S$. Then, for all  $t\in \lambda S$ and all $x\in \mathbb R^{2}$, we have 
by the extra assumption, 
$$
|\psi_{2,2,in}(t,x)|\leq |\psi_{2,2}(t,x)|
=\Big| \int_{0}^{t_{1}} \int_{0}^{t_{2}} \partial_{t_{1}}\partial_{t_{2}} \psi_{2,2}(r,x)\, dr\Big| 
$$
$$
\le |t_{1}||t_{2}|\| \partial_{t_{1}}\partial_{t_{2}}\psi_{2,2}(\cdot ,x)\|_{\infty }
\leq \lambda_{1}\lambda_{2} |S|\| \partial_{t_{1}}\partial_{t_{2}}\psi_{2,2}(\cdot ,x)\|_{\infty }
$$
By definition of bump function
$$
| \partial_{t_1}\partial_{t_2} \psi_{2,2}(r,x)|\leq 
C|R|^{-3/2}\prod_{i=1,2}(1+|R_i|^{-1}|r_{i}-c(R_i)|)^{-N_{2}}
\varphi_{S} (x)
$$
where $\varphi_{S} $ is an $L^2$-normalized bump function adapted to $S$
because $\psi_{2,2} (r,\cdot )$ is 
adapted to $S$. 
Now, from the particular choice of $\Phi_{\frac{8}{32}\lambda S}$, we have
$$
|r_{i}|\leq |t_{i}|\leq 1/2\lambda_{i}|S_{i}|
$$
$$
=1/2|S_{i}|^{1-\theta}\diam(R_{i}\cup S_{i})^{\theta}\leq 1/2\diam(R_{i}\cup S_{i})
$$
because $|S_{i}|\leq \diam(I_{i}\cup J_{i})$, and so, 
$$
1+|R_{i}|^{-1}|r_{i}-c(R_{i})|\geq 1+|R_{i}|^{-1}|c(R_{i})|-|R_{i}|^{-1}|r_{i}|
$$
$$
\geq |R_{i}|^{-1}\diam(R_{i}\cup S_{i})-1/2|R_{i}|^{-1}\diam(R_{i}\cup S_{i})
=1/2|R_{i}|^{-1}\diam(R_{i}\cup S_{i})
$$
Therefore, 
\begin{equation}\label{decaypsi0}
\| \partial_{t_1}\partial_{t_2} \psi_{2,2}(\cdot ,x)\|_{\infty }\leq 
C|R|^{-3/2}\prod_{i=1,2}(|R_i|^{-1}{\rm diam}(R_i\cup S_i))^{-N_{2}}
\varphi_{S} (x)
\end{equation}

Moreover, by lemma \ref{bump4} in the variables $x_{1}, x_{2}$, we have that $\varphi_{S} $ is a bump function
adapted to $\lambda S$ with constant $C\lambda_{1}^{N_{2}+1/2}\lambda_{2}^{N_{2}+1/2}$ and so
$$
|\varphi_{S} (x)|
\leq C\lambda_{1}^{N_{2}+1/2}\lambda_{2}^{N_{2}+1/2}
\prod_{i=1,2}|\lambda_{i}S_{i}|^{-1/2}(1+|\lambda_{i}S_{i}|^{-1}|x-c(S_{i})|)^{-N_{2}}
$$
This  
implies that  we can bound in the following way
$$
|\psi_{2,2}(t,x)|
\leq C \lambda_{1}\lambda_{2}|S||R|^{-3/2}\prod_{i=1,2}(|R_i|^{-1}{\rm diam}(R_i\cup S_i))^{-N}
$$
$$
\lambda_{1}^{N_{2}+1/2}\lambda_{2}^{N_{2}+1/2}
\prod_{i=1,2}|\lambda_{i}S_{i}|^{-1/2}(1+|\lambda_{i}S_{i}|^{-1}|x-c(S_{i})|)^{-N_{2}}
$$
$$
=C |\lambda S|^{-1/2}\left(\frac{|S|}{|R|}\right)^{3/2}\prod_{i=1,2}(|R_i|^{-1}{\rm diam}(R_i\cup S_i))^{-N}
$$
$$
\lambda_{1}^{N_{2}+2}\lambda_{2}^{N_{2}+2}
\prod_{i=1,2}|\lambda_{i}S_{i}|^{-1/2}(1+|\lambda_{i}S_{i}|^{-1}|x-c(S_{i})|)^{-N_{2}}
$$

Therefore, 
$$
|\psi_{2,2,in}(t,x)|= |\psi_{2,2}(t,x)| |\Phi_{\frac{8}{32}\lambda S}(t)|
$$
$$
\leq C \left(\frac{|S|}{|R|}\right)^{3/2}\prod_{i=1,2}(|R_i|^{-1}{\rm diam}(R_i\cup S_i))^{-N}
\lambda_{1}^{N_{2}+2}\lambda_{2}^{N_{2}+2}
$$
$$
\prod_{i=1,2}|\lambda_{i}S_{i}|^{-1/2}(1+|\lambda_{i}S_{i}|^{-1}|x-c(S_{i})|)^{-N_{2}}
\cdot |\lambda S|^{-1/2}|\Phi_{\frac{8}{32}\lambda S}(t)|
$$
and notice that $|\lambda S|^{-1/2}\Phi_{\frac{8}{32}\lambda S}(t)$ is $L^{2}$-adapted to $\lambda S$. This already shows that $\psi_{2,2,in}(t,x)$ is adapted to $\lambda S\times \lambda S$. But we know simplify the constant. 
Since
$$
\lambda_{i}^{N_{2}+2}=(|S_{i}|^{-1}\diam(R_{i}\cup S_{i}))^{\theta (2+N_{2})}
$$
$$
\leq \Big(\frac{|S_{i}|}{|R_{i}|}\Big)^{-\theta (2+N_{2})}(|R_{i}|^{-1}\diam(R_{i}\cup S_{i}))^{\theta (2+N_{2})}
$$
we bound the constant by 
$$
C \left(\frac{|S|}{|R|}\right)^{3/2-\theta (2+N_{2})}
\prod_{i=1,2}(|R_i|^{-1}{\rm diam}(R_i\cup S_i))^{-(N-\theta (2+N_{2}))}
$$
$$
\leq C \left(\frac{|S|}{|R|}\right)^{3/4}
\prod_{i=1,2}(|R_i|^{-1}{\rm diam}(R_i\cup S_i))^{-7}
$$
because, from the choice of $\theta$ and $N$, we have 
$\theta^{2} N<\frac{1}{16}$ and so, $\theta N_{2}<\frac{1}{4}$. This implies 
$\theta (2+N_{2})\leq 2\theta (1+N_{2})<2\theta +1/2<3/4$. 


Then, by the weak boundedness property of $\Lambda $ we get
$$
|\Lambda(\psi_{2,2,in})|
\leq C \left(\frac{|S|}{|R|}\right)^{3/4}\prod_{i=1,2}(|R_i|^{-1}{\rm diam}(R_i\cup S_i))^{-7}
$$

c2) We now work with $\psi_{2,2,out}$.
When $\frac{8}{32}\lambda_{i} |S_{i}|>\frac{1}{2}\mu_{i}|R_{i}|$, we have that 
$\psi_{2,2,out}$ is the zero function and so, the right decay of $|\Lambda (\psi_{2,2,out})|$ is 
trivially satisfied. Hence, we only need to work in the case when 
$\frac{8}{32}\lambda_{i}|S_{i}|\leq \frac{1}{2}\mu_{i}|R_{i}|$.

By the extra assumption again, we have as before
$$
|\psi_{2,2, out}(t,x)|\leq |\psi_{2,2}(t,x)|
=\Big| \int_{0}^{t_{1}} \int_{0}^{t_{2}}\partial_{t_{1}}  \partial_{t_{2}}\psi_{2,2}(r,x)\, dr\Big| 
$$
$$
\le C|t_{1}||t_{2}|
\|\partial_{t_{1}}  \partial_{t_{2}}\psi_{2,2}(\cdot ,x)\|_{\infty }
$$
Again, by the definition of a bump function we have
$$
| \partial_{t_1}\partial_{t_2} \psi_{2,2}(r,x)|\leq 
C|R|^{-3/2}\prod_{i=1,2}(1+|R_i|^{-1}|r_{i}-c(R_i)|)^{-N}
\varphi_{S} (x)
$$
where $\varphi_{S} $ is an $L^2$-normalized bump function adapted to $S$
because $\psi_{2,2} (r,\cdot )$ is 
adapted to $S$. 
But now, from the particular choice of $\Phi_{\frac{1}{4}\mu_{i}\tilde{R}_{i}}$, we have
$$
|r_{i}|\leq |t_{i}|\leq  \mu_{i}|\tilde{R}_{i}|=1/4\mu_{i}|R_{i}|
=1/2\diam(R_{i}\cup S_{i})
$$
and so, 
$$
1+|R_{i}|^{-1}|r_{i}-c(R_{i})|\geq 1+|R_{i}|^{-1}|c(R_{i})|-|R_{i}|^{-1}|r_{i}|
$$
$$
\geq |R_{i}|^{-1}\diam(R_{i}\cup S_{i})-1/2|R_{i}|^{-1}\diam(R_{i}\cup S_{i})
=1/2|R_{i}|^{-1}\diam(R_{i}\cup S_{i})
$$
Therefore, 
$$
\| \partial_{t_1}\partial_{t_2} \psi_{2,2}(\cdot ,x)\|_{\infty }\leq 
C|R|^{-3/2}\prod_{i=1,2}(|R_i|^{-1}{\rm diam}(R_i\cup S_i))^{-N}
\varphi_{S} (x)
$$
and so, 
\begin{equation}\label{decaypsi0}
|\psi_{2,2, out}(t,x)|\leq C|t_{1}||t_{2}|
|R|^{-3/2}\prod_{i=1,2}(|R_i|^{-1}{\rm diam}(R_i\cup S_i))^{-N}\varphi_{S} (x)
\end{equation}
%

On the support of $\psi_{2,2}$, we have $|t_{i}|\leq \frac{1}{2}\mu_{i}|R_{i}|$.
Moreover, on the support of 
$\psi_{2,2,out}$, we have that 
$|t_{i}|>\frac{8}{32}\lambda_{i}|S_{i}|$ while
$|x_{i}|\leq \frac{2}{32}\lambda_{i}|S_{i}|$. 
Then, $2|x_{i}|<|t_{i}|$
and so,
we can use the Calder\'on Zygmund kernel representation
$$
\Lambda (\psi_{2,2,out})=\int \psi_{2,2,out}(t,x) K(x,t)\, dtdx
$$
Using the mean zero of $\psi_{2, 2, out}$ in the variables $x_{i}$, we can rewrite the above integral as
$$
\int \psi_{2, 2, out}(t,x) (K(x,t)-K(0,t))\, dtdx
$$
Since $2|x_{i}|<|t_{i}|$, we can use the smooth property of a product Calder\'on-Zygmund kernel 
and the decay of the bump function of $\psi_{2, 2, in}$ calculated in (\ref{decaypsi0}), 
to estimate 
$$
|\Lambda (\psi_{2,2,out})|\leq 
C \int 
|t_{1}| |t_{2}|||R|^{-3/2}
\prod_{i=1,2}(|R_i|^{-1}{\rm diam}(R_i\cup S_i))^{-N}
|\varphi_{S} (x)| 
\frac{|x_{1}|^\delta}{|t_{1}|^{1+\delta}}\frac{|x_{2}|^\delta}{|t_{2}|^{1+\delta}}\, dtdx
$$
Now, using the restrictions $|x_{i}|<\frac{2}{32}\lambda_{i}|S_{i}|$ and  $\frac{8}{32}|\lambda_{i}S_{i}|<|t_{i}|\leq \frac{1}{2}\mu_{i}|R_{i}|$, we can estimate previous expression by 
$$
C
|R|^{-3/2}\prod_{i=1,2}(|R_i|^{-1}{\rm diam}(R_i\cup S_i))^{-N}
|S|^{-1/2}
$$
$$
\int_{|x_{i}|<\frac{2}{32}\lambda_{i}|S_{i}|} |x_{1}|^\delta |x_{2}|^\delta dx 
\int_{\frac{8}{32}|\lambda_{i}S_{i}|<|t_{i}|\leq \frac{1}{2}\mu_{i}|R_{i}|} \frac{1}{|t_{1}|^{\delta}}\frac{1}{|t_{2}|^{\delta}}dt
$$
Now, in order to compute the integrals, we assume that $\delta <1$. The case $\delta =1$ can be treated in the same way leading to a similar estimate. Then, the product of the integrals is bounded by a constant times 
$$
\prod_{i=1,2}(\frac{2}{32}\lambda_{i}|S_{i}|)^{1+\delta }
((\frac{1}{2}\mu_{i}|R_{i}|)^{1-\delta}-(\frac{8}{32}\lambda_{i}|S_{i}|)^{1-\delta})
\lesssim (\lambda_{i}|S_{i}|)^{1+\delta }(\mu_{i}|I|)^{1-\delta}
$$
$$
=\prod_{i=1,2}|S_{i}|^{(1-\theta)(1+\delta )}
\diam(R_{i}\cup S_{i})^{\theta(1+\delta )}\diam(R_{i}\cup S_{i})^{1-\delta }
$$
$$
=\prod_{i=1,2}|S_{i}|^{1+\delta -\theta(1+\delta )}|R_{i}|^{1-\delta +\theta(1+\delta )}
(|R_{i}|^{-1}\diam(R_{i}\cup S_{i}))^{1-\delta +\theta(1+\delta )}
$$

Therefore, we can bound 
%
%
$$
|\Lambda (\psi_{2,2,out})|
\leq |R|^{-3/2}\prod_{i=1,2}(|R_i|^{-1}{\rm diam}(R_i\cup S_i))^{-N}|S|^{-1/2}
$$
$$
\prod_{i=1,2}|S_{i}|^{1+\delta -\theta(1+\delta )}|R_{i}|^{1-\delta +\theta(1+\delta )}
(|R_{i}|^{-1}\diam(R_{i}\cup S_{i}))^{1-\delta +\theta(1+\delta )}
$$
$$
=C \Big(\frac{|S|}{|R|}\Big)^{1/2+\delta -\theta (\delta +1)} \prod_{i=1,2}(|R_i|^{-1}{\rm diam}(R_i\cup S_i))^{-(N-1)} 
$$
$$
\leq C \Big(\frac{|S|}{|R|}\Big)^{1/2+\delta'} \prod_{i=1,2}(|R_i|^{-1}{\rm diam}(R_i\cup S_i))^{-7} 
$$

\vskip10pt
Once finished the three "pure" cases, we move to the proof of the "mixed" ones. 
\vskip10pt 

{\bf d) WB-CZ. } We work 
$
\psi_{1,3}(t,x)=\psi(t,x) (1-\Phi_{\alpha \lambda_1 S_1}(x_1))\Phi_{\alpha \lambda_2 S_2}(x_2)(1-\Phi_{\beta \mu_2 \tilde{R}_2}(t_2))
$
which will be dealt by using the weak boundedness in the variables $t_1,x_1$ and the kernel representation in the variables $t_2,x_2$.

By lemma \ref{bump1} in the variable $x_1$, we have that $\psi_{1,3}$ is adapted to $(R_1\times R_2)\times (R_1\times S_2)$ with a gain in the constant of at least
$$
C \left(\frac{|S_1|}{|R_1|}\right)^{3/2}\left(|R_1|^{-1}{\rm diam}(R_1\cup S_1)\right)^{-4}
$$
This way we can assume
$\psi_{1,3}(t,x)=\phi_{R_1}(t_1)\phi_{R_2}(t_2)\varphi_{R_1}(x_1)\psi_{S_2}(x_2)$ with $\psi_{S_2}$ of zero mean.
Moreover, on the support of $\psi_2$ we have
that $|t_2|>\diam (S_2\cup R_2)$ while $|x_2|<\frac{2}{32}|S_2|^{1-\theta }\diam(R_2\cup S_2)^{\theta}$. This implies $|x_2|<|t_2|$ and so by
the integral representation of the restricted operator $T^2_{t_2,x_2}$, we have
$$
\Lambda(\psi_{1,3})
=\int \phi_{R_2}(t_2)\psi_{S_2}(x_2) \langle T^2_{t_2,x_2}(\phi_{R_1}),\varphi_{R_1}\rangle \, dt_2dx_2
$$
Using the mean zero of $\psi_{S_2}$ we obtain for the above integral
$$
\int \phi_{R_2}(t_2)\psi_{S_2}(x_2) \langle (T^2_{t_2,x_2}-T^2_{t_2,0})(\phi_{R_1}),\varphi_{R_1}\rangle \, dt_2dx_2
$$
Since $2|x_2|<|t_2|$, by the mixed WB-CZ property of $\Lambda $ and the gain in the constant, we estimate the integral by
$$
\int |\psi_{R_2}(t_2)||\psi_{S_2}(x_2)| C\Big( \frac{|S_1|}{|R_1|}\Big)^{3/2}
(|R_1|^{-1}{\rm diam}(R_1\cup S_1))^{-4}
\frac{|x_2|^\delta }{|t_2|^{1+\delta}}\, dt_2dx_2
$$
which, by using the restriction on the variables, we bound by
$$
C \Big( \frac{|S_1|}{|R_1|}\Big)^{3/2}(|R_1|^{-1}{\rm diam}(R_1\cup S_1))^{-4}
$$
$$
|S_2|^{(1-\theta )\delta }\diam(R_2\cup S_2)^{\theta \delta }\diam (S_2\cup R_2)^{-(1+\delta )}
\| \psi_{R_2}\psi_{S_2}\|_{L^{1}(\mathbb R^2)}
$$
$$
\leq C \Big( \frac{|S_1|}{|R_1|}\Big)^{3/2}\prod_{i=1,2}(|R_1|^{-1}{\rm diam}(R_1\cup S_1))^{-4}
$$
$$
|S_2|^{(1-\theta )\delta }\diam(R_2\cup S_2)^{\theta \delta }\diam (S_2\cup R_2)^{-(1+\delta )}
|R_2|^{1/2}|S_2|^{1/2}
$$
$$
=C \Big( \frac{|S_1|}{|R_1|}\Big)^{3/2}(|R_1|^{-1}{\rm diam}(R_1\cup S_1))^{-4}
\Big(\frac{|S_2|}{|R_2|}\Big)^{(1/2+(1-\theta )\delta )}(|R_2|^{-1}{\rm diam}(R_2\cup S_2))^{-(1+(1-\theta )\delta )}
$$

\vskip10pt
{\bf e) WB-T(1).} We now consider
$\psi_{1,2}(t,x)=\psi(t,x) (1-\Phi_{\alpha \lambda_1 S_1}(x_1))\Phi_{\alpha \lambda_2 S_2}(x_2)\Phi_{\beta \mu_2 \tilde{R}_2}(t_2)$ which will be bounded by the use of weak
boundedness in $t_1,x_1$ and the special cancellation properties in $t_2,x_2$.

First, we impose the extra assumption $\psi_{1,2}((t_1,0),x)=0$ for any $t_1,x$.
The assumption is possible by the substitution of 
$$
\psi_{1,2}(t,x)-{\mathcal D}_{|R_2|}\Phi(t_2)\psi_{1,2}((t_1,0),x)
$$
and we first need to show that the subtracted term also satisfies the stated bounds.

We denote 
$
\tilde{\psi}_{1,2}(t_{1},x)=\psi_{1,2}((t_1,0),x)
$. 
Applying Lemma \ref{bump1} in the variable $x_1$ and similar reasoning as in previous case c) in the variable $t_{2}$, we have 
$$
|\tilde{\psi}_{1,2}(t_{1},x)|\leq C
\Big(\frac{|S_1|}{|R_1|}\Big)^{3/2}(|R_1|^{-1}\diam(R_1,S_1))^{-4}
$$
$$
|R_{2}|^{-1/2}(|R_2|^{-1}\diam(R_2,S_2))^{-N}
|\varphi_{R_1}(t_1)|  |\varphi_{S}(x)|
$$
and so ${\mathcal D}_{|R_2|}\Phi(t_2)\tilde{\psi}_{1,2}(t_1,x)$ is adapted to $(S_1\times R_2)\times (S_1\times S_2)$ with constant 
$$
|R_2|^{-1/2}(|R_2|^{-1}\diam(R_2,S_2))^{-N}\Big(\frac{|S_1|}{|R_1|}\Big)^{3/2}(|R_1|^{-1}\diam(R_1,S_1))^{-4}
$$

Moreover, since $\psi_{1,2}$ is a function of tensor product type, then we can 
write that $\tilde{\psi}_{1,2}(t_{1},x)=\phi_{S_1}(t_1)\psi_{S_1}(x_1)\psi_{S_2}(x_2)$ with each bump function adapted to the interval used as subindex. Notice that 
$\psi_{S_1}$ is $L^{2}$-adapted to the interval $S_{1}$ with constant 
$C(|S_1|/|R_1|)^{3/2}(|R_1|^{-1}\diam(R_1,S_1))^{-4}$
while $\psi_{S_2}$ is adapted to $S_{2}$ with constant 
$C|R_2|^{-1/2}(|R_2|^{-1}\diam(R_2,S_2))^{-N}$.

%
%
%

Now, we make use of the special cancellation condition $\Lambda (\phi_{S_1} \otimes 1,\varphi_{S_1}\otimes \psi_{S_2})=0$ for all
bump functions $\phi_{R_1},\varphi_{R_1}$ and all bump functions $\psi_{S_2}$ of mean zero.
Then, we can use the explicit estimate of Lemma \ref{definebmo2}
to get
$$
|\Lambda (\phi_{S_1}\otimes {\mathcal D}_{|R_{2}|}\Phi ,\psi_{S_1}\otimes \psi_{S_2})|
$$
$$
=|\Lambda(\phi_{S_1}\otimes {\mathcal D}_{\frac{|R_2|}{|S_2|}|S_2|}\Phi ,\psi_{S_1}\otimes \psi_{S_2})-\Lambda(\phi_{S_1}\otimes 1,\psi_{S_1}\otimes \psi_{S_2})|
$$
$$
\leq C
\Big(\frac{|S_1|}{|R_1|}\Big)^{3/2}
(|R_1|^{-1}\diam(R_1\cup S_1))^{-4}
\left(\frac{|R_2|}{|S_2|}\right)^{-\delta}\| \psi_{S_2} \|_1
$$
$$
\leq C\Big(\frac{|S_1|}{|R_1|}\Big)^{3/2}(|R_1|^{-1}\diam(R_1\cup S_1))^{-4}
\left(\frac{|S_2|}{|R_2|}\right)^{\delta}
|R_2|^{-1/2}(|R_2|^{-1}\diam(R_2\cup S_2))^{-N}|S_2|^{1/2}
$$
$$
=C\left(\frac{|S_1|}{|R_1|}\right)^{3/2}(|R_1|^{-1}\diam(R_1\cup S_1))^{-4}
\left(\frac{|S_2|}{|R_2|}\right)^{\delta+1/2}(|R_2|^{-1}\diam(R_2\cup S_2))^{-N}
$$
$$
\leq C\left(\frac{|S|}{|R|}\right)^{\delta+1/2}\prod_{i=1,2}(|R_i|^{-1}\diam(R_i\cup S_i))^{-4}
$$
which is no larger than the desired bound.

Now, with the assumption, we further decompose $\psi_{1,2}$ with respect the $t_{2}$ variable in 
similar way as we did before:
$$
\psi_{1,2}=\psi_{1,2,in}+\psi_{1,2,out}
$$
$$
\psi_{1,2,in}(t,x)= \psi_{1,2}(t,x)\Phi_{\frac{8}{32}\lambda S_{2}}(t_{2})
$$

e1) We first prove that $\psi_{1,2,in}$ is adapted to 
$(R_{1}\times \lambda_{2} S_{2})\times (R_{1}\times \lambda_{2} S_{2})$ with constant
$$
C\left(\frac{|S_{1}|}{|R_{1}|}\right)^{3/2}(|R_1|^{-1}{\rm diam}(R_1\cup S_1))^{-4}
\left(\frac{|S_{2}|}{|R_{2}|}\right)^{3/4}
(|R_2|^{-1}{\rm diam}(R_2\cup S_2))^{-7}
$$ 
and order $N_{2}=[N^{1/2}]$. 
Again, we prove this bound for the function, 
while the work for the derivatives follows in a similar way as in \cite{PPV}. 

The support of $\psi_{1,2,in}$ in the variable $t_{2}$ is in $\lambda_{2} S_{2}$. Then, for all  
$t_{2}\in \lambda_{2} S_{2}$ and all $x$, we have 
by the extra assumption, 
$$
|\psi_{1,2,in}(t,x)|\leq |\psi_{1,2}(t,x)|
=\Big| \int_{0}^{t_{2}} \partial_{t_{2}} \psi_{1,2}(t_{1},r,x)\, dr\Big| 
$$
$$
\le C|t_{2}|\| \partial_{t_{2}}\psi_{1,2}(t_{1,}\cdot ,x)\|_{\infty }
\leq C\lambda_{2} |S_{2}|\| \partial_{t_{2}}\psi_{1,2}(t_{1,}\cdot ,x)\|_{\infty }
$$
As before in (\ref{decaypsi0}), we can show that 
$$
\| \partial_{t_2} \psi_{1,2}(t_{1,}\cdot ,x)\|_{\infty }\leq C
|R_{2}|^{-3/2}(|R_2|^{-1}{\rm diam}(R_2\cup S_2))^{-N}
\varphi (t_{1,}x)
$$
where
$\varphi $ is an $L^2$-normalized bump function adapted to $R_{1}\times (S_{1}\times S_{2})$
because $\psi (\cdot ,t_{2},\cdot )$ is 
adapted to $R_{1}\times S$. 

Now, by the factor $1-\Phi_{\lambda_1 S_1}$ of $\psi_{1,2}$ and Lemma \ref{bump1} in the variable $x_{1}$,
we have that $\varphi $ is adapted to $R_1\times (R_1\times \lambda_{2}S_2)$ with constant
$$
C \Big(\frac{|S_1|}{|R_1|}\Big)^{3/2}\left(|R_1|^{-1}{\rm diam}(R_1\cup S_1)\right)^{-4}
$$

Also, by lemma \ref{bump4} in the variable $x_{2}$, we have that $\varphi $ is a bump function
adapted to $R_{1}\times (S_{1}\times \lambda_{2}S_{2})$ with constant $C\lambda_{2}^{N_{2}+1/2}$. 
and so,  
$$
|\varphi (t_{1},x)|
\leq C\lambda_{2}^{N_{2}+1/2}
|\lambda_{2}S_{2}|^{-1/2}(1+|\lambda_{2}S_{2}|^{-1}|x-c(S_{2})|)^{-N_{2}}
$$

This  
implies that  we can bound in the following way
$$
|\psi_{1,2}(t,x)|
\leq C \lambda_{2}|S_{2}|\left(\frac{|S_{1}|}{|R_{1}|}\right)^{3/2}
\left(|R_1|^{-1}{\rm diam}(R_1\cup S_1)\right)^{-4}
$$
$$
|R_{2}|^{-3/2}(|R_2|^{-1}{\rm diam}(R_2\cup S_2))^{-N}
\lambda_{2}^{N_{2}+1/2}|\lambda_{2}S_{2}|^{-1/2}
(1+|\lambda_{2}S_{2}|^{-1}|x-c(S_{2})|)^{-N_{2}}
$$
Therefore, 
$$
|\psi_{1,2,in}(t,x)|= |\psi_{1,2}(t,x)| |\Phi_{\frac{8}{32}\lambda S_{2}}(t)|
$$
$$
\leq C\left(\frac{|S_{1}|}{|R_{1}|}\right)^{3/2}
\left(|R_1|^{-1}{\rm diam}(R_1\cup S_1)\right)^{-4}
\left(\frac{|S_{2}|}{|R_{2}|}\right)^{3/2}(|R_2|^{-1}{\rm diam}(R_2\cup S_2))^{-N}\lambda_{2}^{N_{2}+2}
$$
$$
|\lambda_{2}S_{2}|^{-1/2}(1+|\lambda_{2}S_{2}|^{-1}|x-c(S_{2})|)^{-N_{2}}
\cdot |\lambda S_{2}|^{-1/2}| |\Phi_{\frac{8}{32}\lambda_{2} S_{2}}(t)|
$$
This already shows that $\psi_{1,2,in}$ is adapted to $(R_{1}\times \lambda_{2} S_{2})\times (R_{1}\times \lambda_{2} S_{2})$ and the constant can 
be simplified in a similar manner as we did in the case c.1).
 
%

Then, by the weak boundedness property of $\Lambda $, we get
$$
|\Lambda(\psi_{1,2,in})|
\leq C \left(\frac{|S_{1}|}{|R_{1}|}\right)^{3/2}\left(|R_1|^{-1}{\rm diam}(R_1\cup S_1)\right)^{-4}
\left(\frac{|S_{2}|}{|R_{2}|}\right)^{3/4}(|R_2|^{-1}{\rm diam}(R_2\cup S_2))^{-7}
$$

e2) We now work with $\psi_{1,2,out}$. As in case c.2), we only need to work in the case when 
$\frac{8}{32}\lambda_{2}|S_{2}|\leq \frac{1}{2}\mu_{2}|R_{2}|$.
By the extra assumption again, we have  
$$
|\psi_{1,2, out}(t,x)|\leq |\psi_{1,2}(t,x)|
=\Big| \int_{0}^{t_{2}} \partial_{t_{2}}\psi_{1,2}(t_{1},r,x)\, dr\Big| 
$$
\begin{equation}\label{decayofpsiplus}
\le C|t_{2}|
\|  \partial_{t_{2}}\psi_{1,2}(t_{1},\cdot ,x)\|_{\infty }
\leq C|t_{2}||R_{2}|^{-3/2}(|R_2|^{-1}{\rm diam}(R_2\cup S_2))^{-N}
\varphi (t_{1},x)
\end{equation}
obtained with similar calculations as before and where $\varphi $ is a bump function adapted to 
$R_{1}\times (S_{1}\times S_{2})$.

Furthermore, 
because of the factor $1-\Phi_{\lambda_1 S_1}$ of $\psi_{1,2}$, we can apply Lemma \ref{bump1} in the variable $x_1$,
to deduce that $\varphi $ is adapted to $R_1\times (R_1\times S_2)$ with a constant
of 
$C (|S_1|/|R_1|)^{3/2}\left(|R_1|^{-1}{\rm diam}(R_1\cup S_1)\right)^{-4}$.

Putting everything together, we can then assume the representation
$$
\psi_{1,2, out}(t,x)=\phi_{R_1}(t_1)\phi_{R_2}(t_2)\varphi_{R_1}(x_1)\psi_{S_2}(x_2)
$$ 
with   
$\varphi_{R_1}$ adapted to $R_{1}$ with constant  
$C (|S_1|/|R_1|)^{3/2}\left(|R_1|^{-1}{\rm diam}(R_1\cup S_1)\right)^{-4}$,
$$
|\phi_{R_2}(t_2)|
\leq C |t_{2}||R_{2}|^{-3/2}(|R_2|^{-1}{\rm diam}(R_2\cup S_2))^{-N}
(1-D_{\lambda_{2}|S_2|}\Phi (t_{2}))
$$
and $\psi_{S_2}$ bump function with mean zero adapted to 
$S_{2}$ with constant $C$.

On the support of $\psi_{1,2}$, we have $|t_{2}|\leq \frac{1}{2}\mu_{2}|\tilde{R_{2}}|=\frac{1}{2}\mu_{2}|R_{2}|$.
Moreover, on the support of
$\psi_{1,2,out}$ we have
that $|t_2|>\frac{8}{32}\lambda_{2}|S_{2}|$, while $|x_2|<\frac{2}{32}\lambda_{2}|S_2|$. 
This implies $2|x_2|<|t_2|$ and so by
the integral representation of the restricted operator $T^2_{t_2,x_2}$, we have
$$
\Lambda(\psi_{1,2})
=\int \phi_{R_2}(t_2)\psi_{S_2}(x_2) \langle T^2_{t_2,x_2}(\phi_{R_1}),\varphi_{R_1}\rangle \, dt_2dx_2
$$
Using the mean zero of $\psi_{S_2}$ we obtain for the above integral
$$
\int \phi_{R_2}(t_2)\psi_{S_2}(x_2) \langle (T^2_{t_2,x_2}-T^2_{t_2,0})(\phi_{R_1}),\varphi_{R_1}\rangle \, dt_2dx_2
$$
Since $2|x_2|<|t_2|$, by the mixed WB-CZ property of $\Lambda $ and the decays calculated, 
we can estimate the integral by
\begin{equation}\label{mixedwbcz}
\int |\psi_{R_2}(t_2)||\psi_{S_2}(x_2)| C\Big( \frac{|S_1|}{|R_1|}\Big)^{3/2}(|R_1|^{-1}{\rm diam}(R_1\cup S_1))^{-4}
\frac{|x_2|^\delta }{|t_2|^{1+\delta}}\, dt_2dx_2
\end{equation}
Now, using previous estimates for $\phi_{R_2}(t_2)\psi_{S_2}(x_2)$, we bound (\ref{mixedwbcz}) by
$$
C \Big( \frac{|S_1|}{|R_1|}\Big)^{3/2}(|R_1|^{-1}{\rm diam}(R_1\cup S_1))^{-4}
$$
$$
|R_{2}|^{-3/2}(|R_2|^{-1}{\rm diam}(R_2\cup S_2))^{-N}
\int_{\tiny \begin{array}{c}|x_{2}|<\frac{2}{32}\lambda_{2}|S_{2}|\\ \frac{8}{32}|\lambda_{2}S_{2}|<|t_{2}|\leq \frac{1}{2}\mu_{2}|R_{2}|\end{array}} 
|t_{2}|
\lambda_{2}^{-1/2}|\psi_{S_{2}}(x_{2})|\frac{|x_{2}|^\delta}{|t_{2}|^{1+\delta}}\, dt_{2}dx_{2}
$$
$$
\leq C \Big( \frac{|S_1|}{|R_1|}\Big)^{3/2}(|R_1|^{-1}{\rm diam}(R_1\cup S_1))^{-4}
\lambda_{2}^{1/2} |R_{2}|^{-3/2}(|R_2|^{-1}{\rm diam}(R_2\cup S_2))^{-N}
$$
$$
|S_{2}|^{-1/2}\int_{|x_{2}|<\frac{2}{32}\lambda_{2}|S_{2}|} |x_{2}|^\delta dx_{2} 
\int_{\frac{8}{32}|\lambda_{2}S_{2}|<|t_{2}|\leq \frac{1}{2}\mu_{2}|R_{2}|} \frac{1}{|t_{2}|^{\delta}}dt_{2}
$$
$$
\leq C \Big( \frac{|S_1|}{|R_1|}\Big)^{3/2}(|R_1|^{-1}{\rm diam}(R_1\cup S_1))^{-4}
|R_{2}|^{-3/2}(|R_2|^{-1}{\rm diam}(R_2\cup S_2))^{-N}
$$
$$
|S_{2}|^{-1/2}(\frac{2}{32}\lambda_{2}|S_{2}|)^{\delta +1}
((\frac{1}{2}\mu_{2}|R_{2}|)^{1-\delta}-(\frac{8}{32}\lambda_{2}|S_{2}|)^{1-\delta})
$$
$$
\leq C \Big( \frac{|S_1|}{|R_1|}\Big)^{3/2}(|R_1|^{-1}{\rm diam}(R_1\cup S_1))^{-4}
C \Big(\frac{|S_{2}|}{|R_{2}|}\Big)^{1/2+\delta'} (|R_2|^{-1}{\rm diam}(R_2\cup S_2))^{-7} 
$$
with similar calculations as before.

\vskip10pt
{\bf f) T(1)-CZ.} The last term is
$
\psi_{2,3}(t,x)=\psi(t,x)\Phi_{\alpha \lambda S}(x)\Phi_{\beta \mu_1 \tilde{R}_1}(t_1)(1-\Phi_{\beta \mu_2 \tilde{R}_2}(t_2))
$
for which we will use the special cancellation and the kernel representation.


As before, we first need to justify the the extra assumption $\psi_{2,3}((0,t_2),x)=0$ for any $t_2\in \mathbb R$, $x\in \mathbb R^2$. In order to do this, we divide the study into two cases: when $2|R_{2}|<\diam(R_{2}\cup S_{2})$ and when 
$\diam(R_{2}\cup S_{2})\leq 2|R_{2}|$.

In the first case, the assumption
comes from the substitution of $\psi_{2,3}(t,x)$ by 
$$
\psi_{2,3}(t,x)-{\mathcal D}_{|R_1|}\Phi(t_1)\psi_{2,3}((0,t_2),x)
$$
and we again prove that the subtracted term satisfies the stated bounds.

We denote $\tilde{\psi}_{2,3}(t_{2},x)=\psi_{2,3}((0,t_2),x)$. Applying 
similar reasoning as in case c) in the variable $t_{1}$ 
we have 
$$
|\tilde{\psi}_{2,3}(t_{2},x)|\leq C
|R_{1}|^{-1/2}(|R_1|^{-1}\diam(R_1\cup S_1))^{-N}
|\varphi_{R_1}(t_2)|  |\varphi_{S}(x)|
$$

Moreover, since $\psi_{2,3}$ is a function of tensor product type, we are allow  
write that $\tilde{\psi}_{2,3}(t_{2},x)=\phi_{R_2}(t_2)\psi_{S_1}(x_1)\psi_{S_2}(x_2)$ with each bump function adapted to the interval used as subindex. Notice that 
$\psi_{S_1}$ has mean zero and it is adapted to $S_{1}$ with constant 
$C|R_1|^{-1/2}(|R_1|^{-1}\diam(R_1\cup S_1))^{-N}$.

Now, due to the support of $\psi_{2,3}$ in the variables $t_{2}$ and $x_{2}$, we have
that $|t_2|>\mu_{2}|R_{2}|/4$ while $|x_2|<\frac{2}{32}\lambda_{2}|S_{2}|$. 
This implies as before $2|x_2|<|t_2|$ and so, since $|R_{2}|<\diam(R_{2}\cup S_{2})$, we can use Lemma 
 \ref{definebmo3}. 
Actually, by the special cancellation condition $\Lambda (1\otimes \phi_{S_2},\psi_{S_1}\otimes \psi_{S_2})=0$,
the explicit estimate of Lemma \ref{definebmo3}
and the decay of $\psi_{S_{2}}$, we have 
$$
|\Lambda ({\mathcal D}_{|R_1|}\Phi\otimes \phi_{R_2}\otimes ,\psi_{S_1}\otimes \psi_{S_2})|
$$
$$
=|\Lambda({\mathcal D}_{\frac{|R_1|}{|S_1|}|S_1|}\Phi \otimes \phi_{R_2},\psi_{S_1}\otimes \psi_{S_2})-\Lambda(1\otimes \phi_{R_2},\psi_{S_1}\otimes \psi_{S_2})|
$$
$$
\leq C\left(\frac{|S_2|}{|R_2|}\right)^{1/2+\delta }\Big(|R_{2}|^{-1}\diam(R_2\cup S_2)\Big)^{-(1+\delta )}
\left(\frac{|R_1|}{|S_1|}\right)^{-\delta}\| \psi_{S_1}\|_1
$$
$$
\leq C\left(\frac{|S_2|}{|R_2|}\right)^{1/2+\delta }\Big(|R_{2}|^{-1}\diam(R_2\cup S_2)\Big)^{-(1+\delta )}
\left(\frac{|S_1|}{|R_1|}\right)^{\delta}
|R_1|^{-1/2}(|R_1|^{-1}\diam(R_1\cup S_1))^{-N}|S_1|^{1/2}
$$
$$
=C\Big(\frac{|S_2|}{|R_2|}\Big)^{1/2+\delta }(|R_2|^{-1}{\rm diam}(R_2\cup S_2))^{-(1+\delta )}
\left(\frac{|S_1|}{|R_1|}\right)^{1/2+\delta}(|R_1|^{-1}\diam(R_1\cup S_1))^{-N}
$$
which is no larger than the desired bound.

On the other hand, in the case $\diam(R_{2}\cup S_{2})\leq 2|R_{2}|$, 
the extra assumption comes the substitution of $\psi_{2,3}(t,x)$ by 
$$
\psi_{2,3}(t,x)-{\mathcal D}_{\tau|R_1|}\Phi(t_1)\psi_{2,3}((0,t_2),x)
$$
with $\tau=(|R_{2}|/|S_{2}|)^{\frac{N+1+\delta}{\delta }}$ 
and we again need to prove that the subtracted term satisfies the stated bounds.

As before,
$\tilde{\psi}_{2,3}(t_{2},x)=\phi_{R_2}(t_2)\psi_{S_1}(x_1)\psi_{S_2}(x_2)$ with each bump function adapted to the interval used as subindex and such that 
$\psi_{S_1}$ has mean zero and it is adapted to $S_{1}$ with constant 
$C|R_1|^{-1/2}(|R_1|^{-1}\diam(R_1\cup S_1))^{-N}$.
But now, $\varphi_{R_{2}}$ and 
$\psi_{S_{2}}$ can be considered to be adapted to the same interval $R_{2}$. Obvioulsy, $\varphi_{R_{2}}$ will 
be adapted to $R_{2}$ with constant $C$. However, due to the difference on size, by Lemma \ref{bump4}, 
$\psi_{S_{2}}$ will be adapted to $R_{2}$ with constant $C(|R_{2}|/|S_{2}|)^{N+1/2}$. 

Then, by the special cancellation condition $\Lambda (1\otimes \phi_{S_2},\psi_{S_1}\otimes \psi_{S_2})=0$,
and the decay of $\psi_{S_{2}}$ and
the explicit estimate this time given by Lemma \ref{definebmo2},
we have 
$$
|\Lambda ({\mathcal D}_{\tau |R_1|}\Phi\otimes \phi_{R_2}\otimes ,\psi_{S_1}\otimes \psi_{R_2})|
$$
$$
=|\Lambda({\mathcal D}_{\tau \frac{|R_1|}{|S_1|}|S_1|}\Phi \otimes \phi_{R_2},\psi_{S_1}\otimes \psi_{R_2})-\Lambda(1\otimes \phi_{R_2},\psi_{S_1}\otimes \psi_{R_2})|
$$
$$
\leq C
\left(\frac{|R_2|}{|S_2|}\right)^{N+1/2}
\tau^{-\delta }\left(\frac{|R_1|}{|S_1|}\right)^{-\delta}\| \psi_{S_1}\|_1
$$
$$
\leq C\left(\frac{|R_2|}{|S_2|}\right)^{N+1/2}
\left(\frac{|R_2|}{|S_2|}\right)^{-\delta \frac{N+1+\delta }{\delta }}\left(\frac{|S_1|}{|R_1|}\right)^{\delta}
|R_1|^{-1/2}(|R_1|^{-1}\diam(R_1\cup S_1))^{-N}|S_1|^{1/2}
$$
$$
\leq C\Big(\frac{|S_2|}{|R_2|}\Big)^{1/2+\delta }(|R_2|^{-1}{\rm diam}(R_2\cup S_2))^{-(1+\delta )}
\left(\frac{|S_1|}{|R_1|}\right)^{1/2+\delta}(|R_1|^{-1}\diam(R_1\cup S_1))^{-N}
$$
which is again no larger than the desired bound.

\vskip10pt
Now, with the assumption, we further decompose $\psi_{2,3}$ with respect the variable $t_{1}$ in the following way:
$$
\psi_{2,3}=\psi_{2,3,in}+\psi_{2,3,out}
$$
$$
\psi_{2,3,in}(t,x)= \psi_{2,3}(t,x)\Phi_{\frac{8}{32}\lambda_{1} S_{1}}(t_{1})
$$

f1) We prove first that $\psi_{2,3,in}$, when considered in the variables $t_{1}$, $x_{1}$, is adapted to $\lambda_{1}S_{1}\times \lambda_{1}S_{1}$ with constant
$$
C\left(\frac{|S_{1}|}{|R_{1}|}\right)^{3/4}(|R_1|^{-1}{\rm diam}(R_1\cup S_1))^{-7}
$$ 
and order $N_{2}=[N^{1/2}]$. 
Again, we prove this bound for the function, 
while the work for the derivatives follows in a similar way. 

The support of $\psi_{2,3,in}$ in the variable $t_{1}$ is in $\lambda_{1} S_{1}$. Then, for all  
$t_{1}\in \lambda_{1} S_{1}$ and all 
$t_{2}, x$, we have 
by the extra assumption, 
$$
|\psi_{2,3}(t,x)|
=\Big| \int_{0}^{t_{1}} \partial_{t_{1}}\psi_{2,3}(r,t_{2},x)\, dr\Big| 
$$
$$
\le C|t_{1}|D_{\lambda_{1}|S_1|}\Phi (t_{1})\| \partial_{t_{1}}\psi_{2,3}(\cdot , t_{2},x)\|_{\infty }
\leq C\lambda_{1}|S_{1}|
\| \partial_{t_{1}}\psi_{2,3}(\cdot ,t_{2}, x)\|_{\infty }
$$
$$
\leq C\lambda_{1}|S_{1}||R_{1}|^{-3/2}(|R_1|^{-1}{\rm diam}(R_1\cup S_1))^{-N}
\varphi (t_{2}, x)
$$
with similar calculation as we did in c.1) and 
where $\varphi $ is an $L^2$-normalized bump function adapted to $R_{2}\times S$.

Moreover, by lemma \ref{bump4} in the variable $x_{1}$, we have that 
$\varphi$ is 
adapted to $R_{1}\times (\lambda_{1} S_{1}\times S_{2})$ with constant $C\lambda_{1}^{N_{2}+1/2}$. 
Therefore, 
$$
|\psi_{2,3,in}(t,x)|=|\psi_{2,3}(t,x)\Phi_{\frac{8}{32}\lambda_{1} S_{1}}(t_{1})|
$$
$$
\leq C\lambda_{1}|S_{1}||R_{1}|^{-3/2}(|R_1|^{-1}{\rm diam}(R_1\cup S_1))^{-N}
\lambda_{1}^{N_{2}+1/2}
\Phi_{\frac{8}{32}\lambda_{1} S_{1}}(t_{1})\varphi_{S_{2}} (t_{2})\psi_{\lambda_{1}S_{1}\times S_{2}}(x_{1})
$$
$$
\leq C\Big(\frac{|S_{1}|}{|R_{1}|}\Big)^{3/2}(|R_1|^{-1}{\rm diam}(R_1\cup S_1))^{-N}
\lambda_{1}^{N_{2}+2}
$$
$$
|\lambda_{1}S_{1}|^{-1/2}\Phi_{\frac{8}{32}\lambda_{1} S_{1}}(t_{1})\varphi_{S_{2}} (t_{2})\psi_{\lambda_{1}S_{1}\times S_{2}}(x_{1})
$$
Moreover, similar calculation as the one carried out in c.1), shows that 
$$
\Big(\frac{|S_{1}|}{|R_{1}|}\Big)^{3/2}(|R_1|^{-1}{\rm diam}(R_1\cup S_1))^{-N}
\lambda_{1}^{N_{2}+2}
\leq C\left(\frac{|S_{1}|}{|R_{1}|}\right)^{3/4}(|R_1|^{-1}{\rm diam}(R_1\cup S_1))^{-7}
$$
which is the bound we were searching for. 

Moreover, $\psi_{2,3,in}$ is of tensor product type, and then, 
we can write that 
$\psi_{2,3,in}(t,x)=\phi_{\lambda_{1}S_1}(t_1)\phi_{R_2}(t_2)\psi_{\lambda_{1}S_1}(x_1)\psi_{S_2}(x_2)$ 
where    
$\psi_{S_2}$ is a bump function with mean zero, 
$
\phi_{S_1}
$
is a bump function adapted to $\lambda_{1}S_{1}$ with constant 
\begin{equation}\label{constantadapted}
C\left(\frac{|S_{1}|}{|R_{1}|}\right)^{3/4}(|R_1|^{-1}{\rm diam}(R_1\cup S_1))^{-7}
\end{equation}
and $\psi_{\lambda_{1}S_1}$ has mean zero and it is adapted to 
$\lambda_{1}S_{1}$ with constant $C$.

On the other hand, in the support of $\psi_{2,3}$ we have
that $|t_2|>\diam (S_2\cup R_2)/4$ while $|x_2|<\frac{2}{32}|S_2|^{1-\theta }\diam(R_2\cup S_2)^{\theta }$. This implies $2|x_2|<|t_2|$ and so by
the integral representation of the restricted operator $T^2_{t_2,x_2}$ we have
$$
\Lambda(\psi_{2,3})
=\int \phi_{R_2}(t_2)\psi_{S_2}(x_2) \langle T^2_{t_2,x_2}(\phi_{\lambda_{1}S_1}),\psi_{\lambda_{1}S_1}\rangle \, dt_2dx_2
$$
$$
=\int \phi_{R_2}(t_2)\psi_{S_2}(x_2) \langle (T^2_{t_2,x_2}-T^2_{t_2,0})(\phi_{\lambda_{1}S_1}),\psi_{\lambda_{1}S_1}\rangle \, dt_2dx_2
$$
due to the mean zero of $\psi_{S_2}$.
Since $2|x_2|<|t_2|$, we can apply the mixed WB-CZ property and so, with the bound calculated
in (\ref{constantadapted}), we can bound this by
$$
\int |\psi_{R_2}(t_2)||\psi_{S_2}(x_2)| C\left(\frac{|S_{1}|}{|R_{1}|}\right)^{3/4}(|R_1|^{-1}{\rm diam}(R_1\cup S_1))^{-7}
\frac{|x_2|^\delta }{|t_2|^{1+\delta}}\, dt_2dx_2
$$
Finally, by similar work as in as we did in the case b), we can estimate the integral and obtain the bound
$$
C\Big(\frac{|S_1|}{|R_1|}\Big)^{3/4}
(|R_1|^{-1}{\rm diam}(R_1\cup S_1))^{-7}
\Big(\frac{|S_2|}{|R_2|}\Big)^{1/2+\delta' }(|R_2|^{-1}{\rm diam}(R_2\cup S_2))^{-(1+\delta' )}
$$

%

f2) We now work with $\psi_{2,3,out}$. Similarly as in previous cases c.2) and e.2), we can assume that  
$\frac{8}{32}\lambda_{1}|S_{1}|\leq \frac{1}{2}\mu_{1}|R_{1}|$.
By the extra assumption again, we have  
$$
|\psi_{2,3, out}(t,x)|\leq |\psi_{2,3}(t,x)|
=\Big| \int_{0}^{t_{1}} \partial_{t_{1}}\psi_{1,2}(r, t_{2}, x)\, dr\Big| 
$$
$$
\le C|t_{1}|
\|  \partial_{t_{1}}\psi_{2,3}(\cdot , t_{2}, x)\|_{\infty }
$$
$$
\leq C|t_{1}||R_{1}|^{-3/2}(|R_1|^{-1}{\rm diam}(R_1\cup S_1))^{-N}
\psi_{S_{1}}(x_{1}) \psi_{R_{2}\times S_{2}}(t_{2},x_{2})
$$
where, similar as before, $\psi_{S_{1}}$ is adapted to $S_{2}$
and $\psi_{R_{2}\times S_{2}}$ is adapted to 
$R_{2}\times S_{2}$. 

On one side we have the following:
due to the support of $\psi_{2,3}$, then $|t_{1}|\leq \mu_{1}|R_{1}|/2$; and due to the support of
$\psi_{2,3,out}$, we have that $|x_{1}|\leq \frac{2}{32}\lambda_{1}|S_{1}|$ while 
$|t_{1}|>\frac{8}{32}\lambda_{1}|S_{1}|$.

On the other hand, due to the support of $\psi_{2,3}$, we have  $|t_{2}|\geq \mu_{2}|R_{2}|/4=\diam(R_{2}\cup S_{2})/4$ and $|x_{i}|\leq \frac{2}{32}\lambda_{2}|S_{2}\leq \frac{2}{32}|S_{2}|^{1-\theta }\diam (R_{2}\cup S_{2})^{\theta }$.

Therefore, in both cases $2|x_{i}|<|t_{i}|$
and so,
we can use the Calder\'on Zygmund kernel representation
$$
\int \psi_{2,3,out}(t,x) K(x,t)\, dtdx
$$

Using the mean zero of $\psi_{2, 3, out}$ in the variables $x_{i}$ we can rewrite the above integral as
$$
\int \psi_{2, 3, out}(t,x) (K(x,t)-K((x_{1},0),t)-K((0,x_{2}),t)+K(0,t))\, dtdx
$$

Since $2|x_{i}|<|t_{i}|$, we can also use the smooth property of a product Calder\'on-Zygmund kernel 
and the constant as bump function of $\psi_{2, 3, out}$ and $\phi_{S_{1}}$, $\phi_{R_{2}\times S_{2}}$ 
to estimate the last display by
$$
C \int |t_{1}| |R_{1}|^{-3/2}(|R_1|^{-1}{\rm diam}(R_1\cup S_1))^{-N}
|\psi_{S_{1}}(x_{1})| | \psi_{R_{2}\times S_{2}}(t_{2},x_{2})| 
\frac{|x_{1}|^\delta}{|t_{1}|^{1+\delta}}\frac{|x_{2}|^\delta}{|t_{2}|^{1+\delta}}\, dtdx
$$
Now, 
using $|t_{2}|>\mu_{2}|R_{2}|/4$, $\mu_{1}|R_{1}|/2>|t_{1}|>\frac{8}{32}\lambda_{1}|S_{1}|$ and 
$|x_{i}|<\frac{2}{32}\lambda_{i}|S_{i}|$, we can estimate
the last display by
$$
C|R_{1}|^{-3/2}(|R_1|^{-1}{\rm diam}(R_1\cup S_1))^{-N}|S_{1}|^{-1/2}
\int_{|x_{1}|<\frac{2}{32}\lambda_{1}|S_{1}|} |x_{1}|^\delta dx_{1} 
\int_{\frac{8}{32}|\lambda_{1}S_{1}|<|t_{1}|\leq \frac{1}{2}\mu_{1}|R_{1}|} \frac{1}{|t_{1}|^{\delta}}dt_{1}
$$
$$
|S_2|^{(1-\theta)\delta }{\rm diam}(R_2\cup S_2)^{\theta \delta} 
{\rm \diam}(R_2\cup S_2)^{-1-\delta}
\|\psi_{R_{2}\times S_{2}}\|_{L^{1}(\mathbb R^{2})}
$$
$$
\le C|R_{1}|^{-3/2}(|R_1|^{-1}{\rm diam}(R_1\cup S_1))^{-N}|S_{1}|^{-1/2}
(\lambda_{1}|S_{1}|)^{\delta +1}((\mu_{1}|R_{1}|)^{1-\delta }-(\lambda_{1}|S_{1}|)^{1-\delta})
$$
$$
|S_2|^{(1-\theta)\delta } {\rm diam}(R_2\cup S_2)^{-1-\delta+\theta\delta}|R_2|^{1/2}|S_2|^{1/2}
$$
$$
\leq C \Big(\frac{|S_{1}|}{|R_{1}|}\Big)^{1/2+\delta'} (|R_1|^{-1}{\rm diam}(R_1\cup S_1))^{-7} 
\Big(\frac{|S_{2}|}{|R_{2}|}\Big)^{1/2+\delta' }
(|R_2|^{-1}\diam(R_2\cup S_2))^{-(1+\delta') }
$$
with similar calculations as in the cases b) and c.1). 

\vskip 20pt
We end the proof of Proposition \ref{twobumplemma} by dealing with the error terms.
We only check that the factors $c^5_1(t,x_2)=c_1(t,x_2)\Phi_{\alpha \lambda_1 S_1}(x_1)(1-\Phi_{\alpha \lambda_2 S_2}(x_2))$
and $c^5_3(t,x_1)=c_3(t,x_1)\Phi_{\alpha \lambda_2 S_2}(x_2)$ are small enough being analogous all the other ones.
In a similar way we obtained equality (\ref{c2}), we now have
$$
c_3(t,x_1)=-c|S_2|^{-1}\int \psi_{out}(t,x)(1-\Phi_{\alpha \lambda_2 S_2}(x_2))dx_2
$$
and so $\int c_3(t,x_1)dx_1=0$. Since
also $\int \psi_{1,1}(t,x)dx_1=0$ we have
$$
c_1(t,x_2)(1-\Phi_{\alpha \lambda_2 S_2}(x_2))\int \Phi_{\alpha \lambda_1 S_1}(x_1)dx_1
=-(1-\Phi_{\alpha \lambda_2 S_2}(x_2))\int \psi(t,x) (1-\Phi_{\alpha \lambda_1 S_1}(x_1))dx_1
$$
and thus
$$
|c_1(t,x_2)|\lambda_1 |S_1|
\leq \int_{|x_1|>\alpha \lambda_1|S_1|} |\psi (t,x)|dx_1
$$
$$
\leq C|S_1|^{-1/2}\int_{|x_1|>\alpha \lambda_1|S_1|} (1+|S_1|^{-1}|x_1|)^{-N}dx_1\, \, \phi (t,x_2)
$$
$$
\leq C{\lambda_1}^{-N}|S_1|^{1/2}\, \phi (t,x_2)
$$
that is,
$$
|c_1(t,x_2)|
\leq C{\lambda_1}^{-N}|S_1|^{-1/2}\, \phi (t,x_2)
$$
where $\phi $ is a $L^2$-normalized bump function adapted to $R\times S_2$.
With this we have that
$$
|c_1(t,x_2)\Phi_{\alpha \lambda_1 S_1}(x_1)|\leq C\lambda_1^{-N} \phi (t,x_2)|S_1|^{-1/2}\Phi_{\alpha \lambda_1 S_1}(x_1)
$$
with $|S_1|^{-1/2}\Phi_{\alpha \lambda_1 S_1}$ a bump function $L^2$-adapted to $S_1$. By lemma \ref{bump1}, this shows that
$|c_1(t,x_2)\Phi_{\alpha \lambda_1 S_1}(x_1)(1-\Phi_{\alpha \lambda_2 S_2}(x_2))|$ is adapted to $R\times R$ with a gain of constant of
$$
C\lambda_1^{-N}
\Big(\frac{|S_2|}{|R_2|}\Big)^{\theta N}(|R_2|^{-1}{\rm diam}(R_2\cup S_2))^{-\theta N}
$$
and by the definition of $\lambda_i$ this equals
$$
C(|S_1|^{-1}{\rm diam}(R_1\cup S_1))^{-\theta N}
\Big(\frac{|S_2|}{|R_2|}\Big)^{\theta N}(|R_2|^{-1}{\rm diam}(R_2\cup S_2))^{-\theta N}
$$
$$
=C\Big( \frac{|S|}{|R|}\Big)^{\theta N}\prod_{i=1,2}(|R_i|^{-1}{\rm diam}(R_i\cup S_i))^{-\theta N}
$$
which is smaller than the required bounds.

Symmetrically we have that since
also $\int \psi_{1,1}(t,x)dx_2=0$,
$$
c_3(t,x_1)\int \Phi_{\alpha \lambda_2 S_2}(x_2)dx_2
=-(1-\Phi_{\alpha \lambda_1 S_1}(x_1))\int \psi(t,x) (1-\Phi_{\alpha \lambda_2 S_2}(x_2))dx_2
$$
$$
-\Phi_{\alpha \lambda_1 S_1}(x_1)\int c_1(t,x_2)(1-\Phi_{\alpha \lambda_2 S_2}(x_2))dx_2
$$
The first term in the right hand side  can be treated in a similar way we did before and obtain
$$
|(1-\Phi_{\alpha \lambda_1 S_1}(x_1))\int \psi(t,x) (1-\Phi_{\alpha \lambda_2 S_2}(x_2))dx_2|
\leq C{\lambda_2}^{-N}|S_2|^{1/2}\, (1-\Phi_{\alpha \lambda_1 S_1}(x_1))\phi (t,x_1)
$$
For the second term, we use the definition of $c_1(t,x_2)=-c|S_1|^{-1}\int \psi(t,x)(1-\Phi_{\alpha \lambda_1 S_1}(x_1))dx_1$
to
bound by
$$
\Big| \frac{\Phi_{\alpha \lambda_1 S_1}(x_1)}{\lambda_1|S_1|}\int \int \psi(t,x)(1-\Phi_{\alpha \lambda_1 S_1}(x_1))dx_1(1-\Phi_{\alpha \lambda_2 S_2}(x_2))dx_2\Big|
$$
$$
\leq \frac{\Phi_{\alpha \lambda_1 S_1}(x_1)}{\lambda_1|S_1|}\int_{|x_i|>\alpha \lambda_i|S_i|} |\psi (t,x)|dx_1dx_2
$$
$$
\leq C\frac{\Phi_{\alpha \lambda_1 S_1}(x_1)}{\lambda_1|S_1|}|S_1|^{-1/2}|S_2|^{-1/2}\int_{|x_i|>\alpha \lambda_i|S_i|} (1+|S_2|^{-1}|x_2|)^{-N}(1+|S_2|^{-1}|x_2|)^{-N}dx_1dx_2\, \, \phi (t)
$$
$$
\leq C\frac{\Phi_{\alpha \lambda_1 S_1}(x_1)}{\lambda_1|S_1|}{\lambda_1}^{-N}{\lambda_2}^{-N}|S_1|^{1/2}|S_2|^{1/2}\, \phi (t)
$$
where $\phi $ is a $L^2$-normalized bump function adapted to $R$.

Both things together imply
$$
|c_3(t,x_1)|\lambda_2 |S_2|
\leq C{\lambda_2}^{-N}|S_2|^{1/2}\, (1-\Phi_{\alpha \lambda_1 S_1}(x_1)) \phi (t,x_1)
+C\frac{\Phi_{\alpha \lambda_1 S_1}(x_1)}{\lambda_1|S_1|}{\lambda_1}^{-N}{\lambda_2}^{-N}|S_1|^{1/2}|S_2|^{1/2}\, \phi (t)
$$
and so
$$
|c_3(t,x_1)|
\leq C{\lambda_2}^{-N}|S_2|^{-1/2}\, \phi (t,x_1)
+C{\lambda_1}^{-N}{\lambda_2}^{-N}|S_1|^{-1/2}|S_2|^{-1/2}\, \Phi_{\alpha \lambda_1 S_1}(x_1)\phi (t)
$$

With this we have that
$$
|c_3(t,x_1)\Phi_{\alpha \lambda_2 S_2}(x_2)|\leq C\lambda_2^{-N} \phi (t,x_1)(1-\Phi_{\alpha \lambda_1 S_1}(x_1))|S_2|^{-1/2}\Phi_{\alpha \lambda_2 S_2}(x_2)
$$
$$
+C{\lambda_1}^{-N}{\lambda_2}^{-N}|S_1|^{-1/2}\Phi_{\alpha \lambda_1 S_1}(x_1) |S_2|^{-1/2}\Phi_{\alpha \lambda_2 S_2}(x_2)\phi (t)
$$
with $|S_i|^{-1/2}\Phi_{\alpha \lambda_i S_i}$ a bump function $L^2$-adapted to $S_i$. This shows that the function is adapted to $R\times R$ with constant
$$
C\lambda_2^{-N}\Big(\frac{|S_1|}{|R_1|}\Big)^{\theta N}(|R_1|^{-1}{\rm diam}(R_1\cup S_1))^{-\theta N}
+C\lambda_1^{-N} \lambda_2^{-N}
$$
where we have used Lemma \ref{bump1} for the first term, which is analogous to the previous case. The second one is also all right since by definition of
$\lambda_i$ we have
$$
C\lambda_1^{-N} \lambda_2^{-N}=C\prod_{i=1,2}(|S_i|^{-1}{\rm diam}(R_i\cup S_i))^{-\theta N}
$$
$$
=C\Big( \frac{|S|}{|R|}\Big)^{\theta N}\prod_{i=1,2}(|R_i|^{-1}{\rm diam}(R_i\cup S_i))^{-\theta N}
$$
$$
\hskip300pt \boxempty
$$

\section{Proof of the main result}\label{L2}

\begin{theorem}{($L^2$ boundedness). }\label{L2bounds}
Let
$\Lambda$ be a bilinear Calder\'on-Zygmund form
satisfying the mixed WB-CZ condition.

We also assume that $\Lambda $ satisfies
the weak boundedness condition,
and the special cancellation conditions
$$
\Lambda (1,\psi_{1}\otimes \psi_{2})=\Lambda (\psi_{1}\otimes \psi_{2},1)=\Lambda (\psi_{1}\otimes 1,1\otimes \psi_{2})=\Lambda (1\otimes \psi_{1},\psi_{2}\otimes 1)=0
$$
for all functions $\psi_{i}\in {\mathcal S}(\mathbb R)$ with mean zero and
$$
\Lambda (f_{1} \otimes 1,g_{1}\otimes \psi )=\Lambda (1\otimes f_{2},\psi \otimes g_{2})
=\Lambda (\psi \otimes f_{2},1\otimes g_{2})=\Lambda (f_{1}\otimes \psi ,g_{1}\otimes 1)=0
$$
for all functions $f_{i},g_{i}\in {\mathcal S}(\mathbb R)$ and all bump functions 
$\psi\in {\mathcal S}(\mathbb R)$ with  mean zero.

Then $\Lambda_0, \Lambda_1, \Lambda_2$ are bounded bilinear forms on $L^2$.
\end{theorem}
\proof
Because of the symmetry on the hypothesis it is clear that we only need to prove the result for $\Lambda $.
We decompose the frequency plane in the standard way to obtain first a Littlewood-Paley decompositon and later
a wavelet decomposition.

Let $\phi \in \S (\mathbb R)$ be an even function such that
$\widehat{\phi }$ is supported in $\{ \xi \in \mathbb R: |\xi |\leq 2\} $ and equals $1$ on
$\{ \xi \in \mathbb R: |\xi |\leq 1\} $. Let $\psi $ be the function
$
\psi (x)=\phi (x)-1/2\phi (x/2)
$.
Then $\hat{\psi }$ is supported on the annulus $\{ \xi \in \mathbb R: 2^{-1}\leq |\xi |\leq 2\} $ and moreover
$
\sum_{k\in \mathbb Z}\widehat{\psi }(\xi /2^k)\approx 1
$,
for all $\xi \neq 0$. We define the Littlewood-Paley projection operators in $\mathbb R$ given by
$
P_{k}(f)=f*D^{1}_{2^{-k}}\psi
$
and
$
P_{\leq k}(f)=f*D^{1}_{2^{-k}}\phi
$.
We observe that
$
\lim_{k\rightarrow \infty }P_{\leq k}(f)=f
$
while
$
\lim_{k\rightarrow \infty }P_{\leq -k}(f)=0
$
where in both cases the convergence is understood in the topology of $\S (\mathbb R)$.

We consider now their counterparts in the biparameter case: for $k\in \mathbb Z^2$,
$$
P_{k}(f)=f*(D^{1}_{2^{-k_1}}\psi \otimes D^{1}_{2^{-k_2}}\psi )
\hskip 30pt
P_{\leq k_1,k_2}(f)=f*(D^{1}_{2^{-k_1}}\phi \otimes D^{1}_{2^{-k_2}}\psi )
$$
$$
P_{k_1,\leq k_2}(f)=f*(D^{1}_{2^{-k_1}}\psi \otimes D^{1}_{2^{-k_2}}\phi )
\hskip 30pt
P_{\leq k}(f)=f*(D^{1}_{2^{-k_1}}\phi \otimes D^{1}_{2^{-k_2}}\phi )
$$
which satisfy
$
\lim_{k\rightarrow \infty }P_{\leq k}(f)=f
$
in the topology of $\S (\mathbb R^2)$ while the other three operators tend to zero in the same sense.

For $N\in \mathbb N$, let $\Lambda_N$ be the bilinear form given by
$$
\Lambda_{N}(f,g)=\sum_{|k_i|,|j_i|\leq N}\Lambda (P_{j}f, P_{k}g)
$$
where $k,j\in \mathbb Z^2$.
We see that for all $f,g\in \S (\mathbb R^2)$ we have
$
\Lambda (f,g)=\lim_{N\rightarrow \infty }\Lambda_{N}(f,g)
$:
unfolding the sum in $\Lambda_{N}$, we have
\begin{eqnarray*}
\Lambda_{N}(f,g)&=&\Lambda (P_{\leq (N,N)}f,P_{\leq (N,N)}g)-\Lambda (P_{\leq (-(N-1),-(N-1))}f,P_{\leq (N,N)}g)\\
&&-\Lambda (P_{\leq (N,N)}f,P_{\leq (-(N-1),-(N-1))}g)+\Lambda (P_{\leq -(N-1,N-1)}f,P_{\leq -(N-1,N-1)}g)
\end{eqnarray*}
and by the continuity of $\Lambda $ we have that the first term tends to $\Lambda (f,g)$ while the other three tend to zero.

Let us now consider the family of intervals
$\omega_{k_i}=[-2^{k_i+1},-2^{k_i-1}]\cup [2^{k_i-1},2^{k_i+1}]$. Since $P_k(f)$ has Fourier support in
$\omega_{k}=\omega_{k_1}\times \omega_{k_2}$,
we have by Shannon's sampling theorem that
$$
P_k(f)=\sum_{R}\langle f,\psi_{R,\, \omega_k}\rangle \psi_{R,\, \omega_k}
$$
where the sum runs over all dyadic rectangles $R=R_1\times R_2$ such that $|R_i|=|\omega_{k_i}|^{-1}$
and the convergence is understood in the topology of ${\mathcal S}(\mathbb R^2)$.
Moreover, the functions $\psi_{R,\, \omega_k}$ satisfy that
$\psi_{R,\, \omega_k}=\psi_{R_1,\, \omega_{k_1}}\otimes \psi_{R_2,\, \omega_{k_2}}$ where $\psi_{R_i,\, \omega_{k_i}}$ are Schwartz functions such that
$\supp \widehat{\psi_{R_i,\, \omega_{k_i}}}\subset \omega_{k_i}$
and $e^{-2\pi ic(\omega_{k_i})}\psi_{R_i,\, \omega_{k_i}}$ are bump functions adapted to $R_{i}$. From now we drop the index $\omega_k$ in the notation of
$\psi_{R}$.

Then by continuity of $\Lambda $ in $\S (\mathbb R^2)$, we can write
$$
\Lambda (f,g)=\sum_{k,j}\Lambda (P_{j}f, P_{k}g)
=\sum_{R,S}\langle f,\psi_{R}\rangle \langle g,\psi_{S}\rangle \Lambda (\psi_R,\psi_S)
$$
where now the sums run over the whole family of dyadic rectangles in $\mathbb R^2$.
From now we work to obtain bounds of the last expression when the
sum runs over finite families of dyadic rectangles in such way that the bounds are independent of the  particular families of rectangles.
Because of the rate of decay of Corollary \ref{symmetricspecialcancellation},
we parametrize the sums according to eccentricities and relative positions of the rectangles:
$$
\sum_{R,S}\langle f,\psi_{R}\rangle \langle g,\psi_{S}\rangle \Lambda (\psi_R,\psi_S)
=\sum_{i=1,2}\sum_{e_i\in \mathbb Z}\sum_{m_i\in \mathbb N}\sum_{R}\sum_{\tiny \begin{array}{c}S\in R_{e,m}\end{array}}
\langle f,\psi_{R}\rangle \langle g,\psi_{S}\rangle \Lambda (\psi_R,\psi_S)
$$
where for fixed eccentricities $e_i$, relative distances $m_i$ and every given rectangle $R$,
we define the family
$$
R_{e,m}=\{ S:|R_i|=2^{e_i}|S_i|, m_i\leq \max(|R_i|,|S_i|)^{-1}\diam(R_i\cup S_i)< m_i+1 \, \, \mbox{ for } i=1,2\}
$$
Notice that by symmetry the product family $\{ (R,S): S\in R_{e,m}\}$ can also be parameterized as $\{ (R,S): R\in S_{-e,m}\}$ with analogous definition for
$S_{-e,m}$.

We denote by $\sum_{P}$ the three first sums over parameters.
By Lemma \ref{twobumplemma} and Cauchy's inequality, we bound the previous quantity by
$$
\sum_{P}\sum_{R}\sum_{\tiny \begin{array}{c}S\in R_{e,m}\end{array}}
|\Lambda (\psi_R,\psi_S)|
|\langle f,\psi_{R}\rangle ||\langle g,\psi_{S}\rangle |
$$
$$
\lesssim \sum_{P}\sum_{R}\sum_{\tiny \begin{array}{c}S\in R_{e,m}\end{array}}2^{-(|e_1|+|e_2|)(1/2+\delta )}(m_1m_2)^{-(1+\delta )}
|\langle f,\psi_{R}\rangle | |\langle g,\psi_{S}\rangle |
$$
$$
\leq \sum_{P}2^{-|e_1+e_2|(1/2+\delta )}(m_1m_2)^{-(1+\delta )}
\Big( \sum_{R}\hspace{-.3cm} \sum_{\tiny \begin{array}{c}S\! \in \! R_{e,m}\end{array}}|\langle f,\psi_{R}\rangle |^2\Big)^{1/2}
\Big( \sum_{S}\hspace{-.3cm} \sum_{\tiny \begin{array}{c}R\! \in \! S_{-e,m}\end{array}}|\langle g,\psi_{S}\rangle |^2\Big)^{1/2}
$$

Now, for every fixed $R_i$ and each $m_i\in \mathbb N$ there are $2^{\max(e_i,0)}$ dyadic intervals $S_i$ such that
$|R_i|=2^{e_i}|S_i|$ and  $m_i\leq \max(|R_i|,|S_i|)^{-1}\diam(R_i\cup S_i)<m_i+1$. This implies that the cardinality of
$R_{e_i,m_i}$ is $2^{\max(e_1,0)}2^{\max(e_2,0)}$.
For the same reason, the cardinality of $S_{-e,m}$ is $2^{\max(-e_1,0)}2^{\max(-e_2,0)}=2^{-\min(e_1,0)}2^{-\min(e_2,0)}$.
Then, the previous expression coincides with
\begin{eqnarray*}
\sum_{P}
&2^{-|e_1+e_2|(1/2+\delta )}(m_1m_2)^{-(1+\delta )}
&\Big( 2^{\max(e_1,0)}2^{\max(e_2,0)}\sum_{R}|\langle f,\psi_{R}\rangle |^2\Big)^{1/2}\\
&&\Big( 2^{-\min(e_1,0)}2^{-\min(e_2,0)}\sum_{\tiny \begin{array}{c}S\end{array}}|\langle g,\psi_{S}\rangle |^2\Big)^{1/2}
\end{eqnarray*}
$$
\leq \prod_{i=1,2}\sum_{e_i\in \mathbb Z}2^{-|e_i|(1/2+\delta )}2^{\max(e_i,0)/2}2^{-\min(e_i,0)/2}\sum_{m_i\in \mathbb N}m_i^{-(1+\delta )}
\| f\|_2\| g\|_2
$$
$$
=\Big( \sum_{e\in \mathbb Z}2^{-|e|\delta }\sum_{m\in \mathbb N}m^{-(1+\delta )}\Big)^{2}\| f\|_2\| g\|_2
$$
since $2^{\max(e_i,0)}2^{-\min(e_i,0)}=2^{|e_i|}$.

\section{Extension to $L^{p}$ spaces}\label{Lp}

As said in the introduction, the weak $L^1$ estimates are no longer true
in the multi-parameter case. So, in order to prove
$L^p$ bounds we cannot apply to our operator the classical method of interpolating between $L^2$ and the weak $L^1$ estimates.
Instead, we follow the steps of
the previous proof and perform again a decomposition of the dual pair which will be controlled
by multi-parameter square functions whose $L^p$ boundedness follows from weak $L^1$ bounds in the uni-parameter case.

\begin{definition} Given a $L^2(\mathbb R^2)$-normalized basis $(\psi_{R})_{R}$, we define the double square function by
$$
SS(f)=\Big( \sum_{R}\frac{|\langle f,\psi_{R}\rangle |^2}{|R|}\chi_{R}\Big)^{1/2}
$$
where the sum runs over all dyadic rectangles $R$ in $\mathbb R^2$.
\end{definition}

See \cite{CF} and specially \cite{FS} for a proof of boundedness of $SS$ on $L^{p}(\mathbb R^2)$ with $1<p<\infty $.

\vskip15pt
We also need to consider the following modified double square function
\begin{definition}
Let $k\in \mathbb Z^2$, $n\in \mathbb N^2$. 
We define
$$
SS_{k,n}(f)(x)=\big( \sum_{R}\sum_{S\in R_{k,n}}\frac{\langle f,\psi_{R}\rangle ^2}{|S|}\chi_{S}(x)\Big)^{1/2}
$$
where for every dyadic rectangle 
$R$, $R_{k,n}$ is the set of 
dyadic rectangles $S$ such that $|R_i|=2^{k_i}|S_i|$
and $n_i\leq \frac{\diam(R_i\cup S_i)}{\max(|R_i|,|S_i|)}<n_i+1$.
\end{definition}

We notice that the double square function corresponds to the values $k_i=0$, $n_i=1$.

We also notice that whenever $k_{i}\geq 0$ we have $|S_{i}|\leq |R_{i}|$ and so, 
\begin{equation}\label{relationsquares}
SS_{k,n}(f)(x)=\big( \sum_{R}\sum_{S\in R_{k,n}}\frac{\langle f,\psi_{R}\rangle ^2}{|S|}\chi_{S}(x)\Big)^{1/2}
=\big( \sum_{R}\langle f,\psi_{R}\rangle ^2
\sum_{S\in R_{k,n}}\frac{\chi_{S}(x)}{|S|}\Big)^{1/2}
\end{equation}
$$
=\big( \sum_{R}\langle f,\psi_{R}\rangle ^2\frac{\chi_{\tilde{S}}(x)}{2^{-k}|\tilde{S}|}\Big)^{1/2}
=2^{(k_1+k_2)/2}SS_{0,n}(f)(x)
$$
where
$\tilde{S}$ is the dyadic rectangle such that $|\tilde{S}_i|=|R_i|$ and
$|R_i|^{-1}\diam(\tilde{S}_i\cup R_i)=n_i$.

We state in the proposition below boundedness of this modified square function. Its proof follows directly from
the analogous result in the uni-parameter case and so, for the sake of completeness, at the end of the paper
we include an appendix in which a proof of this result in the uni-parameter case can be found
(see Proposition \ref{modifiedsquarefunction}).

\begin{proposition}\label{modifieddoublesquarefunction}For every $1<p<\infty$,
$$
\| SS_{k,n}(f)\|_{L^{p}(\mathbb R^2)}
\leq C_{p}\prod_{i=1,2} (2^{-k_{i}\, \sign(\frac{2}{p}-1)}\log(n_i)+1)^{|\frac{2}{p}-1|}\| f\|_{L^{p}(\mathbb R^2)}
$$
\end{proposition}
\proof
Given $k\in \mathbb Z^2$ and $n\in \mathbb N^2$, let
$TT_{k,n}$ be the operator defined by
$$
TT_{k,n}(f)(x)
=\sum_{R}\sum_{S\in R_{k,n}}\langle f,\psi_{R}\rangle \psi_{S}(x)
$$
$$
=\sum_{S}\sum_{R\in S_{-k,n}}\langle f,\psi_{R}\rangle \psi_{S}(x)
$$
where the relationship between $R$ and $S$ is the same one given in the definition of the modified square function.
Now we see how the double square function of $TT_{k,n}(f)$ relates with $SS_{k,n}(f)$. On one side, 
$$
SS(TT_{k,n}(f))(x)
=\Big( \sum_{S}\Big( \sum_{R\in S_{-k,n}}\langle f,\psi_{R}\rangle \Big)^2\frac{\chi_{S}(x)}{|S|}\Big)^{1/2}
$$
$$
\geq \Big( \sum_{S}\sum_{R\in S_{-k,n}}\langle f,\psi_{R}\rangle ^2\frac{\chi_{S}(x)}{|S|}\Big)^{1/2}
=SS_{k,n}(f)(x)
$$

On the other side, because the cardinality of $S_{-k_{i},n_{i}}$ is  $2^{-\min (k_i,0)}$, 
we have
$$
SS(TT_{k,n}(f))(x)
=\Big( \sum_{S}\Big( \sum_{R\in S_{-k,n}}\langle f,\psi_{R}\rangle \Big)^2\frac{\chi_{S}(x)}{|S|}\Big)^{1/2}
$$
$$
\leq \Big( \sum_{S}2^{-(\min (k_1,0)+\min (k_2,0))}\sum_{R\in S_{-k,n}}\langle f,\psi_{R}\rangle ^2\frac{\chi_{S}(x)}{|S|}\Big)^{1/2}
$$
$$
=2^{-(\min (k_1,0)+\min (k_2,0))/2}SS_{k,n}(f)(x)
$$

Moreover, by linearity
$$
TT_{k,n}(f)(x)=\sum_{S}\langle f,\psi_{R}\rangle \psi_{S}(x)
$$
$$
=\sum_{S_1}\sum_{R_{1}\in S_{-k_{1},n_{1}}}\langle \sum_{S_2}\sum_{R_{2}\in S_{-k_{2},n_{2}}}\langle f,\psi_{R_2}\rangle \psi_{S_2}(x_2),\psi_{R_1}\rangle \psi_{S_1}(x_1)
$$
$$
=T_{k_1,n_1}(T_{k_2,n_2}(f)(\cdot,x_2))(x_1)
$$
where $T_{k_i,n_i}(g)$ is defined in the obvious way
$$
T_{k_i,n_i}(g)(x_i)=\sum_{S_i}\sum_{R_{i}\in S_{-k_{i},n_{i}}}\langle g,\psi_{R_i}\rangle \psi_{S_i}(x_i)
$$
while $T_{k_2,n_2}(f)(x_1,x_2)=T_{k_2,n_2}(f_{x_1})(x_2)$
and $f_{y_1}(y_2)=f(y_1,y_2)$.

By the first inequality above, we have 
$$
\| SS_{k,n}(f)\|_{L^{p}(\mathbb R^2)}\leq
\| SS(TT_{k,n}(f))\|_{L^{p}(\mathbb R^2)}
$$
$$
\approx \|TT_{k,n}(f)\|_{L^{p}(\mathbb R^2)}
$$
and we just need to bound the last norm. Now, 
Corollary of Proposition \ref{modifiedsquarefunction} gives us the boundedness result in the uni-parameter case. If
we denote by $C_{p,k,n}=2^{-\min (k,0)/2}2^{k/2}(\log(n)+1)^{|\frac{2}{p}-1|}$ for $k\geq 0$ and 
$C_{p,k,n}=2^{-\min (k,0)/2}2^{-k/2}(2^{-k_i\sign(\frac{2}{p}-1)}+\log(n)+1)^{|\frac{2}{p}-1|}$ for $k\leq 0$, then we have
$$
\| T_{k_i,n_i}(f)\|_{L^{p}(\mathbb R)}\leq C_p \| S(T_{k_i,n_i}(f))\|_{L^{p}(\mathbb R)}
$$
$$
\leq C_{p}2^{-\min (k_i,0)/2}\| S_{k_i,n_i}(f)\|_{L^{p}(\mathbb R)}\lesssim C_{p,k_{i},n_{i}}\| f\|_{L^{p}(\mathbb R)}
$$
Therefore, we finally obtain
$$
\|TT_{k,n}(f)\|_{L^{p}(\mathbb R^2)}=\Big( \int_{\mathbb R} \| T_{k_1,n_1}(T_{k_2,n_2}(f)(\cdot ,x_2))\|_{L^{p}(\mathbb R)}^{p}dx_2\Big)^{1/p}
$$
$$
\lesssim \Big( \int_{\mathbb R} C_{p,k_{1},n_{1}}\| T_{k_2,n_2}(f)(\cdot ,x_2)\|_{L^{p}(\mathbb R)}^{p}dx_2\Big)^{1/p}
$$
$$
=C_{p,k_{1},n_{1}}\Big( \int_{\mathbb R} \| T_{k_2,n_2}(f)(x_1,\cdot )\|_{L^{p}(\mathbb R)}^{p}dx_1\Big)^{1/p}
$$
$$
=C_{p,k_{1},n_{1}} \Big( \int_{\mathbb R} \| T_{k_2,n_2}(f_{x_1})\|_{L^{p}(\mathbb R)}^{p}dx_1\Big)^{1/p}
$$
$$
\lesssim C_{p,k_{1},n_{1}} \Big( \int_{\mathbb R} C_{p,k_{2},n_{2}}\| f_{x_1}\|_{L^{p}(\mathbb R)}^{p}dx_1\Big)^{1/p}
$$
$$
=C_{p,k_{1},n_{1}}C_{p,k_{2},n_{2}}\| f\|_{L^{p}(\mathbb R^2)}
$$

\vskip10pt
Now we turn to the main result of this section.
\begin{theorem}{($L^p$ boundedness). }\label{Lpbounds}
Let
$\Lambda$ be a bilinear Calder\'on-Zygmund form
satisfying the mixed WB-CZ condition.

We also assume that $\Lambda $ satisfies
the weak boundedness condition,
and the special cancellation conditions
and the special cancellation conditions
$$
\Lambda (1,\psi_{1}\otimes \psi_{2})=\Lambda (\psi_{1}\otimes \psi_{2},1)=\Lambda (\psi_{1}\otimes 1,1\otimes \psi_{2})=\Lambda (1\otimes \psi_{1},\psi_{2}\otimes 1)=0
$$
for all functions $\psi_{i}\in {\mathcal S}(\mathbb R)$ with mean zero and
$$
\Lambda (f_{1} \otimes 1,g_{1}\otimes \psi )=\Lambda (1\otimes f_{2},\psi \otimes g_{2})
=\Lambda (\psi \otimes f_{2},1\otimes g_{2})=\Lambda (f_{1}\otimes \psi ,g_{1}\otimes 1)=0
$$
for all functions $f_{i},g_{i}\in {\mathcal S}(\mathbb R)$ and all bump functions 
$\psi\in {\mathcal S}(\mathbb R)$ with  mean zero.

Then $\Lambda_0, \Lambda_1, \Lambda_2$ are bounded bilinear forms on $L^p$.
\end{theorem}
\proof
Again, we prove the result only for $\Lambda $.
As in previous theorem, we use a $L^2$-normalized wavelet expansion of the functions appearing in the dual pair, we
parametrize the terms accordingly with eccentricity and relative distances in exactly the same way and we apply
Lemma \ref{symmetricspecialcancellation} to obtain
$$
|\Lambda (f,g)|\leq \sum_{R,S}|\Lambda (\psi_R,\psi_S)| |\langle f,\psi_R\rangle  | |\langle g,\psi_S \rangle  |
$$
$$
\leq C\sum_{i=1,2}\sum_{e_i\in \mathbb Z}\sum_{m_i\in \mathbb N}
2^{-|e_1+e_2|(1/2+\delta )}(m_1m_2)^{-(1+\delta )}
\sum_{R}\sum_{S\in R_{e,m}}
|\langle f,\psi_{R}\rangle  | |\langle g,\psi_{S} \rangle  |
$$
where $R_{e,m}$ is the set of dyadic rectangles $S$ such that
$|R_i|=2^{e_i}|S_i|$ and
$n_i\leq \frac{\diam(R_i\cup S_i)}{\max(|R_i|,|S_i|)}<n_i+1$.
Notice that $S\in R_{e,m}$ if and only if $R\in S_{-e,m}$.

Now,
we denote by $K=K_1\times K_2$ the rectangle minimum, that is, such that $K_i=R_i$ if $|R_i|\leq |S_i|$ and $K_i=S_i$ otherwise.
This way, the inner sum can be rewritten as
$$
\sum_{\tiny \begin{array}{c}(R,S)\in P_{e,m}\end{array}}\int_{\mathbb R^2}\frac{|\langle f,\psi_{R}\rangle  |}{|K|^{1/2}}
\frac{|\langle g,\psi_{S}\rangle |}{|K|^{1/2}}\chi_{K}(x)dx
$$
$$
\leq \int_{\mathbb R^2}
\Big( \sum_{R}\sum_{S\in R_{e,m}}
\frac{|\langle f,\psi_{R}\rangle  |^2}{|K|}\chi_{K}(x)\Big)^{1/2}
\Big( \sum_{S}\sum_{R\in S_{-e,m}} 
\frac{|\langle g,\psi_{S}\rangle |^2}{|K|}\chi_{K}(x)\Big)^{1/2}dx
$$

In order to build up the modified square functions,
we denote by $k,k'\in \mathbb Z^2$ the scale parameters $k_i=\max(e_i,0)$, $k'_i=-\min(e_i,0)$ and
by $n,n'\in \mathbb Z^2$ the translation parameters
$n_i=m_i$, $n_{i}'=1$ if $e_i\geq 0$ while $n_i=1$, $n_{i}'=m_i$ if $e_i\leq 0$.
Notice that $2^{\max(e_i,0)}2^{-\min(e_i,0)}=2^{|e_i|}$ and $n_in_{i}'=m_i$. Therefore, previous expression can 
be rewritten as 
\begin{equation}\label{4modsquare}
\int_{\mathbb R^2}
SS_{k,n}(f)(x)SS_{k',n'}( g)(x)dx
\end{equation}

We show how to bound the first factor.
By the choice of $K$ we have that $|R_i|,|S_i|\geq |K_i|$ which is why  
$k_i\geq 0$ and $k_{i}'\geq 0$. 
Then, as we saw before in (\ref{relationsquares})
$$
SS_{k,n}(f)(x)
=2^{(k_1+k_2)/2}SS_{0,n}(f)(x)
$$

This implies that expression (\ref{4modsquare}) is equal to
$$
\int_{\mathbb R^2}2^{\max(e_1,0)/2}2^{\max(e_2,0)/2}SS_{0,n}(f)(x)2^{-\min(e_1,0)/2}2^{-\min(e_2,0)/2}SS_{0,n'}( g)(x)dx
$$
$$
\leq 2^{(|e_1|+|e_2|)/2}\| SS_{0,n}(f)\|_{L^{p}(\mathbb R^2)} \| SS_{0,n'}(g)\|_{L^{p'}(\mathbb R^2)}
$$
According to the boundedness of
the modified
double square functions given by Proposition \ref{modifieddoublesquarefunction}, this can be bounded by
$$
C_p\, 2^{(|e_1|+|e_2|)/2} \prod_{i=1,2}(\log(n_i)+1)^{|\frac{2}{p}-1|}
(\log(n'_{i})+1)^{|\frac{2}{p}-1|}
\| f\|_{L^{p}(\mathbb R^2)}\| g\|_{L^{p'}(\mathbb R^2)}
$$
$$
= C_{p}\, 2^{(|e_1|+|e_2|)/2}\prod_{i=1,2}(\log(m_{i})+1)^{|\frac{2}{p}-1|}
\| f\|_{L^{p}(\mathbb R^2)}\| g\|_{L^{p'}(\mathbb R^2)}
$$
$$
\leq C_{p}\, 2^{(|e_1|+|e_2|)/2}\prod_{i=1,2}m_{i}^{\epsilon |\frac{2}{p}-1|}
\| f\|_{L^{p}(\mathbb R^2)}\| g\|_{L^{p'}(\mathbb R^2)}
$$

Then, putting everything back together, we have
$$
|\Lambda (f,g)|
\leq C_p \sum_{i=1,2}\sum_{e_i\in \mathbb Z}\sum_{m_i\in \mathbb N}2^{-(|e_1|+|e_2|)(1/2+\delta )}(m_1m_{2})^{-(1+\delta )}
2^{(|e_1|+|e_2|)/2} (m_1m_{2})^{\epsilon |\frac{2}{p}-1|}
$$
$$
=C_{p}\prod_{i=1,2}\sum_{e_i\in \mathbb Z}2^{-|e_i|\delta }\sum_{m_i\in \mathbb N}m_{i}^{-(1+\delta -\epsilon |\frac{2}{p}-1|)}
\| f\|_{L^{p}(\mathbb R^2)} \| g\|_{L^{p'}(\mathbb R^2)}
\leq C_{p}
\| f\|_{L^{p}(\mathbb R^2)} \| g\|_{L^{p'}(\mathbb R^2)}
$$
as long as  $\epsilon |\frac{2}{p}-1|<1+\delta $.

\section{The general case: different types of paraproducts}\label{para}

We devote this last section to the extension of the previous theorems to the general case, that is, the proof of boundedness for
singular integral operators that do not satisfy the special cancellation properties.
As in the classical case, this is done by constructing appropriate paraproducts.
But in the multiparametric case, the process is more involved not only because we need more paraproducts (three different types in 
total) but also because these paraproducts cannot be independent each other. 

In particular, let $b_i$ with $i=1,\ldots ,4$,
be four functions in $\BMO(\mathbb R^2)$ and
$b_{i}$ with $i=5,\ldots ,8$, be four functions in $\BMO(\mathbb R)$.
Let also $\Lambda $ be a bilinear form satisfying the hypotheses of Theorem \ref{biparT1} such that $\Lambda (1\otimes 1,\cdot )=b_1$
$\Lambda (\cdot ,1 \otimes 1)=b_2$, $\Lambda (\cdot \otimes 1, 1\otimes \cdot )=b_3$, 
$\Lambda (1\otimes \cdot ,\cdot \otimes 1)=b_4$
and so on.

In order to prove boundedness of $\Lambda $, we construct eight
bilinear forms $\Lambda_{i}$ organized in three different groups
in such a way that their associated linear operators are bounded and moreover they recover the functions $b_i$, in the sense that
for example
$\Lambda_{b_1}(1\otimes 1,\cdot  )=b_1$, while the bilinear form vanishes in all other possible cases, namely,
$\Lambda_{b_1}(\cdot ,1\otimes 1)
=\Lambda_{b_1}(\cdot \otimes 1 ,1\otimes \cdot )=0$ 
and so on. 

This way the bilinear form $\tilde{\Lambda}=\Lambda -\sum_{i}\Lambda_{b_i}$ satisfies the eight
special cancellation hypotheses
of Corollary \ref{symmetricspecialcancellation} and so, by applying the corollary, we
deduce that $\tilde{\Lambda}$ is bounded.
Moreover, since every $\Lambda_{b_{i}}$ is also bounded by construction, we finally obtain boundedness of the
initial form $\Lambda $.

Before we start with the construction of paraproducts, we present a lemma
that in some way shows that the sufficient conditions we have used in the main theorem are the right ones, while it
also justifies the paraproducts we will define later on. For the sake of simplicity, we write the proof only for operators that preserve the space support,
since then the error terms are zero and then the expression can be written by means of the Haar basis.

Let $(h_{I})_{I}$ the Haar basis in $\mathbb R$ defined by $h_{I}=|I|^{-1/2}\chi_{I_l}-|I|^{-1/2}\chi_{I_r}$ where $I_l$ and
$I_r$ are the children intervals of $I$. Let
$(h_{R})_{R}$ the Haar basis in $\mathbb R^2$ defined by $h_{R}=h_{R_1}\otimes h_{R_2}$.

\begin{lemma}\label{representation}
Let $T:{C}_{0}^{\infty }(\mathbb R^2)\rightarrow \mathbb C$ be a linear mapping continuous with respect to the topology in ${C}_{0}^{\infty }(\mathbb R^2)$,
such that $\supp T(f)\subset \supp f$.
Then
$$
\langle T(f),g\rangle
=\sum_{R}\langle f,h_{R}\rangle \langle g,h_{R}\rangle \langle T(h_{R}),h_{R}\rangle
$$
$$
+\big\langle \sum_{R}\langle f, h_{R}^2\rangle \langle g, h_{R}\rangle h_{R}\, , \, T(1)\big\rangle
+\big\langle \sum_{R}\langle f, h_{R}\rangle \langle g, h_{R}^2\rangle h_{R}\, , \,T^{*}(1)\big\rangle
$$
$$
+\big\langle \sum_{R}\langle f, h_{R_1}h_{R_2}^2\rangle \langle g, h_{R_1}^2h_{R_2}\rangle h_{R}\, , \,T_1(1) \big\rangle
+\big\langle \sum_{R}\langle f, h_{R_1}^2h_{R_2}\rangle \langle g, h_{R_1}h_{R_2}^2\rangle h_{R}\, , \,T_1^{*}(1) \big\rangle
$$
$$
+\sum_{R}\langle f, h_{R_1}h_{R_2}^2\rangle \langle g, h_{R_1}h_{R_2}\rangle \langle T(h_{R_1}\otimes 1),h_{R}\rangle
+\sum_{R}\langle f, h_{R_1}h_{R_2}\rangle \langle g, h_{R_1}h_{R_2}^2\rangle \langle T^{*}(h_{R_1}\otimes 1),h_{R}\rangle
$$
$$
+\sum_{R}\langle f, h_{R_1}^2h_{R_2}\rangle \langle g, h_{R_1}h_{R_2}\rangle \langle T(1\otimes h_{R_2}),h_{R}\rangle
+\sum_{R}\langle f, h_{R_1}h_{R_2}\rangle \langle g, h_{R_1}^2h_{R_2}\rangle \langle T^{*}(1\otimes h_{R_2}),h_{R}\rangle
$$
\end{lemma}
\begin{remark}
The formula for more general operators includes some error terms whose contribution is smaller than the one described in the previous statement.

Let $(\psi_{R})_{R}$ be a wavelet basis in $\mathbb R^2$ and for every rectangle $R$, let $\psi_R^2$ be a bump function $L^1$-adapted to $R$ and of mean one.
Then, such a general formula can be stated in the following way:
\begin{eqnarray*}
\langle T(f),g\rangle
&={\displaystyle \sum_{n\in \mathbb Z^2}}\Big( &\sum_{R}\langle f,\psi_{R}\rangle \langle g,\psi_{R^n}\rangle \langle T(\psi_{R}),\psi_{R^n}\rangle \\
&+&\big\langle \sum_{R}\langle f, \psi_{R}^2\rangle \langle g, h_{R^n}\rangle \psi_{R}\, , \, T(1)\big\rangle
+\ldots \\
&+&\sum_{R}\langle f, \psi_{R_1}\psi_{R_2}^2\rangle \langle g, \psi_{R_1^n}\psi_{R_2^n}\rangle \langle T(\psi_{R_1}\otimes 1),\psi_{R^n}\rangle
+\ldots \Big)
\end{eqnarray*}
where $R_i^n=R_i+n_i|R_i|$. The leading term
is associated with $n=0$, which is the one appearing in the statement of the lemma.
\end{remark}
\proof
Since $\supp T(h_{R})\subset \supp h_{R}=R$ we have
$$
\langle T(f),g\rangle = \sum_{R\cap S\neq \emptyset }f_Rg_S\langle T(h_R),h_{S}\rangle
$$
where $f_R=\langle f,h_R\rangle $ and the same for the function $g$.

Now, given two dyadic rectangles $R,S$ such that $R\cap S\neq \emptyset $
there are only nine different possibilities, namely,
\begin{enumerate}
\item[1)] $R=S$, which leads to $\langle T(h_R),h_R\rangle $
\item[2)] $R\subset S$, which gives $T^{*}(1)$
\item[3)] $S\subset R$, which analogously gives $T(1)$
\item[4)] $R<S$, meaning $R_1\subset S_1$ and $R_2\subset S_2$, which leads to $T_1(1)=T_2^{*}(1)$
\item[5)] $S<R$, meaning $R_1=S_1$ and $R_2\subset S_2$, which leads to $T_1^{*}(1)=T_2(1)$
\item[6)] $R_1=S_1$ and $S_2\subset R_2$, which leads to $\langle T(h_{R_1}\otimes 1),h_{R_1}\otimes h_{R_2}\rangle $
\item[7)] $R_1=S_1$ and $R_2\subset S_2$, which leads to $\langle T(h_{R_1}\otimes h_{R_2}),h_{R_1}\otimes 1\rangle $
\item[8)] $S_1\subset R_1$ and $R_2=S_2$, which leads to $\langle T(1\otimes h_{R_2}),h_{S_1}\otimes h_{R_2}\rangle $
\item[9)] $R_1\subset S_1$ and $R_2=S_2$, which leads to $\langle T(h_{R_1}\otimes h_{R_2}),1\otimes h_{R_2}\rangle $
\end{enumerate}

Then the decomposition of $\langle T(f),g\rangle $ is obtained as follows: from 1) we get directly the first term
$$
\sum_{R}f_Rg_R\langle T(h_R),h_R\rangle
$$
From 2) and 3) we get the two following terms (we only write the second one)
$$
\sum_{R}\sum_{R\subset S}f_Rg_S\langle T(h_R),h_{S}\rangle
=\sum_{R}f_R\langle T(h_R),\big( \sum_{R\subset S}g_S\, h_{S}\big) \chi_{R})\rangle
$$
$$
=\sum_{R}f_R\langle T(h_R),m_R(g)\rangle
=\sum_{R}f_R\, m_R(g)\langle T(h_R),1\rangle
=\sum_{R}f_R\, m_R(g)\langle h_R,T^{*}(1)\rangle
$$
$$
=\Big\langle \sum_{R}f_R\, m_R(g)h_R,T^{*}(1)\Big\rangle
$$
where we write $m_R(g)=\big( \sum_{R\subset S}g_S\, h_{S}\big) \chi_{R}=\langle g,h_R^2\rangle =|R|^{-1}\int_{R}g(x)dx$.

On the other hand, from 4) and 5) and using the partial adjoints $T_i$ we get the two following ones (we only write the fourth one)
$$
\sum_{R}\sum_{R<S}f_Rg_S\langle T(h_R),h_{S}\rangle
=\sum_{R}\sum_{R<S}f_Rg_S\langle T_2(h_{R_1}\otimes h_{S_2}),h_{S_1}\otimes f_{R_2}\rangle
$$
$$
=\sum_{R_1,S_2}\Big\langle T_2(h_{R_1}\otimes h_{S_2}),
\big( \sum_{R_1\subset S_1, S_2\subset R_2}f_Rg_S h_{S_1}\otimes h_{R_2}\big) \chi_{R_1\times S_2}\Big\rangle
$$
$$
=\sum_{R_1,S_2}\Big\langle T_2(h_{R_1}\otimes h_{S_2}),
\big(\sum_{R_1\subset S_1}(g_{S_2})_{S_1} h_{S_1}\big)
\otimes \big( \sum_{S_2\subset R_2}(f_{R_1})_{R_2}h_{R_2}\big) \chi_{R_1\times S_2}\Big\rangle
$$
$$
=\sum_{R_1,S_2}\langle T_2(h_{R_1}\otimes h_{S_2}),m_{R_1}(g_{S_2})m_{S_2}(f_{R_1}) \rangle
=\sum_{R_1,S_2}m_{R_1}(g_{S_2})m_{S_2}(f_{R_1})\langle T_2(h_{R_1}\otimes h_{S_2}),1 \rangle
$$
$$
=\langle \sum_{R_1,S_2}m_{S_2}(f_{R_1})m_{R_1}(g_{S_2})h_{R_1}\otimes h_{S_2},T_2^{*}(1) \rangle
=\langle \sum_{R}m_{R_2}(f_{R_1})m_{R_1}(g_{R_2})h_{R},T_1(1) \rangle
$$

From 6) 7), 8) and 9) we get the remaining terms (we only write the sixth one)
$$
\sum_{S}\sum_{R_1=S_1, S_2\subset R_2}f_Rg_S\langle T(h_R),h_{S}\rangle
=\sum_{S}\sum_{R_1=S_1, S_2\subset R_2}f_Rg_S\langle h_R,T^{*}(h_{S})\rangle
$$
$$
=\sum_{S}g_S\Big\langle \big( \sum_{S_2\subset R_2}f_{S_1\times R_2}h_{S_1\times R_2}\big)\chi_{S_2},T^{*}(h_{S})\Big\rangle
=\sum_{S}g_S\Big\langle \big( \sum_{S_2\subset R_2}\langle f,h_{S_1\times R_2}\rangle h_{S_1}\otimes h_{R_2}\big)\chi_{S_2}
,T^{*}(h_{S})\Big\rangle
$$
$$
=\sum_{S}g_S\Big\langle h_{S_1}\otimes \big(\sum_{S_2\subset R_2}\langle f,h_{S_1\times R_2}\rangle h_{R_2}\big)\chi_{S_2}
,T^{*}(h_{S})\Big\rangle
=\sum_{S}g_S\langle h_{S_1}\otimes m_{S_2}(f_{S_1}),T^{*}(h_{S})\rangle
$$
$$
=\sum_{S}m_{S_2}(f_{S_1})g_S\langle h_{S_1}\otimes 1,T^{*}(h_{S})\rangle
=\sum_{S}m_{S_2}(f_{S_1})g_S\langle T(h_{S_1}\otimes 1),h_{S}\rangle
$$

We start now with the construction of paraproducts. We need up to eight of such operators but by symmetry it will be enough to show only three of them. In particular,
we construct the paraproducts associated with $T(1)$, $T_1(1)$ and $\langle T(\psi_{R_1}\otimes 1),\psi_{S_1}\rangle $.
\begin{lemma}{(Classical paraproducts).}\label{paraproducts1}
Given a function $b$ in $\BMO(\mathbb R^2)$, there exists a bounded bilinear form $\Lambda_b^1$
such that 
\begin{enumerate}
\item $\Lambda_b^1(1\otimes 1,\psi)=b$,
$\Lambda_b^1(\psi, 1\otimes 1)=0$ for every $\psi \in {\mathcal S}(\mathbb R^{2})$ with mean zero
\item $\Lambda_b^1(\psi_{1} \otimes 1, 1\otimes \psi_{2} )
=\Lambda_b^1(1\otimes \psi_{1} ,\psi_{2}\otimes 1)=0$ for every $\psi_{i} \in {\mathcal S}(\mathbb R)$ with mean zero

\end{enumerate}
\end{lemma}
\proof

Let $(\psi_I)_I$ be a wavelets basis on $L^2(\mathbb R)$
Let $(\psi_R)_R$ be the wavelets basis on $L^2(\mathbb R^2)$ defined by $\psi_R=\psi_{R_1}\otimes \psi_{R_2}$.
We denote by $\varphi_I$ a bump function adapted to $I$ such that $\widehat{\varphi_{I}}$ is adapted to an interval  of measure comparable with
$|I|^{-1}$ and center the origin. Let finally $\varphi_R=\varphi_{R_1}\otimes \varphi_{R_2}$.

We define the bilinear form
$$
\Lambda_b^1(f,g)=\sum_R \langle b, \psi_{R}\rangle \langle f, \varphi_{R}\rangle \langle g,\psi_{R}\rangle
$$
which, to simplify notation, we will just denote by $\Lambda_b$ during the proof of the lemma.

At least formally, $\Lambda_b$ satisfies
$$
\Lambda_b(1,g)=\langle g,b\rangle
$$
$$
\Lambda_b(f,1)=\Lambda_b(f_1\otimes 1,1\otimes g_2)=\Lambda_b(1\otimes f_2,g_1\otimes 1)=0
$$
being the proof trivial in all these cases.

\vskip10pt
For the proof of their boundedness we proceed by using the duality $H^1(\mathbb R^2)-\BMO(\mathbb R^2)$. Since
$$
\Lambda_b(f,g)
=\sum_R \langle b, \psi_{R}\rangle \langle f, \varphi_{R}\rangle \langle g,\psi_{R}\rangle
$$
$$
=\Big\langle b, \sum_R \langle f, \varphi_{R}\rangle \langle g,\psi_{R}\rangle \psi_{R}\Big\rangle
$$
we have
$$
|\Lambda_b(f,g)|\leq \| b\|_{\BMO(\mathbb R^2)}
\| \sum_R \langle f, \varphi_{R}\rangle \langle g,\psi_{R}\rangle \psi_{R}\|_{H^1(\mathbb R^2)}
$$
Just assuming the sum is finite, we get that $\sum_R \langle f, \psi_{R}^2\rangle \langle g,\psi_{R}\rangle \psi_{R}\in H^1(\mathbb R^2)$ and then
$$
\| \sum_R \langle f, \varphi_{R}\rangle \langle g,\psi_{R}\rangle \psi_{R}\|_{H^1(\mathbb R^2)}
\approx  \| S(\sum_R \langle f, \varphi_{R}\rangle \langle g,\psi_{R}\rangle \psi_{R})\|_{L^1(\mathbb R^2)}
$$
with implicit constants independent of the number of terms in the sum. Now
$$
S\Big(\sum_R \langle f, \varphi_{R}\rangle \langle g,\psi_{R}\rangle \psi_{R}\Big)^2
$$
$$
=\sum_R \langle f, \varphi_{R_1}\otimes \varphi_{R_2}\rangle \langle g,\psi_{R_1}\otimes \psi_{R_2}\rangle^2
\frac{\chi_{R_1}}{|R_1|}\otimes \frac{\chi_{R_2}}{|R_2|}
$$
$$
\leq \sup_{R} |\langle f, \varphi_{R_1}\otimes \varphi_{R_2}\rangle |^2 \hskip5pt \sum_R\langle g,\psi_{R_1}\otimes \psi_{R_2}\rangle^2
\frac{\chi_{R_1}}{|R_1|}\otimes \frac{\chi_{R_2}}{|R_2|}
$$
$$
=(M\otimes M)(f)^2(S\otimes S)(g)^2
$$
where $M\otimes M$ and $S\otimes S$ are defined by the two previous expressions and are known to be bounded operators on $L^p(\mathbb R^2)$.
Then finally
$$
\| S(\sum_R \langle f, \varphi_{R}\rangle \langle g,\psi_{R}\rangle \psi_{R})\|_{L^1(\mathbb R^2)}
$$
$$
\leq \| (M\otimes M)(f)\|_{L^p(\mathbb R^2)}\| (S\otimes S)(g)\|_{L^{p'}(\mathbb R^2)}
\leq C\| f\|_{L^p(\mathbb R^2)}\| g\|_{L^{p'}(\mathbb R^2)}
$$

We still need to prove that this family of operators also belong to the class of operators for which the
theory applies.
In particular, we show that
they have integral representations like the ones stated in definition \ref{intrep} with
kernels satisfying the definition of a product Calder\'on-Zygmund kernel \ref{prodCZ}. We also prove
the operator satisfies the weak boundedness Calder\'on-Zygmund
condition stated in \ref{WB-CZ}. From
$$
\Lambda_b(f,g)
=\sum_{R} \langle b, \psi_{R}\rangle
\langle f,\varphi_{R} \rangle
\langle g,\psi_{R}\rangle
=\int f(t)g(x)
\sum_{R} \langle b, \psi_{R}\rangle
\varphi_{R}(t)\psi_{R}(x)dtdx
$$
we obtain the integral representation regardless disjointness of supports of the argument functions. Moreover,
$$
K(x,t)=\Big\langle b, \sum_{R}
\varphi_{R}(t)\psi_{R}(x)\psi_{R}\Big\rangle
$$
and we check the two properties of a product C-Z kernel:
$$
|K(x,t)|\leq \| b\|_{\BMO(\mathbb R^2)} \Big\| \sum_{R} \varphi_{R}(t)\psi_{R}(x)\psi_{R}\Big\|_{H^1(\mathbb R^2)}
$$
As before, the $H^1$-norm is equivalent to
$$
\Big\| S\Big(\sum_{R} \psi_{R}^2(t)\psi_{R}(x)\psi_{R}\Big)\Big\|_{L^1(\mathbb R^2)}
=\int \Big(\sum_{R} \varphi_{R}(t)^2\psi_{R}(x)^2\frac{\chi_{R}(y)}{|R|}\Big)^{1/2}dy
$$
\begin{equation}\label{kernelclass}
=\prod_{i=1,2}\int \Big(\sum_{R_i} \varphi_{R_i}(t_i)^2\psi_{R_i}(x_i)^2\frac{\chi_{R_i}(y_i)}{|R_i|}\Big)^{1/2}dy_i
\end{equation}

Now we are going to use the fact that 
$|\varphi_{R_i}|\leq |R_i|^{-1}\phi_{R_i}$ and $|\psi_{R_i}|\leq |R_i|^{-1/2}\phi_{R_i}$,
where $\phi_{R_i}$ is a bump function $L^{\infty }$-adapted to $R_{i}$ of order $N$. 

If $x_{i},t_{i}\in R_{i}$, we take 
$I_{x_i,t_i}$ to be the smallest dyadic interval such that $x_i,t_i\in I_{x_i,t_i}$ and let $(I_k)_{k\geq 0}$ the family of dyadic intervals such that
$I_{x_i,t_i}\subset I_k$ with $|I_k|=2^{k}|I_{x_i,t_i}|$.
Then, since $|\varphi_{R_i}|\lesssim |R_i|^{-1}$ and $|\psi_{R_i}|\lesssim |R_i|^{-1/2}$,
we can bound (\ref{kernelclass}) by
$$
\sum_{k\geq 0}\int_{I_{k+1}\backslash I_{k}}\Big(\sum_{j\geq k} \frac{1}{|I_j|^2}\frac{1}{|I_j|}\frac{1}{|I_j|}\Big)^{1/2}dy
=\sum_{k\geq 0}\Big(\sum_{j\geq k}\frac{1}{|I_j|^4}\Big)^{1/2}|I_{k+1}\backslash I_{k}|
$$
$$
=\sum_{k\geq 0}\Big(\sum_{j\geq k}\frac{1}{2^{4j}|I_{x_i,t_i}|^4}\Big)^{1/2}2^{k}|I_{x_i,t_i}|
\lesssim \frac{1}{|I_{x_i,t_i}|}\sum_{k\geq 0}\frac{1}{2^{2k}}2^{k}
\lesssim \frac{1}{|x_i-t_i|}
$$
On the other hand, if $x_{i}\notin R_{i}$, we decompose into dyadic pieces of the form 
$2^{j}|R_{i}|\leq |x_{i}-c(R_{i})|\leq 2^{j+1}|R_{i}|$ for $j\geq 1$ and use, in a similar way as before, that
$|\varphi_{R_i}|\leq |R_i|^{-1}2^{-jN}$. Actually, previous calculation shows that
in this case we can bound (\ref{kernelclass}) by 
$$
\sum_{j\geq 1}\sum_{k\geq 0}\int_{I_{k+1}\backslash I_{k}}\Big(\sum_{j\geq k} \frac{1}{|I_j|^2}\frac{2^{-2jN}}{|I_j|}\frac{1}{|I_j|}\Big)^{1/2}dy
\lesssim \frac{1}{|I_{x_i,t_i}|}\sum_{j\geq 1}2^{-jN}\sum_{k\geq 0}\frac{1}{2^{2k}}2^{k}
\lesssim \frac{1}{|x_i-t_i|}
$$
ending the first condition. 

For the second one, we prove that
$$|\partial_{t_1}\partial_{t_2}K(x,t)|
+|\partial_{t_1}\partial_{x_2}K(x,t)|+|\partial_{x_1}\partial_{t_2}K(x,t)|
+|\partial_{x_1}\partial_{x_2}K(x,t)|\leq C \prod_{i=1,2}|x_i-t_i|^{-2}$$
For expository reasons we deal only with the second term
$$
\partial_{t_1}\partial_{x_2}K(x,t)=
\Big\langle b, \sum_{R} \partial_{t_1}\varphi_{R_1}(t_1)\varphi_{R_2}(t_2)\psi_{R_1}(x_1)\partial_{x_2}\psi_{R_2}(x_2)\psi_{R}
\Big\rangle
$$
The four possible terms are not really symmetric since the averaging function $\psi_{R_i}^2$ only appear in
the $t_i$ variables. So, sometimes the derivatives hit an averaging function while some other times they do not.
However, it is the presence of derivatives of wavelets what produces the final estimates,
regardless whether it is
$\partial_{t_i}\varphi_{R_i}(t_i)$ or $\partial_{x_i}\psi_{R_i}(x_i)$. Actually, in all cases the derivatives
increase by one the
degree of the powers in $|R|$ involved and so, they have the same impact in all four terms. Let's see this point.
As before
$$
|\partial_{t_1}\partial_{x_2}K(x,t)|\leq \| b\|_{\BMO(\mathbb R^2)}
\Big\| \sum_{R} \partial_{t_1}\varphi_{R_1}(t_1)\varphi_{R_2}(t_2)\psi_{R_1}(x_1)\partial_{x_2}\psi_{R_2}(x_2)
\psi_{R}\Big\|_{H^1(\mathbb R^{2})}
$$
with the $H^1$-norm equivalent to
$$
\Big\| S\Big(\sum_{R} \partial_{t_1}\varphi_{R_1}(t_1)\varphi_{R_2}(t_2)\psi_{R_1}(x_1)\partial_{x_2}\psi_{R_2}(x_2)
\psi_{R}\Big)\Big\|_{L^1(\mathbb R^2)}
$$
$$
=\int_{\mathbb R^2} \Big(\sum_{R} (\partial_{t_1}\varphi_{R_1}(t_1))^2\varphi_{R_2}(t_2)^2
\psi_{R_1}(x_1)^2(\partial_{x_2}\psi_{R_2}(x_2))^2\frac{\chi_{R}(y)}{|R|}\Big)^{1/2}dy
$$
$$
=\int \Big(\sum_{R_1} (\partial_{t_1}\varphi_{R_1}(t_1))^2\psi_{R_1}(x_1)^2\frac{\chi_{R_1}(y_1)}{|R_1|}\Big)^{1/2}dy_1
\int \Big(\sum_{R_2} \varphi_{R_2}(t_2)^2(\partial_{x_2}\psi_{R_2}(x_2))^2\frac{\chi_{R_2}(y_2)}{|R_2|}\Big)^{1/2}dy_2
$$
Let $I_{x_1,t_1}$ be the smallest dyadic interval such that $x_1,t_1\in I_{x_1,t_1}$ and let $(I_k)_{k\geq 0}$ the family of dyadic intervals such that
$I_{x_1,t_1}\subset I_k$ with $|I_k|=2^{k}|I_{x_1,t_1}|$.
Since
$|\partial_{t_1}\varphi_{R_1}(t_1)|\leq |R_1|^{-2}\phi_{R_1}(t_1)$
and $|\partial_{t_2}\psi_{R_2}(t_2)|\leq |R_2|^{-3/2}\phi_{R_2}(t_2)$
with $\phi_{R_{i}}$ a bump function $L^{\infty }$ adapted to $R_{i}$, a similar argument as before
allow us to bound the first of the previous integrals by
$$
\sum_{k\geq 0}\int_{I_{k+1}\backslash I_{k}}\Big(\sum_{j\geq k} \frac{1}{|I_j|^4}\frac{1}{|I_j|}\frac{1}{|I_j|}\Big)^{1/2}dy
=\sum_{k\geq 0}\Big(\sum_{j\geq k}\frac{1}{|I_j|^6}\Big)^{1/2}|I_{k+1}\backslash I_{k}|
$$
$$
=\sum_{k\geq 0}\Big(\sum_{j\geq k}\frac{1}{2^{6j}|I_{x_i,t_i}|^6}\Big)^{1/2}2^{k}|I_{x_i,t_i}|
\lesssim \frac{1}{|I_{x_i,t_i}|^2}\sum_{k\geq 0}\frac{1}{2^{3k}}2^{k}
\lesssim \frac{1}{|x_i-t_i|^2}
$$
while the second one is bounded by
$$
\sum_{k\geq 0}\int_{I_{k+1}\backslash I_{k}}\Big(\sum_{j\geq k} \frac{1}{|I_j|^2}\frac{1}{|I_j|^3}\frac{1}{|I_j|}\Big)^{1/2}dy
=\sum_{k\geq 0}\Big(\sum_{j\geq k}\frac{1}{|I_j|^6}\Big)^{1/2}|I_{k+1}\backslash I_{k}|
\lesssim \frac{1}{|x_i-t_i|^2}
$$

This way,
$$
|\partial_{t_1}\partial_{x_2}K(x,t)|\lesssim \| b\|_{\BMO(\mathbb R^2)}\frac{1}{|x_1-t_1|^2}\frac{1}{|x_2-t_2|^2}
$$

On the other hand, we also have
$$
\Lambda_b(f_1\otimes f_2,g_1\otimes g_2)
=\int f_1(t_1)g_1(x_1)
\sum_{R} \langle b, \psi_{R}\rangle \langle f_2,\varphi_{R_2} \rangle \langle g_2,\psi_{R_2} \rangle
\varphi_{R_1}(t_1)\psi_{R_1}(x_1)dt_1dx_1
$$
and we obtain the integral representation regardless disjointness of supports of the argument funtions. Moreover,
$$
\Lambda_{x_1,t_1}^2(f_2,g_2)=\Big\langle b, \sum_{R} \langle f_2,\varphi_{R_2} \rangle \langle g_2,\psi_{R_2} \rangle
\varphi_{R_1}(t_1)\psi_{R_1}(x_1)\psi_{R}\Big\rangle
$$
and we check the two properties of a WB-CZ condition: for all bump functions $f_2,g_2$ which are
$L^2$-adapted to
the same interval,
$$
|\Lambda_{x_1,t_1}^2(f_2,g_2)|\leq \| b\|_{\BMO(\mathbb R^2)}
\Big\| \sum_{R} \langle f_2,\varphi_{R_2} \rangle \langle g_2,\psi_{R_2} \rangle
\varphi_{R_1}(t_1)\psi_{R_1}(x_1)\psi_{R}\Big\|_{H^1(\mathbb R^2)}
$$
As before, the $H^1$-norm is equivalent to
$$
\Big\| S\Big(\sum_{R} \langle f_2,\varphi_{R_2} \rangle \langle g_2,\psi_{R_2} \rangle
\varphi_{R_1}(t_1)\psi_{R_1}(x_1)\psi_{R}\Big)\Big\|_{L^1(\mathbb R^2)}
$$
$$
=\int \Big(\sum_{R} \langle f_2,\varphi_{R_2} \rangle^2 \langle g_2,\psi_{R_2} \rangle^2
\varphi_{R_1}(t_1)^2\psi_{R_1}(x_1)^2\frac{\chi_{R}(y)}{|R|}\Big)^{1/2}dy
$$
$$
=\int \Big(\sum_{R_1} \varphi_{R_1}(t_1)^2\psi_{R_1}(x_1)^2\frac{\chi_{R_1}(y_1)}{|R_1|}\Big)^{1/2}dy_1
\int \Big(\sum_{R_2} \langle f_2,\varphi_{R_2} \rangle^2 \langle g_2,\psi_{R_2} \rangle^2 \frac{\chi_{R_2}(y_2)}{|R_2|}\Big)^{1/2}dy_2
$$
$$
\lesssim \frac{1}{|x_1-t_1|}\int M(f_2)(y_2)S(g_2)(y_2)dy_2
$$
$$
\leq C\| f_2\|_{L^2(\mathbb R)} \| g_2\|_{L^2(\mathbb R)} \frac{1}{|x_1-t_1|}
\leq C\frac{1}{|x_1-t_1|}
$$
Finally, we need to prove the analog estimates for
$(\Lambda_{x_1,t_1}^2-\Lambda_{x_1',t_1'}^2)(f_2,g_2)$ which will be deduced from
$|\partial_{t_1}\Lambda_{x_1,t_1}^2(f_2,g_2)|
+|\partial_{x_1}\Lambda_{x_1,t_1}^2(f_2,g_2)|\leq C|x_1-t_1|^{-2}$. By symmetry, we work only with one of such
terms:
$$
\partial_{t_1}\Lambda_{x_1,t_1}^2(f_2,g_2)=\Big\langle b, \sum_{R} \langle f_2,\varphi_{R_2} \rangle \langle g_2,\psi_{R_2} \rangle
\partial_{t_1}\varphi_{R_1}(t_1)\psi_{R_1}(x_1)\psi_{R}\Big\rangle
$$
and therefore, as before,
$$
|\partial_{t_1}\Lambda_{x_1,t_1}(f_2,g_2)|\leq \| b\|_{\BMO(\mathbb R^2)}
\Big\| \sum_{R} \langle f_2,\varphi_{R_2} \rangle \langle g_2,\psi_{R_2} \rangle
\partial_{t_1}\varphi_{R_1}(t_1)\psi_{R_1}(x_1)\psi_{R}\Big\|_{H^1(\mathbb R^2)}
$$
being the $H^1$-norm is equivalent to
$$
\Big\| S\Big(\sum_{R} \langle f_2,\varphi_{R_2} \rangle \langle g_2,\psi_{R_2} \rangle
\partial_{t_1}\varphi_{R_1}(t_1)\psi_{R_1}(x_1)\psi_{R}\Big)\Big\|_{L^1(\mathbb R^2)}
$$
$$
=\int \Big(\sum_{R_1} (\partial_{t_1}\varphi_{R_1}(t_1))^2\psi_{R_1}(x_1)^2\frac{\chi_{R_1}(y_1)}{|R_1|}\Big)^{1/2}dy_1
\int \Big(\sum_{R_2} \langle f_2,\varphi_{R_2} \rangle^2 \langle g_2,\psi_{R_2} \rangle^2 \frac{\chi_{R_2}(y_2)}{|R_2|}\Big)^{1/2}dy_2
$$
$$
\lesssim \frac{1}{|x_1-t_1|^2}\int M(f_2)(y_2)S(g_2)(y_2)dy_2
\leq C\frac{1}{|x_1-t_1|^2}
$$

We end this lemma by 
noticing that despite we trivially have, 
$$
\Lambda_b (f,1\otimes g_{2})=\Lambda_b (f,g_{1}\otimes 1)=0
$$
the following two equalities
$$
\Lambda_b(f_{1}\otimes 1,g_{1}\otimes \psi )=\Lambda_b(1\otimes f_{2},\psi \otimes g_{2})=0
$$
do not hold in general (as one can easily check in the particular case of $b=b_{1}\otimes b_{2}$). 
We will deal with this issue when constructing the third type of paraproducts. 
Because of this, we prove now a property which we will later need, that is, that 
$\Lambda_b(\psi_{R_{1}}\otimes 1,\psi_{S_{1}}\otimes \cdot )\in \BMO(\mathbb R)$ with norms uniformly 
bounded by $\| b\|_{\BMO(\mathbb R^{2})}$. For every 
$\psi \in H^{1}(\mathbb R)$ we have
$$
\Lambda_b(\psi_{R_{1}}\otimes 1,\psi_{S_{1}}\otimes \psi )
=\sum_{R'} \langle b, \psi_{R'}\rangle \langle \psi_{R_{1}}\otimes 1, \varphi_{R'}\rangle 
\langle \psi_{S_{1}}\otimes \psi ,\psi_{R'}\rangle
$$
and due to orthogonality $R'_{1}=S_{1}$ and previous expression equals 
$$
\sum_{R'_{2}} \langle b, \psi_{S_{1}\otimes R'_{2}}\rangle \langle \psi_{R_{1}}, \varphi_{S_{1}}\rangle 
\langle \psi ,\psi_{R'_{2}}\rangle
=\langle \psi_{R_{1}}, \varphi_{S_{1}}\rangle \langle b, \psi_{S_{1}}\otimes \psi \rangle 
$$
Therefore, since $\langle \psi_{R_{1}}, \varphi_{S_{1}}\rangle =|R_{1}|^{-1/2}$ when $S_{1}\subset R_{1}$ and 
zero otherwise, we have
$$
|\Lambda_b(\psi_{R_{1}}\otimes 1,\psi_{S_{1}}\otimes \psi )|
=|\langle \psi_{R_{1}}, \varphi_{S_{1}}\rangle | |\langle b, \psi_{S_{1}}\otimes \psi \rangle |
$$
$$
\leq \frac{1}{|R_{1}|^{1/2}} \| b\|_{\BMO(\mathbb R^{2})} \| \psi_{S_{1}}\|_{H^{1}(\mathbb R)}
\| \psi \|_{H^{1}(\mathbb R)}
$$
and since $\psi_{S_{1}}$ is an $L^{2}$ atom, we can bound by
$$
\frac{1}{|R_{1}|^{1/2}} \| b\|_{\BMO(\mathbb R^{2})} |S_{1}|^{1/2}\| \psi \|_{H^{1}(\mathbb R)}
\leq \| b\|_{\BMO(\mathbb R^{2})}\| \psi \|_{H^{1}(\mathbb R)}
$$
proving the statement. 

\vskip15pt
We continue with the so-called mixed paraproduct, that is, the one associated with $T_1(1)$.
\begin{lemma}{(Mixed paraproducts).}\label{paraproducts2}
Given a function $b$ in $\BMO(\mathbb R^2)$, there exists a bounded bilinear form $\Lambda_b^2$ such that
\begin{enumerate}
\item $\Lambda_b^2(1\otimes \psi_{1},\psi_{2} \otimes 1)=b$, and $\Lambda_b^2(\psi_{1} \otimes 1, 1\otimes \psi_{2} )=0$  for every $\psi_{i} \in {\mathcal S}(\mathbb R)$ with mean zero
\item $\Lambda_b^2(1\otimes 1,\psi)=\Lambda_b^2(\psi , 1\otimes 1)=0$ for every $\psi \in {\mathcal S}(\mathbb R^{2})$ with mean zero
\end{enumerate}

\end{lemma}
\proof
Using the same wavelet basis as in the previous lemma, we define
$$
\Lambda_b^2(f,g)=\sum_R \langle b, \psi_{R}\rangle \langle f, \varphi_{R_1}\otimes \psi_{R_2}\rangle \langle g,\psi_{R_1}\otimes \varphi_{R_2}\rangle
$$
which in this proof we will just denote by $\Lambda_b$.

At least formally, $\Lambda_b$ satisfies
$$
\Lambda_b(1\otimes f_2,g_1\otimes 1)=\langle b,f_2\otimes g_1\rangle
$$
$$
\Lambda_b(1,g)=\Lambda_b(f,1)=\Lambda_b(f_1\otimes 1,1\otimes g_2)=0
$$
and
$$
\Lambda_b(f_1\otimes 1,g)
=\Lambda_b (f,1\otimes g_2)=0
$$

On the other hand, as before, we have that the equalities 
$$
\Lambda_b(1\otimes f_{2},g)
=\Lambda_b (f,g_1\otimes 1)=0
$$
do not hold in general. Also as before, we can prove that 
$\Lambda_b(1\otimes \psi_{R_{2}}\otimes 1,\cdot \otimes \psi_{S_{2}})\in \BMO(\mathbb R)$ with norms uniformly 
bounded by $\| b\|_{\BMO(\mathbb R^{2})}$.

\vskip10pt
For the proof of their boundedness we proceed as before.
$$
\Lambda_b(f,g)
=\sum_R \langle b, \psi_{R}\rangle \langle f, \varphi_{R_1}\otimes \psi_{R_2}\rangle \langle g,\psi_{R_1}\otimes \varphi_{R_2}\rangle
$$
$$
=\Big\langle b, \sum_R \langle f, \varphi_{R_1}\otimes \psi_{R_2}\rangle \langle g,
\psi_{R_1}\otimes \varphi_{R_2}\rangle \psi_{R}\Big\rangle
$$
and then, using the duality $H^1(\mathbb R^2)-\BMO(\mathbb R^2)$ we have
$$
|\Lambda_b(f,g)|\leq \| b\|_{\BMO(\mathbb R^2)}
\| \sum_R \langle f, \varphi_{R_1}\otimes \psi_{R_2}\rangle \langle g,\psi_{R_1}\otimes \varphi_{R_2}\rangle  \psi_{R}\|_{H^1(\mathbb R^2)}
$$
Just assuming the sum is finite, we get that
$\sum_R \langle f, \varphi_{R_1}\otimes \psi_{R_2}\rangle \langle g,\psi_{R_1}\otimes \varphi_{R_2}\rangle  \psi_{R}\in H^1(\mathbb R^2)$ and then
$$
\| \sum_R \langle f, \varphi_{R_1}\otimes \psi_{R_2}\rangle \langle g,\psi_{R_1}\otimes \varphi_{R_2}\rangle  \psi_{R}\|_{H^1(\mathbb R^2)}
\approx  \| S(\sum_R \langle f, \varphi_{R_1}\otimes \psi_{R_2}\rangle \langle g,\psi_{R_1}\otimes 
\varphi_{R_2}\rangle  \psi_{R})\|_{L^1(\mathbb R^2)}
$$
with implicit constants independent of the number of terms in the sum. Since
$$
S\Big(\sum_R \langle f, \varphi_{R_1}\otimes \psi_{R_2}\rangle \langle g,\psi_{R_1}\otimes \varphi_{R_2}\rangle \psi_{R}\Big)^2(x,y)
$$
$$
=\sum_R \langle f, \varphi_{R_1}\otimes \psi_{R_2}\rangle^2 \langle g,\psi_{R_1}\otimes \varphi_{R_2}\rangle^2
\frac{\chi_{R_1}}{|R_1|}(x)\frac{\chi_{R_2}}{|R_2|}(y)
$$
$$
\leq \Big(\sup_{R_1}\sum_{R_2} \langle f, \varphi_{R_1}\otimes \psi_{R_2}\rangle^2 \chi_{R_1}(x) \frac{\chi_{R_2}}{|R_2|}(y)\Big)
\Big(\sum_{R_1} \sup_{R_2}\langle g,\psi_{R_1}\otimes \varphi_{R_2}\rangle^2 \frac{\chi_{R_1}}{|R_1|}(x)\chi_{R_2}(y)\Big)
$$
$$
=M_1S_2(f)^2(x,y)S_1M_2(g)^2(x,y)
$$
where the given expressions are not a composition of operators but just notation. We first prove that those operators are bounded on $L^p(\mathbb R^2)$.
We do so by applying Fefferman-Stein's inequality to $S_iM_j$: by denoting $g_{y}(x)=g(x,y)$, we have
$$
\Big\| \Big( \sum_{R_1} \sup_{R_2}\langle g,\psi_{R_1}\otimes \varphi_{R_2}\rangle^2 \chi_{R_2}\frac{\chi_{R_1}}{|R_1|}\Big)^{1/2}\Big\|_{L^{p'}(\mathbb R^2)}
\leq \Big\| \Big( \sum_{R_1} M(\langle g,\psi_{R_1}\rangle )^2 \frac{\chi_{R_1}}{|R_1|}\Big)^{1/2}\Big\|_{L^{p'}(\mathbb R^2)}
$$
$$
=\Big( \int_{\mathbb R} \Big\| \Big( \sum_{R_1} M(\langle g,\psi_{R_1}\rangle )^2(y) 
\frac{\chi_{R_1}}{|R_1|}\Big)^{1/2}\Big\|_{L^{p'}(\mathbb R)}^{p'}dy\Big)^{1/p'}
$$
$$
\leq C\Big( \int_{\mathbb R}\Big\| \Big( \sum_{R_1} \langle g_{y},\psi_{R_1}\rangle ^2 
\frac{\chi_{R_1}}{|R_1|}\Big)^{1/2}\Big\|_{L^{p'}(\mathbb R)}^{p'}dy\Big)^{1/p'}
$$
$$
=C\Big(\int_{\mathbb R}\| S(g_{y})\|_{L^{p'}(\mathbb R)}^{p'}dy\Big)^{1/p'}
\leq C\Big( \int_{\mathbb R}\| g_{y}\|_{L^{p'}(\mathbb R)}^{p'}dy\Big)^{1/p'}= C\| g\|_{L^{p'}(\mathbb R^2)}
$$
The other operator is easier because of the pointwise inequality $M_iS_j\leq S_jM_i$.

Then, we finally get
$$
\| S(\sum_R \langle f, \varphi_{R}\rangle \langle g,\psi_{R}\rangle \psi_{R})\|_{L^1(\mathbb R^2)}
$$
$$
\leq \| (M_1S_2)(f)\|_{L^p(\mathbb R^2)}\| (S_1M_2)(g)\|_{L^{p'}(\mathbb R^2)}
\leq C\| f\|_{L^p(\mathbb R^2)}\| g\|_{L^{p'}(\mathbb R^2)}
$$

To prove that this operators belong to the class of operators with a product Calder\'on-Zygmund kernel satisfying 
the WB-CZ condition, we apply the same reasoning as in the case of classical paraproducts. We do not write the details. 

\vskip15pt
Now we are going to construct the last class of paraproducts, the ones associated with the terms 
$\langle T(1\otimes \cdot ),\cdot \otimes \cdot \rangle $.

\begin{lemma}{(Third type of paraproducts).}\label{paraproducts3}
Let $b=(b_{R_{2},S_{2}})$ be a sequence of functions in $\BMO(\mathbb R)$,
parametrized by dyadic intervals $R_2,S_2$
in such way that
\begin{equation}\label{brs}
\| b_{R_2,S_2}\|_{\BMO(\mathbb R)}\lesssim \Big(\frac{\min(|R_2|,|S_2|)}{\max(|R_2|,|S_2|)}\Big)^{1/2+\delta }
\Big(\frac{\diam(R_2\cup S_2)}{\max(|R_2|,|S_2|)}\Big)^{-(1+\delta )}
\end{equation}

Let also $(\psi_{I})_{I}$ be a wavelet basis in $L^{2}(\mathbb R)$ such that every function $\psi_{I}$ is a bump function supported in the 
dyadic interval $I$ and adapted to $I$ with constant $C>0$.

Then, there exists a bounded bilinear form $\Lambda_{b}^3$ such that
\begin{enumerate}
\item $\Lambda_{b}^3 (1\otimes 1,\cdot \otimes \cdot )=\Lambda_{b}^3 (\cdot \otimes \cdot, 1\otimes 1)
=\Lambda_{b}^3 (\cdot \otimes 1, 1\otimes \cdot )=\Lambda_{b}^3 (1\otimes \cdot ,\cdot \otimes 1)=0$,
\item $\Lambda_{b}^3(1\otimes \psi_{R_{2}},\cdot \otimes \psi_{S_{2}} )=b_{R_{2},S_{2}}$
\item $\Lambda_{b}^3 (\cdot \otimes 1,\cdot \otimes \cdot )
=\Lambda_{b}^3(\cdot \otimes \cdot , 1\otimes \cdot )=\Lambda_{b}^3(\cdot \otimes \cdot ,\cdot \otimes 1)=0$
\end{enumerate}
\end{lemma}
\proof
We define
$$
\Lambda_{b}^3(f,g)
=\sum_{R_{2},S_{2}}\sum_{R_{1}} \langle b_{R_{2},S_{2}}, \psi_{R_1}\rangle
\langle f, \varphi_{R_1}\otimes \psi_{R_2}\rangle
\langle g, \psi_{R_1}\otimes \psi_{S_2}\rangle
$$
which we will just denote by $\Lambda_b$.

At least formally, $\Lambda_b $ satisfies
$$
\Lambda_b(1,g)=\Lambda_b(f,1)=\Lambda_b(f_1\otimes 1,1\otimes g_2)
=\Lambda_b(1\otimes f_2,g_1\otimes 1)=0
$$
$$
\Lambda_b(f_1\otimes 1,g_1\otimes g_2)=\Lambda_b(f_1\otimes f_2,1\otimes g_2)
=\Lambda_b(f_1\otimes f_2,g_1\otimes 1)=0
$$
being the proof trivial in all cases. It also trivially satisfies
$$
\Lambda_b(1\otimes \psi_{R_{2}},\psi \otimes \psi_{S_{2}})
=\sum_{R_1} \langle b_{R_2,S_2}, \psi_{R_1}\rangle
\langle \psi,\psi_{R_1}\rangle 
$$
$$
=\langle b_{R_2,S_2}, \sum_{R_1}\langle \psi,\psi_{R_1}\rangle \psi_{R_1}\rangle
=\langle b_{R_2,S_2},\psi \rangle
$$

We prove now boundedness of $\Lambda_{b}$. Notice that 
$$
\Lambda_{b} (f,g)
=\sum_{R_2,S_2}\sum_{R_1} \langle b_{R_2,S_2}, \psi_{R_1}\rangle
\langle \langle f, \psi_{R_2}\rangle,\varphi_{R_{1}}\rangle 
\langle \langle g, \psi_{S_2}\rangle ,\psi_{R_{1}}\rangle 
$$
$$
=\sum_{R_2,S_2}\Big\langle {\prod}_{h,l}(b_{R_{2},S_{2}}, \langle f, \psi_{R_2}\rangle),\langle g, \psi_{S_2}\rangle \Big\rangle
$$ 
where $ {\prod}_{h,l}(b_{R_{2},S_{2}}, \langle f, \psi_{R_2}\rangle)
=\sum_{R_1}\langle b_{R_{2},S_{2}},\psi_{R_1} \rangle
\langle \langle f,\psi_{R_{2}}\rangle, \varphi_{R_1}\rangle \psi_{R_1} $ is a classical one-parameter paraproduct. 

Then, by boundedness of classical paraproducts and the property of $b_{R_{2},S_{2}}$ we have 
$$
|\Lambda_{b} (f,g)|
\leq \sum_{R_2,S_2}\| b_{R_2,S_2}\|_{\BMO(\mathbb R)}
\| \langle f, \psi_{R_2}\rangle\|_{L^{p}(\mathbb R)}\| \langle g, \psi_{S_2}\rangle\|_{L^{p'}(\mathbb R)}
$$

 As we have done several times before, we
 parametrize the terms accordingly with eccentricity and relative distances 
and rewrite previous expression as
 $$
 \sum_{e\in \mathbb Z}\sum_{m\in \mathbb N}
 \sum_{\tiny \begin{array}{c}(R_2,S_2)\in P_{e,m}\end{array}}
 \| b_{R_2,S_2}\|_{\BMO(\mathbb R)}
 \| \langle f, \psi_{R_2}\rangle\|_{L^{p}(\mathbb R)}\| \langle g, \psi_{S_2}\rangle\|_{L^{p'}(\mathbb R)}
 $$
 where $P_{e,m}$ is the set of pairs of dyadic intervals $(R_2,S_2)$ such that
 $|R_2|=2^{e}|S_2|$ and
 $n\leq \frac{\diam(R_2\cup S_2)}{\max(|R_2|,|S_2|)}<n+1$.
 Notice that $(R_2,S_2)\in P_{e,m}$ if and only if $(S_2,R_2)\in P_{-e,m}$.
 Moreover, because of the hypothesis (\ref{brs})
 we can bound previous display by
 $$
 C\sum_{e\in \mathbb Z}\sum_{m\in \mathbb N}
 2^{-|e|(1/2+\delta )}m^{-(1+\delta )}
 \sum_{\tiny \begin{array}{c}(R_2,S_2)\in P_{e,m}\end{array}}
 \| \langle f, \psi_{R_2}\rangle\|_{L^{p}(\mathbb R)}\| \langle g, \psi_{S_2}\rangle\|_{L^{p'}(\mathbb R)}
 $$
 
 Now,
 we denote by $K_2$ the interval minimum, that is, such that $K_2=R_2$ if $|R_2|\leq |S_2|$ and $K_2=S_2$ otherwise.
 This way, the inner sum can be rewritten as
 $$
 \sum_{\tiny \begin{array}{c}(R_2,S_2)\in P_{e,m}\end{array}}\int_{\mathbb R}
 \| \langle f, \psi_{R_2}\rangle\|_{L^{p}(\mathbb R)}\| \langle g, \psi_{S_2}\rangle\|_{L^{p'}(\mathbb R)}
 \frac{\chi_{K_2}(x_2)}{|K_2|}dx_2
 $$
 \begin{equation}\label{formodsquare2}
 \leq \int_{\mathbb R}\Big( \sum_{\tiny \begin{array}{c}(R_2,S_2)\in P_{e,m}\end{array}}
 \| \langle f, \psi_{R_2}\rangle\|_{L^{p}(\mathbb R)}^2\frac{\chi_{K_2}(x_2)}{|K_2|}\Big)^{1/2}
 \Big( \sum_{\tiny \begin{array}{c}(R_2,S_2)\in P_{e,m}\end{array}}
 \| \langle g, \psi_{S_2}\rangle\|_{L^{p'}(\mathbb R)}
^2\frac{\chi_{K_2}(x_2)}{|K_2|}\Big)^{1/2}dx_2
 \end{equation}

 In order to build up the modified square functions,
 we denote by $k,k'\in \mathbb Z$ the scale parameter $k=\max(e,0)$, $k'=-\min(e,0)$
 and by $n\in \mathbb Z$ the translation parameters
 $n=m$, $n'=1$ if $e\geq 0$ while $n=1$, $n'=m$ if $e\leq 0$.
 We show how to bound the first factor.
 By the choice of $K_2$ we have that $|R_2|,|S_2|\geq |K_2|$.
 If $k\geq 0$ then $K_2=R_2$ and there is nothing to show. So, we may assume $k\leq 0$ and $K_2=S_2$ in which case
 $$
 \tilde{S}_{k,n}(f)(x_2)=\Big( \sum_{(R_2,S_2)\in {\cal P}_{k,n}}
 \frac{\| \langle f, \psi_{R_2}\rangle\|_{L^{p}(\mathbb R)}^2}{|S_2|}\chi_{S}(x_2)\Big)^{1/2}
 =\Big( \sum_{R_2}\| \langle f, \psi_{R_2}\rangle\|_{L^{p}(\mathbb R)}^2\hspace{-.5cm}
 \sum_{\tiny \begin{array}{c}S_2\\(R_2,S_2)\in {\cal P}_{k,n}\end{array}}
 \hspace{-.4cm}\frac{\chi_{S_2}(x_2)}{|S_2|}\Big)^{1/2}
 $$
 $$
 =\Big( \sum_{R_2}\| \langle f, \psi_{R_2}\rangle\|_{L^{p}(\mathbb R)}^2\frac{\chi_{\tilde{S_2}}(x_2)}{2^{-k}|\tilde{S_2}|}\Big)^{1/2}
 =2^{k/2}S_{0,n}(f)(x_2)
 $$
 where
 $\tilde{S_2}$ is the dyadic interval such that $|\tilde{S_2}|=|S_2|$ and
 $|S_2|^{-1}\diam(\tilde{S}_2\cup S_2)=n$, while
 $$
 S_{0,n}(f)=\Big( \sum_{R_2,S_2}\| \langle f, \psi_{R_2}\rangle\|_{L^{p}(\mathbb R)}^2\frac{\chi_{R_2}}{|R_2|}\Big)^{1/2}
 $$

 This implies that expression (\ref{formodsquare2}) is equal to
 $$
 \int_{\mathbb R^2}2^{\max(e,0)/2}S_{0,n}(f_2)(x)2^{-\min(e,0)/2}S_{0,n'}(g_2)(x)dx
 \leq 2^{|e|/2}\| S_{0,n}(f_2)\|_{L^{p}(\mathbb R)} \| S_{0,n'}(g_2)\|_{L^{p'}(\mathbb R)}
 $$
 As indicated in the remark \ref{vectorvalued}, the vector-valued counterpart of the modified one-parameter square function also satisfies the 
 analog boundedness result of Proposition \ref{modifieddoublesquarefunction}. Therefore, 
 previous display can be bounded by
 $$
 C_p\, 2^{|e|/2} (\log(n)+1)^{|\frac{2}{p}-1|}
 \| f_2\|_{L^{p}(\mathbb R)}(\log(n')+1)^{|\frac{2}{p}-1|}\| g_2\|_{L^{p'}(\mathbb R)}
 $$
 $$
 \leq C_p\, 2^{|e|/2} n^{\epsilon |\frac{2}{p}-1|}(n')^{\epsilon |\frac{2}{p}-1|}
 \| f_2\|_{L^{p}(\mathbb R)}\| g_2\|_{L^{p'}(\mathbb R)}
 $$
 $$
 =C_{p}\, 2^{|e|/2}m^{\epsilon |\frac{2}{p}-1|}
 \| f_2\|_{L^{p}(\mathbb R)}\| g_2\|_{L^{p'}(\mathbb R)}
 $$
 since $nn'=m$. 
 
 Then, putting everything back together, we have
 $$
 |\Lambda_b (f,g)|
 \leq C_p \sum_{e\in \mathbb Z}\sum_{m\in \mathbb N}2^{-|e|(1/2+\delta )}m^{-(1+\delta )}
 2^{|e|/2} m^{\epsilon |\frac{2}{p}-1|}
 $$
 $$
 =C_{p}\sum_{e\in \mathbb Z}2^{-|e|\delta }\sum_{m\in \mathbb N}m^{-(1+\delta -\epsilon |\frac{2}{p}-1|)}
 \| f_2\|_{L^{p}(\mathbb R)} \| g_2\|_{L^{p'}(\mathbb R)}
 \leq C_{p}
 \| f_2\|_{L^{p}(\mathbb R)} \| g_2\|_{L^{p'}(\mathbb R)}
 $$
 as long as  $\epsilon |\frac{2}{p}-1|<1+\delta $.

We prove now that this family of operators also belong to the class of operators for which the theory applies.
The task will require a much finer analysis of the operator than in the case of the two previous paraproducts. 
We need to show that the operators satisfy the integral representation stated in Definition \ref{intrep} with
a kernel satisfying the Definition \ref{prodCZ} of a product Calderon-Zygmund kernel. Then we will also need to prove that the kernel satisfies the WB-CZ condition of Definition \ref{WB-CZ}. However, due to the length of computations, will 
only show explicitly the kernel decay being analogous the proof of the smoothness condition and the WB-CZ condition. Moreover, in order to simplify computations, from now we will assume that the wavelets basis used to define the bilinear form $\Lambda_{b} $ is actually the Haar basis. 

From
$$
\Lambda_b(f,g)
=\sum_{R_{2},S_{2}}\sum_{R_{1}} \langle b_{R_{2},S_{2}}, \psi_{R_1}\rangle
\langle f,\varphi_{R_1} \otimes \psi_{R_2}\rangle
\langle g,\psi_{R_1}\otimes \psi_{S_2}\rangle
$$
$$
=\int f(t)g(x)
\sum_{R_{2},S_{2}}\sum_{R_{1}} \langle b_{R_{2},S_{2}},\psi_{R_1}\rangle
\varphi_{R_1}(t_1)\psi_{R_2}(t_2)
\psi_{R_1}(x_1)\psi_{S_2}(x_2)dtdx
$$
we directly obtain the integral representation regardless disjointness of supports of the argument functions: 
$$
K(x,t)=\Big\langle \sum_{R_{2},S_{2}}b_{R_{2},S_{2}}\psi_{R_2}(t_2)\psi_{S_2}(x_2), 
\sum_{R_1} \varphi_{R_1}(t_1)\psi_{R_1}(x_1)\psi_{R_1}\Big\rangle
$$
and we check the two properties of a product C-Z kernel. By duality,
$$
|K(x,t)|\leq \Big\| \sum_{R_{2},S_{2}}b_{R_{2},S_{2}}\psi_{R_2}(t_2) \psi_{S_2}(x_2) \Big\|_{\BMO(\mathbb R)} 
\Big\| \sum_{R_1} \varphi_{R_1}(t_1)\psi_{R_1}(x_1)\psi_{R_1}\Big\|_{H^{1}(\mathbb R)}
$$

As we have seen before, the $H^{1}(\mathbb R)$ norm is equivalent to 
$$
\Big\| S\Big( \sum_{R_1} \varphi_{R_1}(t_1)\psi_{R_1}(x_1)\psi_{R_1}\Big)\Big\|_{L^{1}(\mathbb R)}
\lesssim \frac{1}{|x_{1}-t_{1}|}
$$
and therefore 
$$
|K(x,t)|\leq \frac{1}{|x_{1}-t_{1}|}
\Big\| \sum_{R_{2},S_{2}}b_{R_{2},S_{2}}\psi_{R_2}(t_2) \psi_{S_2}(x_2) \Big\|_{\BMO(\mathbb R)} 
$$

Let $I_{x_{2},t_{2}}$ be the smallest dyadic interval such that $x_{2},t_{2}\in I_{x_{2},t_{2}}$, which satisfies 
$|I_{x_{2},t_{2}}|\approx |x_{2}-t_{2}|$. 
Let also $K_{2}$ be the  interval maximum, that is, such that $K_2=R_2$ if $|S_2|\leq |R_2|$ ($e\geq 0$) and $K_2=S_2$ otherwise ($e\leq 0$). Then, we break the analysis into two different cases: when $|K_{2}|\geq |x_{2}-t_{2}|$ and when $|K_{2}|<|x_{2}-t_{2}|$.

In the first case, we bound crudely 
$$
\Big\| \sum_{R_{2},S_{2}}b_{R_{2},S_{2}}\psi_{R_2}(t_2) \psi_{S_2}(x_2) \Big\|_{\BMO(\mathbb R)} 
\leq  \sum_{R_{2},S_{2}}\| b_{R_{2},S_{2}}\|_{\BMO(\mathbb R)} |\psi_{R_2}(t_2)| |\psi_{S_2}(x_2)| 
$$
Then, we use the properties $|\psi_{R_{2}}(t_{2})|\leq |R_{2}|^{-1/2}\chi_{R_{2}}(t_{2})$, 
$|\psi_{S_{2}}(x_{2})|\leq |S_{2}|^{-1/2}\chi_{S_{2}}(x_{2})$
and parametrize the sum by eccentricity and relative distance to get the bound
\begin{equation}\label{para}
\sum_{e\in \mathbb Z}\sum_{m\in \mathbb N}2^{-|e|(1/2+\delta )}(1+m)^{-(1+\delta )}
\sum_{(R_{2},S_{2})\in P_{e,m}}\frac{1}{|R_{2}|^{1/2}}\chi_{R_{2}}(t_{2})
\frac{1}{|S_{2}|^{1/2}}\chi_{S_{2}}(x_{2})
\end{equation}


Since $|K_{2}|\geq |I_{x_{2},t_{2}}|$, we have $I_{x_{2},t_{2}}\subset K_{2}$ and $R_{2}\cap S_{2}\neq \emptyset $ which implies $m\leq 2$. 

 In the case we are considering, this family can be parametrized by the 
size, $|K_{2,k}|=2^{k}|I_{x_{2},t_{2}}|$ with $k\geq 0$.

Then, we bound (\ref{para}) by
$$
C\sum_{e\in \mathbb Z}
2^{-|e|(1/2+\delta )}
\sum_{k\geq 0}2^{|e|/2}2^{-k}\frac{1}{|x_{2}-t_{2}|}
$$
$$
=C\sum_{e\in \mathbb Z}2^{-|e|\delta }\sum_{k\geq 0}2^{-k}
\frac{1}{|x_{2}-t_{2}|}
\leq C \frac{1}{|x_{2}-t_{2}|}
$$

In the second case, when $|K_{2}|<|x_{2}-t_{2}|$, we use the fact that the sequence of functions $b_{R_{2},S_{2}}$ is generated through a bilinear form $\Lambda $.
Actually, given a $\psi \in H^{1}(\mathbb R)$ we need to bound  
$$
\Big\langle \sum_{R_{2},S_{2}\subset I_{x_{2},t_{2}}}b_{R_{2},S_{2}}\psi_{R_2}(t_2) \psi_{S_2}(x_2),\psi \Big\rangle 
=\sum_{R_{2},S_{2}\subset I_{x_{2},t_{2}}}\langle b_{R_{2},S_{2}},\psi \rangle \psi_{R_2}(t_2) \psi_{S_2}(x_2) 
$$
$$
=\sum_{R_{2},S_{2}\subset I_{x_{2},t_{2}}}
\Lambda (1\otimes \psi_{R_2}, \psi \otimes \psi_{S_2})\psi_{R_2}(t_2) \psi_{S_2}(x_2)
$$
$$
=\Lambda (1\otimes \sum_{R_{2}\subset I_{x_{2},t_{2}}}\psi_{R_2}(t_2) \psi_{R_2}, 
\psi \otimes \sum_{S_{2}\subset I_{x_{2},t_{2}}}\psi_{S_2}(x_2)\psi_{S_2})
$$

From now, we assume that both sums run over finite families of dyadic intervals and 
we work to obtain bounds of the previous expression that are independent of the cardinality of the particular families of intervals.
In this case, the functions $F_{t_{2}}=\sum_{R_{2}\subset I_{x_{2},t_{2}}}\psi_{R_2}(t_2) \psi_{R_2}$ and 
$G_{x_{2}}=\sum_{S_{2}\subset I_{x_{2},t_{2}}}\psi_{S_2}(x_2) \psi_{S_2}$ have disjoint support and mean zero. Then, by the integral representation of the bilinear form, we have that the inner term of previous expression equals 
$$
\int F_{t_{2}}(t_{2}')G_{x_{2}}(x_{2}')\Lambda_{x_{2}',t_{2}'}(1,\psi )dt_{2}'dx_{2}'
$$
Let $K$ be the support of $\psi $ and 
$\Phi$ be a bump function $L^{\infty }$-adapted and supported in $K$. We denote $c_K=c(K)$. Then, previous expression can be decomposed into
$$
\int_{\mathbb R^2} F_{t_{2}} (t_2')G_{x_{2}}(x_2') \Lambda_{x_2',t_2'}(\Phi ,\psi )dt_2'dx_2'
$$
$$
+\int_{\mathbb R^2} F_{t_{2}} (t_2')G_{x_{2}}(x_2') \Lambda_{x_2',t_2'}(1-\Phi ,\psi )dt_2'dx_2'
$$

We use the mean zero of $G_{x_{2}}$ to rewrite the first term in the following way
$$
\int_{\mathbb R^2} F_{t_{2}}(t_2')G_{x_{2}}(x_{2}') 
(\Lambda_{x_2',t_2'}(\Phi ,\psi )-\Lambda_{x_{2},t_2'}(\Phi ,\psi ))dt_2'dx_2'
$$

Now, we have $|x_2-t_2|\geq  2|x_2-x_{2}'|$.
Then, since $|K|^{-1/2}\Phi $ and $|K|^{1/2}\psi $ are $L^2$-adapted to the same interval, by the WB,
we can bound the previous expression by
$$
\int_{|x_2-t_2|>2(|x_2-x_{2}'|+|t_2-t_{2}'|)}  |F_{t_{2}}(t_2')| |G_{x_{2}}(x_{2}') | 
|(\Lambda_{x_2',t_2'}-\Lambda_{x_{2},t_2'})(|K|^{-1/2}\Phi ,|K|^{1/2}\psi )|dt_2'dx_2'
$$
$$
\lesssim \int_{|x_2-t_2|>2(|x_2-x_{2}'|+|t_2-t_{2}'|)} |F_{t_{2}}(t_2')| |G_{x_{2}}(x_{2}') |
\frac{|x_2'-x_{2}|^{\delta }}{|x_2-t_2|^{1+\delta }}dt_2'dx_2'
$$
$$
\leq \| F_{t_{2}}\|_1\| G_{x_{2}} \|_1\frac{|x_2-t_{2}|^{\delta }}{|x_2-t_2|^{1+\delta }}
=\| F_{t_{2}}\|_1^{2}\frac{1}{|x_2-t_2|}
$$

Now, we compute the $L^{1}$ norm. 
We normalize and apply Khintchine's inequaltiy to obtain for any $0<p<1$, 
$$
\int \Big| \sum_{R_{2}\subset I_{x_{2},t_{2}}}\psi_{R_2}(t_2) \psi_{R_2}(t_{2}')\Big| dt_{2}'
=|x_{2}-t_{2}|\int \Big| \sum_{R_{2}\subset I_{x_{2},t_{2}}}\psi_{R_2}(t_2) \psi_{R_2}(t_{2}')\Big| 
\frac{dt_{2}'}{|x_{t}-t_{2}|}
$$
$$
\leq |x_{2}-t_{2}|C_{p}\Big( \int \Big| \sum_{R_{2}\subset I_{x_{2},t_{2}}}\psi_{R_2}(t_2) \psi_{R_2}(t_{2}')\Big|^{p} 
\frac{dt_{2}'}{|x_{t}-t_{2}|}\Big)^{1/p}
$$
$$
=|x_{2}-t_{2}|^{1-1/p}C_{p}\Big(\sum_{k\leq 0}\int_{R_{2,k}\backslash R_{2,k-1}} 
\Big| \sum_{j=k}^{0} \psi_{R_{2,j}}(t_2) \psi_{R_{2,j}}(t_{2}')\Big|^{p} dt_{2}'\Big)^{1/p}
$$
$$
\leq |x_{2}-t_{2}|^{1-1/p}C_{p}\Big( \sum_{k\leq 0}\int_{R_{2,k}\backslash R_{2,k-1}} 
\Big|\sum_{j=k}^{0} \frac{1}{|R_{2,j}|}\Big|^{p} dt_{2}'\Big)^{1/p}
$$
$$
=|x_{2}-t_{2}|^{1-1/p}C_{p}\Big(\sum_{k\leq 0}
\Big|\sum_{j=k}^{0} \frac{1}{|R_{2,j}|}\Big|^{p} |R_{2,k}\backslash R_{2,k-1}|\Big)^{1/p}
$$
$$
=|x_{2}-t_{2}|^{1-1/p}C_{p}\Big( \sum_{k\leq 0}
\Big|\sum_{j=k}^{0} \frac{1}{2^{j}|x_{2}-t_{2}|}\Big|^{p} 2^{k-1}|x_{2}-t_{2}|\Big)^{1/p}
$$
$$
\lesssim |x_{2}-t_{2}|^{1-1/p}C_{p}\Big(\sum_{k\leq 0} \frac{1}{2^{kp}} 2^{k}\Big)^{1/p}|x_{2}-t_{2}|^{1/p-1}
=C_{p}\Big(\sum_{k\leq 0} 2^{k(1-p)}\Big)^{1/p}
\lesssim C_{p}
$$

To deal with the second term, we notice that $(1-\Phi)\otimes F_{t_{2}}$ and $\psi \otimes G_{x_{2}}$ have disjoint support
and so we can use the integral representation
$$
\int_{\mathbb R^4} (1-\Phi)(t_1')F_{t_{2}}(t_{2}')\psi(x_{1}') G_{x_{2}}(x_2') K(x',t')dt'dx'
$$
Now, because of the mean zero of both $G_{x_{2}}$ and $\psi $ we can rewrite the integral as
$$
\int_{\mathbb R^4} (1-\Phi)(t_1')F_{t_{2}}(t_{2}')\psi(x_{1}') G_{x_{2}}(x_2') 
(K(x',t')-K((c_K,x_2'),t')-K((x_1',x_{2}),t')+K((c_K,x_{2}'),t'))dt'dx'
$$
Notice that  $2|x_1'-c_K|\leq 2|K|<|x_1'-t_1'|$ and $2|x_2'-x_{2}|<|x_2-t_2|$. Then,
by the property of product C-Z kernel, we can bound by
$$
\int_{\mathbb R^4} |(1-\Phi)(t_1')| |F_{t_{2}}(t_{2}')| |\psi(x_{1}')| |G_{x_{2}}(x_2')| 
\frac{|x_1'-c_K|^{\delta }}{|x_1'-t_1'|^{1+\delta }}\frac{|x_2'-x_{2}|^{\delta }}{|x_2-t_2|^{1+\delta }}dt'dx'
$$
$$
\lesssim \frac{1}{|x_2-t_2| }\| F_{t_{2}} \|_1\| G_{x_{2}} \|_1 |K|^{\delta}
\int_{\mathbb R}|\psi(x_1')|\int_{|x_1'-t_1'|>|K|}\frac{1}{|x_1'-t_1'|^{1+\delta }}dt_1'dx_1'
$$
$$
\lesssim  C_{p}^{2} |K|^{\delta}
\|\psi \|_1\frac{1}{|K|^{\delta }}\frac{1}{|x_2-t_2| }
$$
$$
\leq C\frac{1}{|x_2-t_2| }
$$
which ends the proof of the decay of a Calder\'on-Zygmund kernel.

\vskip30pt
To apply the just constructed operators to the problem of reduction to the special cancellation 
we first need to prove the following lemma. The computations in the proof are very similar to the work we did on proving that the kernel of the third class of paraproducts are standard Calder\'on-Zygmund kernels, but with a different 
outcome. 

\begin{lemma}\label{fromWBtoBMO}
Let $\Lambda$ be a bilinear Calder\'on-Zygmund form with associated
kernel $K$ and satisfying the mixed WB-CZ conditions.

We also assume that
$\Lambda $ satisfies
the weak boundedness condition
and the cancellation conditions:
$$
\langle T(\phi_{I}\otimes 1),\varphi_{I}\otimes \cdot \rangle ,\langle T(1\otimes \phi_{I}),\cdot \otimes \varphi_{I}\rangle ,
\langle T^{*}(\phi_{I}\otimes 1),\varphi_{I}\otimes \cdot \rangle ,
\langle T^{*}(1\otimes \phi_{I}),\cdot \otimes \varphi_{I}\rangle \in \BMO (\mathbb R)
$$
for all $\phi_{I}$, $\varphi_{I}$ bump functions adapted to $I$
with norms uniformly bounded in $I$.

Then,
$$
\langle T(\phi_{I}\otimes 1),\psi_{J}\otimes \cdot \rangle ,\langle T(1\otimes \phi_{I}),\cdot \otimes \psi_{J}\rangle ,
\langle T^{*}(\phi_{I}\otimes 1),\psi_{J}\otimes \cdot \rangle ,
\langle T^{*}(1\otimes \phi_{I}),\cdot \otimes \psi_{J}\rangle \in \BMO (\mathbb R)
$$
for all $\phi_{I}$, $\psi_{J}$ bump functions adapted to $I$ and $J$ respectively such that the function associated with the interval of smallest length has mean zero. Moreover, the norms satisfy
$$
\| \langle T(\phi_{I}\otimes 1),\psi_{J}\otimes \cdot \rangle \|_{\BMO(\mathbb R)}
\lesssim \Big(\frac{\min(|I|,|J|)}{\max(|I|,|J|)}\Big)^{1/2+\delta }
\Big(\frac{\diam(I\cup J)}{\max(|I|,|J|)}\Big)^{-(1+\delta )}
$$

\end{lemma}
\proof
We assume that $\phi_I$ and $\psi_J$ are supported in $I$ and $J$ respectively. The general case needs 
some extra decomposition of the argument functions involved in the same way we did in the proof of the bump Lemma \ref{twobump}, but we will not give the full details here. 
We assume that $|J|\leq |I|$ and so that $\psi_J$ has mean zero.
 
If $|I|^{-1}\diam(I\cup J)\leq 2$ then $\phi_I$ and $\psi_J$ are both adapted to $I$ with the same constant and thus,
by hypothesis ($WB\otimes T(1)$)
$$
|\Lambda (1\otimes \phi_I,f\otimes \psi_J)|\leq C\| f\|_{H^1}
$$
for any atom $f$.

If $|I|^{-1}\diam(I\cup J)>2$, we reason as follows.
Let $f$ be an atom supported in $K$.
Since $\phi_I$ and $\psi_J$ have disjoint support, we have the following integral representation
$$
\Lambda (1\otimes \phi_I,f\otimes \psi_J)
=\int_{\mathbb R^2} \phi_I (t_2)\psi_J (x_2) \Lambda_{x_2,t_2}(1,f)dt_2dx_2
$$

Now, let $\Phi$ a bump function $L^{\infty }$-adapted and supported in $K$. We denote $c_J=c(J)$ and $c_K=c(K)$. Then,
$$
\Lambda (1\otimes \phi_I,f\otimes \psi_J)
=\int_{\mathbb R^2} \phi_I (t_2)\psi_J (x_2) \Lambda_{x_2,t_2}(\Phi ,f)dt_2dx_2
$$
$$
+\int_{\mathbb R^2} \phi_I (t_2) \psi_J (x_2) \Lambda_{x_2,t_2}(1-\Phi ,f)dt_2dx_2
$$

We use the mean zero of $\psi_J$ to rewrite the first term in the following way
$$
\int_{\mathbb R^2} \phi_I (t_2)\psi_J (x_2) (\Lambda_{x_2,t_2}(\Phi ,f)-\Lambda_{c_J,t_2}(\Phi ,f))dt_2dx_2
$$

Now, we have $|x_2-t_2|\geq \diam(I\cup J)>2|I|\geq 2|J|\geq 2|x_2-c_J|$.
Then, since $|K|^{-1/2}\Phi $ and $|K|^{1/2}f$ are $L^2$-adapted to the same interval, by the WB,
we can bound the previous expression by
$$
\int_{|x_2-t_2|>\diam (I\cup J)} |\phi_I (t_2)| |\psi_J (x_2)| |(\Lambda_{x_2,t_2}-\Lambda_{c_J,t_2})
(|K|^{-1/2}\Phi ,|K|^{1/2}f)|dt_2dx_2
$$
$$
\lesssim \int_{|x_2-t_2|>\diam (I,J)} |\phi_I (t_2)| \psi_J (x_2)|
\frac{|x_2-c_J|^{\delta }}{|x_2-t_2|^{1+\delta }}dt_2dx_2
$$
$$
\leq \|\phi_I \|_1\| \psi_J \|_1|J|^{\delta }\diam (I\cup J)^{-(1+\delta )}
\lesssim |I|^{1/2}|J|^{1/2}|J|^{\delta }\diam (I\cup J)^{-(1+\delta )}
$$
$$
=\Big( \frac{|J|}{|I|}\Big)^{1/2+\delta }(|I|^{-1}\diam (I\cup J))^{-(1+\delta )}
$$

To deal with the second term, we notice that $\phi_I\otimes (1-\Phi)$ and $\psi_J\otimes f$ have disjoint support
and so we can use the integral representation
$$
\int_{\mathbb R^4} (1-\Phi)(t_1)f(x_1)\phi_I (t_2) \psi_J (x_2) K(x,t)dtdx
$$
Now, because of the mean zero of both $\psi_J$ and $f$ we can rewrite the integral as
$$
\int_{\mathbb R^4} (1-\Phi)(t_1)f(x_1)\phi_I (t_2) \psi_J (x_2) (K(x,t)-K((c_K,x_2),t)-K((x_1,c_J),t)+K((c_K,c_J),t))dtdx
$$
Notice that  $2|x_1-c_K|\leq 2|K|<|x_1-t_1|$ and $2|x_2-c_J|<|x_2-t_2|$. Then,
by the property of product C-Z kernel, we can bound by
$$
\int_{\mathbb R^4} |(1-\Phi)(t_1)| |f(x_1)| |\phi_I (t_2)| |\psi_J (x_2)| 
\frac{|x_1-c_K|^{\delta }}{|x_1-t_1|^{1+\delta }}\frac{|x_2-c_J|^{\delta }}{|x_2-t_2|^{1+\delta }}dtdx
$$
$$
\lesssim |J|^{\delta }{\rm diam}(I\cup J)^{-(1+\delta)}\|\phi_I \|_1\| \psi_J \|_1 |K|^{\delta}
\int_{\mathbb R}|f(x_1)|\int_{|x_1-t_1|>|K|}
\frac{1}{|x_1-t_1|^{1+\delta }}dt_1dx_1
$$
$$
\lesssim |J|^{\delta }{\rm diam}(I\cup J)^{-(1+\delta)} |I|^{1/2}|J|^{1/2}|K|^{\delta}
\|f\|_1\frac{1}{|K|^{\delta }}
$$
$$
=\Big( \frac{|J|}{|I|}\Big)^{1/2+\delta }(|I|^{-1}\diam (I\cup J))^{-(1+\delta )}\|f\|_1
$$
which ends the proof of this lemma.

\begin{corollary} Let $\psi_{I}$ and $\psi_{J}$ bump functions $L^{2}$ adapted to $I$ and $J$ respectively 
with constant $C$ and mean zero. Then, for any $f$ atom, 
$$
|\Lambda (1\otimes \phi_I,f\otimes \psi_J)|
\leq C\Big(\frac{\min(|I|,|J|)}{\max(|I|,|J|)}\Big)^{1/2+\delta }
\Big(\frac{\diam(I\cup J)}{\max(|I|,|J|)}\Big)^{-(1+\delta )}\| f\|_{1}
$$
\end{corollary}
\proof
By duality and previous Lemma
$$
|\Lambda (1\otimes \phi_I,f\otimes \psi_J)|
=|\langle T(1\otimes \phi_{I}),f\otimes \psi_{J}\rangle |
$$
$$
\leq \| \langle T(\phi_{I}\otimes 1),\psi_{J}\otimes \cdot \rangle \|_{\BMO(\mathbb R)}\| f\|_{H^{1}(\mathbb R)}
$$
$$
\leq C\Big(\frac{\min(|I|,|J|)}{\max(|I|,|J|)}\Big)^{1/2+\delta }
\Big(\frac{\diam(I\cup J)}{\max(|I|,|J|)}\Big)^{-(1+\delta )}\| f\|_{L^{1}(\mathbb R)}
$$

\vskip20pt
Once previous lemma is proven, we proceed as follows. 
We first consider the functions defined by 
$$
\langle b_1, \psi \rangle =\Lambda (1\otimes 1, \psi)
\hskip 60pt  \langle b_2, \psi \rangle =\Lambda (\psi, 1\otimes 1)
$$  
$$
\langle b_3, \psi \rangle =\Lambda (\psi_1\otimes 1, 1\otimes \psi_2)
\hskip60pt
\langle b_4, \psi \rangle =\Lambda (1\otimes \psi_2, \psi_1\otimes 1)
$$
for every $\psi \in {\mathcal S}(\mathbb R^{2})$ of mean zero. By hypothesis, all of them 
are functions in $\BMO(\mathbb R^2)$ and so we can construct the bounded paraproducts $\Lambda_{b_1}^1$,
$\Lambda_{b_2}^1$ and 
$\Lambda_{b_3}^2$, $\Lambda_{b_4}^2$. We also denote 
$\Lambda_{1,2,3,4}=\sum_{i=1,2}\Lambda_{b_i}^1+\sum_{i=3,4}\Lambda_{b_i}^2$.

Moreover, in the remark int the proofs of lemmata \ref{paraproducts1} and \ref{paraproducts2}, we showed that the functions 
$\Lambda_{b_i}^j(\psi_{R_{1}}\otimes 1, \psi_{S_{1}}\otimes \cdot)$ are uniformly in $\BMO(\mathbb R)$. Then, 
we consider the sequences of functions in $\BMO(\mathbb R)$, ${\bold b}_{i}=(b^{i}_{I,J})_{I,J}$, for $i=5,6,7,8$, defined by 
$$
\langle b_{R_{1},S_{1}}^{5}, \psi \rangle =(\Lambda -\Lambda_{1,2,3,4})
(\psi_{R_{1}}\otimes 1, \psi_{S_{1}}\otimes \psi)
$$
$$
\langle b_{R_{2},S_{2}}^{6}, \psi \rangle =(\Lambda -\Lambda_{1,2,3,4}) 
(1\otimes \psi_{R_{2}}, \psi \otimes \psi_{S_{2}})
$$ 
$$
\langle b_{R_{1},S_{1}}^{7}, \psi \rangle =(\Lambda -\Lambda_{1,2,3,4})
(\psi_{R_{1}}\otimes \psi, \psi_{S_{1}}\otimes 1)
$$
$$
\langle b_{R_{2},S_{2}}^{8}, \psi \rangle =(\Lambda -\Lambda_{1,2,3,4})
(\psi \otimes \psi_{R_{2}}, 1 \otimes \psi_{S_{2}})
$$ 
for every $\psi \in {\mathcal S}(\mathbb R)$ of mean zero.
Then, we construct the paraproducts $\Lambda_{{\bold b}_i}^3$ for $i=5,6,7,8$.

We first notice that, previous equalities imply that 
for any $f_{2},g_{2}\in{\mathcal S}(\mathbb R)$, we also have 
 $
 \Lambda_{{\bold b}_{i}}(1\otimes f_{2},\psi \otimes g_{2})=\Lambda (1\otimes f_{2},\psi \otimes g_{2})
 $ 
 as we see: since $(\psi_I)_I$ is a basis, we have
 $$
 \Lambda_{{\bold b}_{i}}(1\otimes f_{2},\psi \otimes g_{2})
 =\sum_{R_{2},S_{2}}\langle f_{2},\psi_{R_{2}}\rangle \langle g_{2},\psi_{S_{2}}\rangle
 \Lambda_{{\bold b}_{i}}(1\otimes \psi_{R_{2}},\psi \otimes \psi_{S_{2}})
 $$
 $$
 =\sum_{R_{2},S_{2}}\langle f_{2},\psi_{R_{2}}\rangle \langle g_{2},\psi_{S_{2}}\rangle 
 \langle b_{R_{2},S_{2}}^{i},\psi \rangle
 =\sum_{R_{2},S_{2}}\langle f_{2},\psi_{R_{2}}\rangle \langle g_{2},\psi_{S_{2}}\rangle 
(\Lambda -\Lambda_{1,2,3,4}) (1\otimes f_{2},\psi \otimes g_{2})
 $$
 $$
 =\Lambda(1\otimes f_{2},\psi \otimes g_{2})
 -\sum_{i=1,2}\Lambda_{b_i}^1(1\otimes f_{2},\psi \otimes g_{2})
 -\sum_{i=3,4}\Lambda_{b_i}^2(1\otimes f_{2},\psi \otimes g_{2})
 $$

Moreover, because of Lemma \ref{fromWBtoBMO}, we have that the crucial hypothesis (\ref{brs}) of Lemma \ref{paraproducts3} is satisfied and then, the constructed paraproducts $ \Lambda_{{\bold b}_{i}}$ are bounded. 



Finally then, we define the bilinear form 
$$
\tilde{\Lambda }
=\Lambda -\sum_{i=1,2}\Lambda_{b_i}^1 -\sum_{i=3,4}\Lambda_{b_i}^2
-\sum_{i=5,6,7,8}\Lambda_{{\bold b}_i}^3
$$
which clearly satisfies all the required cancellation conditions
$$
\tilde{\Lambda}(1\otimes 1, \psi)
=\tilde{\Lambda }(\psi, 1\otimes 1)
=\tilde{\Lambda }(\psi_1\otimes 1, 1\otimes \psi_2)
=\tilde{\Lambda }(1\otimes \psi_2, \psi_1\otimes 1)=0
$$
and also
$$
\tilde{\Lambda }(f_{1}\otimes 1,g_{1}\otimes \psi )
=\tilde{\Lambda }(1\otimes f_{2},\psi \otimes g_{2})
=\tilde{\Lambda }(\psi \otimes f_{2},1\otimes g_{2})
=\tilde{\Lambda }(f_{1}\otimes \psi,g_{1}\otimes 1)=0
$$ 

\vskip20pt
With the three previous lemmata and the other five symmetrical statements which come from all the possible permutations of the argument functions,
we can finally prove boundedness
of product singular integrals in the
general case and finish the proof of Theorem \ref{biparT1} and also of Theorem \ref{Lpbounds} in this way.

$$
\hskip400pt \boxempty
$$

\section{Appendix}\label{umsf}
In the proof of the extension to $L^p$ spaces (see theorem \ref{Lpbounds}), we used a bi-parameter modified square function whose boundedness properties
are a direct consequence of their uni-parameter counterparts. Now,
in this appendix, we prove boundedness of such uni-parameter modified square functions.

\begin{definition}
Given $k\in \mathbb Z$, $n\in \mathbb N$, we consider
the following operator
$$
S_{k,n}(f)(x)=\big( \sum_{I}\sum_{J\in I_{k,n}}\frac{\langle f,\psi_{I}\rangle ^2}{|J|}\chi_{J}(x)\Big)^{1/2}
$$
where $I_{k,n}$ is the family of dyadic intervals $J$ satisfying $|I|=2^{k}|J|$ and
$n\leq \frac{\diam(I\cup J)}{\max(|I|,|J|)}<n+1$.
\end{definition}

We will prove bounds of such operators by means of the following modified square functions:

\begin{definition}
Given $k\in \mathbb Z$, $n\in \mathbb N$, we consider the following variant of square function
$$
\tilde{S}_{k,n}(f)(x)=\big( \sum_{I}\frac{\langle f,\psi_{I}\rangle ^2}{|J|}\chi_{J}(x)\Big)^{1/2}
$$
where $I$ and $J$ are two dyadic intervals satisfying $|I|=2^{k}|J|$ and
$n\leq \frac{\diam(I\cup J)}{\max(|I|,|J|)}<n+1$, chosen in such a way that for every interval $I$ there is a unique interval $J\in I_{k,n}$.
\end{definition}

This way we actually define a family of operators that depends on the particular choice of intervals but whose bounds do not depend
on such choice, as we will soon prove.
Notice that the particular choice does not depend on the point $x$.

We see now the
the reason why this modified square function helps to control boundedness of the previous ones. Fixed a dyadic interval $I$, we denote
$I_{k,n}$ the family of dyadic intervals $J$ such that $|I|=2^{k}|J|$ and
$n\leq \max(|I|,|J|)^{-1}\diam(I\cup J)<n+1$. We also denote by $\tilde{I}$ the dyadic interval such that
$|\tilde{I}|=|I|$ and
$|I|^{-1}\diam(I\cup \tilde{I})=n$.
Then, for all $k\geq 0$
$$
S_{k,n}(f)(x)=\big( \sum_{I}\sum_{J\in I_{k,n}}\frac{\langle f,\psi_{I}\rangle ^2}{|J|}\chi_{J}(x)\Big)^{1/2}
$$
$$
=\big( \sum_{I}\langle f,\psi_{I}\rangle ^2\sum_{J\in I_{k,n}}\frac{\chi_{J}(x)}{|J|}\Big)^{1/2}
=\big( \sum_{I}\langle f,\psi_{I}\rangle ^2\frac{\chi_{\tilde{I}}(x)}{2^{-k}|\tilde{I}|}\Big)^{1/2}
$$
$$
=2^{k/2}S_{0,n}(f)(x)=2^{k/2}\tilde{S}_{0,n}(f)(x)
$$
Meanwhile, when $k\leq 0$ we have
$$
S_{k,n}(f)(x)=\big(  \sum_{I}\sum_{J\in I_{k,n}}\frac{\langle f,\psi_{I}\rangle ^2}{|J|}\chi_{J}(x)\Big)^{1/2}
$$
$$
=\big( \sum_{J}\Big(\sum_{I\in J_{-k,n}}\langle f,\psi_{I}\rangle^2\Big) \frac{\chi_{J}(x)}{|J|}\Big)^{1/2}
\leq \big( \sum_{J}2^{-k}\langle f,\psi_{I_{J}}\rangle ^2\frac{\chi_{J}(x)}{|J|}\Big)^{1/2}
=2^{-k/2}\tilde{S}_{k,n}(f)(x)
$$
where the interval $I_{J}$ is chosen such that $|\langle f,\psi_{I_{J}}\rangle|=\max_{I\in J_{-k,n}} |\langle f,\psi_{I}\rangle|$

\begin{proposition}\label{modifiedsquarefunction} For every $1<p<\infty $, we have that if $k\geq 0$
$$
\| \tilde{S}_{k,n}f\|_{p} \leq C_{p}(2^{-k\, \sign(\frac{2}{p}-1)}\log(n)+1)^{|\frac{2}{p}-1|}\| f\|_p
$$
while if $k\leq 0$
$$
\| \tilde{S}_{k,n}f\|_{p} \leq C_{p}(2^{-k\, \sign(\frac{2}{p}-1)}+\log(n)+1)^{|\frac{2}{p}-1|}\| f\|_p
$$
with constants $C_p$ independent of $f$, $k$ and $n$.
\end{proposition}

\begin{corollary} For every $1<p<\infty $, we have that if $k\geq 0$
$$
\| S_{k,n}f\|_{p} \leq C_{p}2^{k/2}(\log(n)+1)^{|\frac{2}{p}-1|}\| f\|_p
$$
while if $k\leq 0$
$$
\| S_{k,n}f\|_{p} \leq C_{p}2^{-k/2}(2^{-k\, \sign(\frac{2}{p}-1)}+\log(n)+1)^{|\frac{2}{p}-1|}\| f\|_p
$$
with constants $C_p$ independent of $f$, $k$ and $n$.

\end{corollary}

\begin{remark}\label{vectorvalued}
Before starting with the proof, we notice that a careful read of it reveals that, by means of vector-valued
interpolation,the result also holds for vector-valued modified square function
with values in a Banach space $X$ with the UMD property
of the form
$$
\tilde{S}_{k,n}^{X}(f)(x)=\big( \sum_{I}\frac{\|\langle f,\psi_{I}\rangle \|_X^2}{|J|}\chi_{J}(x)\Big)^{1/2}
$$
for which every $1<p<\infty $, we have that if $k\geq 0$
$$
\| \tilde{S}_{k,n}f\|_{L^p(X)} \leq C_{p,X}(2^{-k\, \sign(\frac{2}{p}-1)}\log(n)+1)^{|\frac{2}{p}-1|}\| f\|_{L^p(X)}
$$
while if $k\leq 0$
$$
\| \tilde{S}_{k,n}f\|_{L^p(X)} \leq C_{p,X}(2^{-k\, \sign(\frac{2}{p}-1)}+\log(n)+1)^{|\frac{2}{p}-1|}\| f\|_{L^p(X)}
$$
Then, in particular for $X=L^p(\mathbb R)$ we get for $k\geq 0$
$$
\| \tilde{S}_{k,n}f\|_{L^p(\mathbb R^2)} \leq C_{p}(2^{-k\, \sign(\frac{2}{p}-1)}\log(n)+1)^{|\frac{2}{p}-1|}\| f\|_{L^p(\mathbb R^2)}
$$
while if $k\leq 0$
$$
\| \tilde{S}_{k,n}f\|_{L^p(\mathbb R^2)} \leq C_{p}(2^{-k\, \sign(\frac{2}{p}-1)}+\log(n)+1)^{|\frac{2}{p}-1|}\| f\|_{L^p(\mathbb R^2)}
$$

\end{remark}

The estimate for $p=2$ is a trivial consequence of Plancherel's inequality.
To extend the result to other exponents $p$ we plan to use interpolation and duality. So, we first prove the
following weak $L^1$ type estimate whose proof comes from a slight modification of the one appearing in \cite{thielelectures}.

\begin{proposition}\label{weakL1}
If $f$ is integrable and $\lambda >0$, then we have
$$
\| \{ x : \tilde{S}_{k,n}f(x)>\lambda \}\| \leq C(2^{-k}n+1)\| f\|_1 \lambda^{-1}
$$
with a constant $C$ independent of $f$ and $\lambda $.
\end{proposition}
\proof
Consider the collection $\cal I$ of maximal dyadic intervals $I$ with respect
to set inclusion such that
$$
|I|^{-1}\int_{I}|f(x)|dx >\lambda
$$
Let $E$ be the union of all $I\in {\cal I}$, the set that contains all intervals where $f$ has large average.
The intervals in $\cal I$ are pairwise disjoint and so
$$
|E|\leq \sum_{I\in \cal I}|I|\leq \lambda^{-1}\sum_{I\in \cal I}\int_{I}|f(x)|dx \leq \|f\|_{1}\lambda^{-1}
$$

We take a classical Calderon-Zymund decomposition $f=g+b$ given by
$$
g=\sum_{I} m_{I}(f)\chi_I+f\chi_{E^{c}}
\hskip 30pt b=\sum_{I} f_{I}
$$
with $m_I(f)=|I|^{-1}\int_{I}f$ and $f_I=(f-m_I(f))\chi_I$.

We see that $g$ is essentially bounded by $2\lambda $. Outside $E$, this follows by Lebesgue's
differentiation theorem. To prove this inside $E$, it suffices to consider each interval
$I\in {\cal I}$ separately. Let $I$ be such an interval and $\tilde{I}$ its parent interval. Then by maximality of $I$ we have
$$
\int_{I}f(x)dx \leq\int_{I}|f|dx \leq \int_{\tilde{I}}|f|dx\leq \lambda |\tilde{I}|=2\lambda |I|
$$
Moreover it is also clear that $$\int |g|\leq \int |f|$$

Because of the $L^2$ boundedness of $S_{k,n}$ we have
$$
\| \tilde{S}_{k,n}g \|_{2}^{2}\leq C\| g\|_{2}^{2}\leq C \int |g|\lambda dx \leq C\lambda \| f\|_1
$$
and so
$$
\{ \tilde{S}_{k,n}g > \lambda/2\}\leq C\| S_{k,n}g\|_{2}\lambda^{-2}
\leq C\| f\|_{1}\lambda^{-1}
$$

We plan to prove the same estimate for $b=\sum_{I}f_I$.
To do so, we
define $\tilde{E}$ as
the union of all $3I$ with $I\in {\cal I}$. We  also define $F$ as the union of all $J$ such that the corresponding
$I$ satisfies $I\subset I'$ for some $I'\in {\cal I}$ and $\tilde{F}$ as the union of all $3J$ with $J\subset F$. Then,
$$
| \{ \tilde{S}_{k,n}b > \lambda/2\} | \leq |\tilde{F}|+\lambda^{-1}\| \tilde{S}_{k,n}b\|_{L^1(\mathbb R\backslash \tilde{F})}
$$

Now we measure $\tilde{F}$ by means of a geometric argument that distinguishes between large and small scales.
Since $|\tilde{F}|\leq 3|F|$ and
$$
F
=\bigcup_{I'\in {\cal I}}\cup \{ J:I\subset I'\}
$$
we fix now $I'\in {\cal I}$.

We also separate between $k\geq 0$ and $k\leq 0$ since the separation in scales is slightly different. We first assume $k\geq 0$ for which
we separate into two different scales: smaller and larger than $\log(n)$.

The family of dyadic intervals $I\subset I'$ such that $|I|=2^{-r}|I'|$ with $0\leq r\leq \log(n)$
has the property that the corresponding intervals $J$ satisfy
$$
\diam (I'\cup J)\leq \diam(I\cup J)\geq n|I|=n2^{-r}|I'|>|I'|
$$
and so the intervals $J$
are disjoint with $I'$. Moreover, their union measures at most $2^{-k}|I'|$ as we see: the intervals $J$ are
pairwise disjoint and so for every $0\leq r\leq \log(n)$ we have
$$
|\bigcup \{ J:I\subset I',|I|=2^{-r}|I'|\}|= \sum_{\tiny \begin{array}{c}J:I\subset I'\\|I|=2^{-r}|I'|\end{array}} |J|
=2^{-k}\sum_{\tiny \begin{array}{c}I:I\subset I'\\|I|=2^{-r}|I'|\end{array}}|I|=2^{-k}|I'|
$$
Then
$$
|\bigcup_{r=0}^{\log(n)} \bigcup \{ J:I\subset I',|I|=2^{-r}|I'|\}|\leq 2^{-k}\log(n)|I'|
$$

On the other hand for smaller scales, that is, intervals $I\subset I'$ such that $|I|=2^{-r}|I'|$ with $r>\log(n)$, we have
that the corresponding intervals $J$ satisfy
$$
\diam (I'\cup J)\leq |I'|/2+|c(I')-c(J)|+|J|/2\leq |I'|/2+|c(I')-c(I)|+|c(I)-c(J)|+|J|/2
$$
$$
\leq |I'|+\diam(I\cup J)\leq |I'|+(n+1)|I|<|I'|+2n2^{-r}|I'|<3|I'|
$$
Then the intervals $J$ are included in $3I'$ and so
$$
|\bigcup_{r>\log(n)} \bigcup \{ J:I\subset I',|I|=2^{-r}|I'| \} |\leq C|I'|
$$

Both things show that
$$
|\cup \{ J:I\subset I' \} |\leq C(2^{-k}\log(n)+1)|I'|
$$
and therefore
$$
|\tilde{F}|\lesssim 3\sum_{I'\in {\cal I}}|\cup \{ J:I\subset I'\} |
$$
$$
\leq C(2^{-k}\log(n)+1)\sum_{I\in \cal{I}}|I'|
\leq C(2^{-k}\log(n)+1)\| f\|_1 \lambda^{-1}
$$

When $k\leq 0$, the computations are similar but we separate into three different scales: smaller than $-k$, between $-k$ and $-k+\log(n)$, and
larger than $-k+\log(n)$.

The subfamily of dyadic intervals $I\subset I'$ such that $|I|=2^{-r}|I'|$ with $0\leq r\leq -k$
has the property that the corresponding intervals $J$ satisfy
$$
\diam(I'\cup J)\geq \diam(I\cup J)\geq n|J|=n2^{-k}|I|=n2^{-k-r}|I'|>|I'|
$$ and so they
are disjoint with $I'$. Moreover,
their union measures at most $2^{-k}|I'|2^{-r}$ as we see: for all intervals $I$ considered, their corresponding intervals $J$
satisfy $|J|=2^{-k-r}|I'|>|I'|$ and so there is a unique interval $J$ corresponding with the different $I$. This way for every $0\leq r\leq -k$ we have
$$
|\bigcup \{ J:I\subset I',|I|=2^{-r}|I'|\}|\leq |J|
=2^{-k}|I|=2^{-k-r}|I'|
$$
and then summing a geometric series we have
$$
|\bigcup_{r=0}^{-k} \bigcup \{ J:I\subset I',|I|=2^{-r}|I'|\}|\leq C2^{-k}|I'|
$$
On the other hand, the subfamily of dyadic intervals $I\subset I'$ such that $|I|=2^{-r}|I'|$ with $-k\leq r\leq -k+\log(n)$
has the property that the corresponding intervals $J$ satisfy $\diam(I\cup J)\geq n2^{-k-r}|I'|>|I'|$ and so they
are still disjoint with $I'$. Moreover, their union measures at most $|I'|$ as we see: now $|J|=2^{-k-r}|I'|<|I'|$ and on varying $I$ there
are $|I'|/|J|=2^{-k-r}$ different disjoint intervals $J$ whose union measures exactly $|I'|$.  Then
for every $-k\leq r\leq \-k+\log(n)$ we have
$$
|\bigcup \{ J:I\subset I',|I|=2^{-r}|I'|\}|\leq |I'|
$$
For different $r$, the intervals $J$ are contained in a different translation of $I'$. Then the intervals $J$ are pairwise disjoint and so
$$
|\bigcup_{r=-k}^{-k+\log(n)} \bigcup \{ J:I\subset I',|I|=2^{-r}|I'|\}|\leq \log(n)|I'|
$$

Finally, for smaller scales, that is, intervals $I\subset I'$ such that $|I|=2^{-r}|I'|$ with $r>-k+\log(n)$, we have
that the corresponding intervals $J$ satisfy
$$
\diam (I'\cup J)\leq |I'|+\diam(I\cup J)\leq (n+1)|J|<|I'|+2n2^{-k-r}|I'|<3|I'|
$$
and so they are included in $3I'$ and so
$$
|\bigcup_{r>-k+\log(n)}\bigcup \{ J:I\subset I',|I|=2^{-r}|I'|\}|\leq C|I'|
$$

The three bounds together show that
$$
|\cup \{ J:I\subset I' \} |\leq C(2^{-k}+\log(n)+1)|I'|
$$
and so
$$
|\tilde{F}|\lesssim 3\sum_{I'\in {\cal I}}|\cup \{ J:I\subset I'\} |
$$
$$
\leq C(2^{-k}+\log(n)+1)\sum_{I\in \cal{I}}|I'|
\leq C(2^{-k}+\log(n)+1)\| f\|_1 \lambda^{-1}
$$

The following step of the proof is to show
$$
\| \tilde{S}_{k,n}b \|_{L^1(\mathbb R \setminus \tilde{F})} \leq C 2^{-k/2}\| f\|_1
$$
and, by sub-linearity, it suffices to prove
$$
\| \tilde{S}_{k,n}f_{I'} \|_{L^1(\mathbb R \setminus \tilde{F})} \leq C 2^{-k/2}\lambda |I'|
$$
for each $I'\in {\cal I}$.
This in turn follows from
$$
\Bigg\| \Bigg( \sum_{J\not \subset \tilde{F}}\frac{|\langle f_{I'}, \psi_{I}\rangle |^2}{|J|}\chi_{J}\Bigg)^{1/2}\Bigg\|_{1}
\leq \Bigg\| \sum_{J\not \subset {\tilde F}}\frac{|\langle f_{I'}, \psi_{I}\rangle |}{|J|^{1/2}}\chi_{J}\Bigg\|_{1}
$$
$$
\leq \sum_{J\not \subset \tilde{F}}|\langle f_{I'}, \psi_{I}\rangle ||J|^{1/2}
\leq 2^{-k/2}\sum_{I\not \subset \tilde{E}}|\langle f_{I'}, \psi_{I}\rangle ||I|^{1/2}
$$
$$
\leq C2^{-k/2}\| f_{I'}\|_{1}\leq C2^{-k/2}\lambda |I'|
$$
where we have used that $J\not \subset \tilde{F}$ implies by definition of $\tilde{F}$ that the corresponding interval $I$ satisfies
$I\not \subset \tilde{E}$. Moreover, the last inequality follows from lemma \ref{lowoscillation} below.

Finally then,
$$
| \{ \tilde{S}_{k,n}b > \lambda/2\} | \leq C(2^{-k}n+1)\| f\|_1 \lambda^{-1}+C2^{-k/2}\| f\|_1\lambda^{-1}
$$
which ends the proof since $2^{-k/2}\leq 2^{-k}+\log(n)+1\leq 2^{-k}\log(n)+1$ for all $k\in \mathbb Z$ and $n\in \mathbb N$.

\vskip 10pt
The following lemma is the technical result needed to prove proposition \ref{weakL1}. For the sake of
completeness we include its proof although is exactly the same one that can be found in \cite{thielelectures}.

\begin{lemma}\label{lowoscillation}
Let $I'$ be some interval and f be an integrable function supported
in $I'$ with mean zero. For each dyadic interval $I$ let $\phi_{I}$ be a bump function adapted
to $I$. Then
$$
\sum_{I:I\not \subset 3I'} |\langle f, \phi_{I}\rangle | |I|^{1/2}\leq C\| f\|_{1}
$$
Here $3I'$ denotes the interval that shares the center with $I'$ and is of length $3|I'|$.
\end{lemma}
\proof
We first consider the sum over all dyadic intervals $I$ such that $I\not \subset 3I'$ with $|I|<|I'|$.
Let $c$ be the midpoint between $c(I)$ and $c(I')$. By symmetry we may assume that
$\supp f\subset (-\infty ,c)$ and then,
$$
|\langle f, \phi_{I} \rangle |
\leq  \| f\|_{L^{1}(-\infty ,c)}\| \phi_{I}\|_{L^{\infty }(-\infty ,c)}
+\| f\|_{L^{\infty }(c,\infty )}\| \phi_{I}\|_{L^{1}(c,\infty )}
$$
$$
\leq \| f\|_1 C |I|^{-1/2} \Big( 1+\frac{|c(I)-c(I')|}{|I|}\Big)^{-N}
$$
which gives
$$
|\langle f, \phi_{I} \rangle | |I|^{1/2}
\leq C\| f\|_1 \Big( 1+\frac{|c(I)-c(I')|}{|I|}\Big)^{-N}
$$

For any two integers $k>0$ and $m>0$ there are at most two intervals $I$ such that
$|I'|/|I|=2^k$ and the integer part of $1+\frac{|c(I)-c(I')|}{|I|}$ is $m$. If $m<2^k$, there are no
such intervals which satisfy $I\not \subset 3I'$.
Thus we can estimate
$$
\sum_{\tiny \begin{array}{c}I:|I|<|I'|\\ I\not \subset 3I'\end{array}}|\langle f, \phi_{I}\rangle | |I|^{1/2}\leq C\| f\|_{1}\sum_{k>0}\sum_{m\geq 2^k}
m^{-N}\leq C\| f\|_1
$$

We now consider the sum over all dyadic I with $I\not \subset 3I'$ and $|I|\geq |I'|$. Let $D$ denote the operator of differentiation and
$D^{-1}$ the antiderivative operator
$$
D^{-1}f(x)=\int_{-\infty }^{x}f(y)dy
$$
Notice that because of the mean zero of $f$, the support of $D^{-1}f$ is also included in $I'$.
Then, by partial integration and the fact that $|I|D\phi_{I}$ is a bump function adapted to $I$, we have
$$
|\langle f, \phi_{I} \rangle |
=|I|^{-1}|\langle D^{-1}f, |I|D\phi_{I} \rangle |
$$
$$
\leq  |I|^{-1}(\| D^{-1}f\|_{L^{1}(-\infty ,c)}\| |I|D\phi_{I}\|_{L^{\infty }(-\infty ,c)}
+\| D^{-1}f\|_{L^{\infty }(c,\infty )}\| |I|D\phi_{I}\|_{L^{1}(c,\infty )})
$$
$$
\leq |I|^{-1}\| D^{-1}f\|_1 C |I|^{-1/2}\Big( 1+\frac{|c(I)-c(I')|}{|I|}\Big)^{1-N}
$$
Now, from $\| D^{-1}f\|_1\leq |I'|\| f\|_1$ we obtain
$$
|\langle f, \phi_{I} \rangle ||I|^{1/2}
\leq C\| f\|_{1}\frac{|I'|}{|I|}\Big( 1+\frac{|c(I)-c(I')|}{|I|}\Big)^{1-N}
$$

For any two integers $k\geq 0$ and $m>0$, there are at most two intervals such that
$|I|/|I'|=2^k$ and the integer part of $1+\frac{|c(I)-c(I')|}{|I|}$ is $m$. Thus we can estimate
$$
\sum_{\tiny \begin{array}{c}I:|I|\geq |I'|\\ I\not \subset 3I'\end{array}}|\langle f, \phi_{I}\rangle | |I|^{1/2}\leq C\| f\|_{1}\sum_{k\geq 0}\sum_{m\geq 1}
2^{-k}m^{1-N}\leq C\| f\|_1
$$
ending the proof of this lemma.

\vskip15pt
Once the weak $L^1$ type inequality is proved, by interpolation we obtain for $1<p\leq 2$
$$
\| \tilde{S}_{k,n}(f)\|_p \leq C_p(2^{-k}n+1)^{\frac{2}{p}-1}\| f\|_p
$$

In order to obtain boundedness for the remaining exponents $2\leq p<\infty $ we consider
the following martingale operator
$$
\tilde{T}_{k,n}(f)(x)
=\sum_{J}\langle f,\psi_{I}\rangle \psi_{J}(x)
$$
where $I$ and $J$ are given by the same relationship that in the definition of $\tilde{S}_{k,n}$. This operator
trivially satisfies $\tilde{T}_{k,n}^{*}(f)=\tilde{T}_{-k,n}(f)$. Actually, the implicit index $j$ does not change and so we may also write
$\tilde{T}_{k,n,j}^{*}(f)=\tilde{T}_{-k,n,j}(f)$.
Moreover, we have that the classical square function of $T_{k,n}(f)$ coincides with $S_{k,n}(f)$:
$$
S(\tilde{T}_{k,n}(f))(x)
=\big( \sum_{J}\frac{\langle f,\psi_{I}\rangle ^2}{|J|}\chi_{J}(x)\Big)^{1/2}
=\big( \sum_{I}\frac{\langle f,\psi_{I}\rangle ^2}{|J|}\chi_{J}(x)\Big)^{1/2}
=\tilde{S}_{k,n}(f)(x)
$$
So, by properties of classical square function and the previous case with $1<p'\leq 2$, we have
$$
\| \tilde{S}_{k,n}\|_{L^p\rightarrow L^p}=\|S(\tilde{T}_{k,n})\|_{L^p\rightarrow L^p}\approx \| \tilde{T}_{k,n}\|_{L^p\rightarrow L^p}
$$
$$
=\|\tilde{T}_{k,n}^{*}\|_{L^{p'}\rightarrow L^{p'}}
=\|\tilde{T}_{-k,n}\|_{L^{p'}\rightarrow L^{p'}}\approx \| \tilde{S}_{-k,n}\|_{L^{p'}\rightarrow L^{p'}}
$$
$$
\leq C_{p'}(2^{k}\log(n)+1)^{(\frac{2}{p'}-1)}
= C_{p}(2^{k}\log(n)+1)^{|\frac{2}{p}-1|}
$$
or
$$
\| \tilde{S}_{k,n}\|_{L^p\rightarrow L^p}
\approx \| \tilde{S}_{-k,n}\|_{L^{p'}\rightarrow L^{p'}}
$$
$$
\leq C_{p'}(2^{k}+\log(n)+1)^{(\frac{2}{p'}-1)}
= C_{p}(2^{k}+\log(n)+1)^{|\frac{2}{p}-1|}
$$

This ends the proof.

\end{document}